\documentclass[11pt]{article}
\usepackage{amsmath,amsfonts}
\usepackage{graphicx}
\usepackage{subfigure}
\usepackage{enumerate}
\usepackage{natbib}
\usepackage{url} 

\usepackage[colorlinks=true,linkcolor=black,citecolor=black,urlcolor=black]{hyperref}

\usepackage{bbm}
\usepackage{amsfonts}
\usepackage{algpseudocode}
\usepackage{algorithm}
\usepackage{algorithmicx}
\usepackage{caption}
\usepackage{subfigure}
\usepackage{soul}
\usepackage{tikz}
\usepackage[final]{pdfpages}
\usepackage{tabularx,ragged2e,booktabs,caption}

\newcolumntype{C}[1]{>{\Centering}m{#1}}

\newcommand{\blind}{0}

\newcommand{\ph}{\varphi}

\def\E{\mathbb{E}}
\def\mR{\mathbb{R}}

\def\tX{{\tilde X}}
\def\tx{{\tilde x}}

\def\definedas{\stackrel{\Delta}{=}}

\newcommand{\cN}{{\mathcal N}}

\newcommand{\cL}{{\mathcal L}}

\newcommand{\cZ}{{\mathcal Z}}

\def\real{\mR}
\def\bc{\begin{center}}
\def\ec{\end{center}}
\def\L{\mathcal L}

\def\smallhalf{\mbox{ $\frac{1}{2}$}}

\newcommand{\beq}{\begin{eqnarray}}
\newcommand{\eeq}{\end{eqnarray}}
\newcommand{\beqq}{\begin{eqnarray*}}
\newcommand{\eeqq}{\end{eqnarray*}}

\begin{document}

\def\spacingset#1{\renewcommand{\baselinestretch}%
{#1}\small\normalsize} \spacingset{1}


\if0\blind
{
  \title{\bf Modelling Persistence Diagrams with Planar Point Processes, and Revealing Topology with  Bagplots}
  \author{Robert J.\ Adler\thanks{
    Research supported in part by \textit{URSAT: Understanding Random Systems via Algebraic Topology}, ERC  Advanced Grant 320422.} \ and 
    Sarit Agami$^*$\hspace{.2cm}\\
    Andrew and Erna Viterbi Faculty of  Electrical Engineering,\\
     Technion -- Israel Institute of Technology\\
}
  \maketitle
} \fi

\if1\blind
{
  \bigskip
  \bigskip
  \bigskip
  \begin{center}
    {\LARGE\bf Title}
\end{center}
  \medskip
} \fi

\bigskip
\begin{abstract}

We introduce a new model for planar point point processes, with the aim of capturing the structure of  point interaction and spread  in persistence diagrams. Persistence diagrams themselves are a key tool of TDA (topological data analysis), crucial for the delineation and estimation of global topological structure in large data sets. To a large extent, the statistical analysis of persistence diagrams has been 
hindered by difficulties in providing replications, a problem that was addressed in an earlier paper, which introduced a procedure called
RST (replicating statistical topology). Here we significantly improve  on the power of RST via the introduction of a more realistic class of models for the persistence diagrams. In addition, we introduce to TDA the idea of bagplotting, a powerful technique from non-parametric statistics well adapted for differentiating between topologically significant points, and noise, in persistence diagrams.

Outside the setting of TDA, our model provides a setting for fashioning point processes, in any dimension, in which both local interactions between the points, along with global restraints on the general point cloud, are important and perhaps  competing.

%
\end{abstract}

\noindent%
\vfill
\newpage
\tableofcontents

\section{Introduction}
\label{sec:intro}

Over the past decade there has been considerable interest and success in the exploitation of topological thinking, particularly that coming from Algebraic Topology, in developing tools  for the analysis of large and complex data sets and networks. Under the brand name of 
`Topological Data Analysis' -- hereafter `TDA', this approach has  been put on a reasonably solid mathematical footing, and 
applications have been as widespread as signal processing, \cite{Barbarossa}, genetic analysis for some breast cancers, \cite{CarlssonCancer}, functional (neural) networks in the brain, \cite{Bendich,Vaccarino}, and cosmology, \cite{PNAS,Sousbie1,Sousbie2,PratyushOld}. In addition, there are  connections between TDA, dimension reduction,  and machine learning that are both mathematically deep and highly practical, as can be seen by the adoption of both of these techniques by the start up {\it Ayasdi}. 

In almost all of the applications of TDA, and in most of the theory, persistence diagrams (or, equivalently, barcodes) arise  as the key topological summary  of the underlying phenomenon, or  mathematical models,   being studied, and provide the basis for all subsequent analysis. We shall describe persistence diagrams briefly, as well as giving some pertinent examples,  in Section \ref{sec:PD-new}.
With relatively few exceptions, notably \cite{chazal,fasy,RobTurner,OmerSayan,bubenik,wasserman} (see  additional citations in the SI Appendix (Sec.\  1.3) to \cite{PNAS})  TDA has not employed statistical methodology as part of its approach, and, as a consequence, has typically been unable to associate clearly defined levels of statistical significance to its discoveries. While there may be a variety of reasons for this, one of the main obstacles to doing so is that  the mathematical challenges involved in computing the statistical distributions of topological quantifiers have so far proven to be intractable. 

There have been a number of attempts to solve this problem by reducing the entire persistence diagram to a lower dimensional summary statistic, usually combining this with some sort of bootstrapping procedure (e.g.\ \cite{chazal,fasy,RobTurner}). In 
 \cite{PNAS} we introduced a new approach, dubbed  `Replicating Statistical Topology', hereafter `RST'. We shall describe all of this below in Section \ref{sec:inference-new}, along with examples in Section  \ref{sec:examples-new}, but the main point, from the point of view of the current paper, was that RST was based on a semi-parametric procedure  for  providing multiple simulations  of an observed persistence diagram. The first step involved fitting a parametric model to a given persistence diagram via a Gibbs distribution, and the current paper is mainly about a developing a significant  improvement to the this step. 
  
 After developing the new  model, and seeing that, from many aspects, it, together with an appropriate MCMC simulation procedure,  does an excellent job of producing (almost) independent persistent diagrams from a single initial example, we introduce a powerful new method for distinguishing between topologically informative points and `topological noise'. This method is based on the statistical procedure known as bagplotting. To the best of our knowledge this has never before been used in the TDA setting (with the possible exception of \cite{Nature}, but in a different role). Nevertheless, as we shall show in Section 
 \ref{sec:bagplots}, the use of bagplots to distinguish between `signal' and `noise' in a persistence diagram is a natural
 application of this technique. On the other hand, without the replication methods we develop in the paper, it would not allow assigning levels of statistical certainty to these discoveries. This, of course, is the role of RST.
 
 In fact, although our motivation comes from TDA and RST, what this paper provides extends beyond these specific settings. Persistence diagrams are (at each homology level) random, planar, point processes about which very little is known. What we do know, as will be explained below by investigating specific cases, is that they are most definitely {\it not} Poisson processes, as the points interact strongly at a local level.  This interaction strongly suggests the adoption of a Gibbsian approach to modelling, along with its associated 
  MCMC technology for later simulation. On the other hand, there is a very definite global structure to the `cloud' of points making up a persistence diagram that is hard to capture  via purely local interactions. Placing these two (competing) requirements into a single, data based, model is the truly novel contribution of this paper.  The application of the new model in  an RST setting, in the examples of Section \ref{sec:examples-new}, then basically follows the recipe of  \cite{PNAS} and \cite{agamiadler}. However, the fact the resultant analysis   significantly improves on the previous ones  is compelling  justification for  applying the new model in the RST setting in particular, and TDA in general.
  
  \vskip0.2truein
\noindent {\bf Acknowledgment:} We are grateful to Katherine Turner who, at a conference in Japan, suggested that we should be able to improve on the model of \cite{PNAS} by incorporating  information  on the   global shape of  the persistence diagrams into the model. It was an insightful and useful suggestion.


\section{Persistence diagrams}
\label{sec:PD-new}

As mentioned in the Introduction, persistence diagrams are undeniably the most single useful tool of TDA. To make this paper self-contained, they should now be carefully defined, but we have decided not to do so. There are three main reasons for this:
\begin{itemize}
\item[(i)] A reader who is interested in this paper from the point of view of its application in TDA will already be so familiar with the notion of persistence that any description which we would provide would simply be skipped over.

\item[(ii)] For the reader interested in our results from the point of view of modelling point processes, the intrinsic interest of this is sufficient that  understanding its connection to persistence diagrams is not necessary.

\item[(iii)] There are already so many good introductions to TDA and persistence  that we would  be hard put to write anything without accusations of plagiarism. For statisticians, there is the comprehensive and up to date  review  \cite{wasserman} which includes close to  100 references.   (See also \cite{RobTurner} which also considers a number of  interesting hypothesis testing issues for persistence diagrams.)  
For the reader wanting a more general approach, there are the 
recent excellent and quite different books and reviews 
\cite{Carlsson-review,CarlssonReview,EdelsShortCourse,EdelsHarerSurvey,EdelsHarerBook,Afra,Oudot} and \cite{ghrist2014elementary}, all of which 
give broad  expositions of  the homology needed for TDA.  (For a description of the history of persistence, see the Introduction in \cite{EdelsHarerSurvey}.)
\end{itemize}

Despite having abrogated our didactic responsibilities to the reader regarding persistence diagrams, we nevertheless do need to set up
some notation and conventions.

Throughout this paper, all the specific  examples we shall treat arise as  persistence diagrams for which the underlying filtration is generated by the excursion, or upper level, sets of a real valued, smooth function $f$ on a nice space $\cZ$. Thus, the filtration is given by 
$\{\cZ_u, u\in\real\}$, where 
\beq
\label{Au-defn}
\cZ_u \ \definedas \{z\in \cZ : f(z) \in [u,\infty)\} \
\equiv\ f^{-1}([u,\infty)).
\eeq
Given this filtration, there will be a persistence diagram for each persistence homology $H_k(\cZ)$ of order $k$, $0\leq k \leq \dim (\cZ)$. Each point 
$(b,d)$ of each diagram represents a generator in   $H_k(\cZ)$ which was `born' at level $b$ and `died' at level $d$.
Since $u<v$ implies $\cZ_u\supseteq \cZ_v$, we always have $b>d$, and since we shall place births on the vertical axis and deaths on the horizontal axis of our persistence diagrams, the points will always lie in the half plane above the line $b=d$.

We now turn to two of the motivating examples behind our work. While each comes from quite a different background, they share a common theme in that they involve excursion set filtrations.

\subsection{Gaussian random field excursions}
\label{subsec:Gauss-new}

By `Gaussian random field' we mean a random, real valued, function $f$ over a topological space $\cZ$, the finite dimensional distributions of which are multivariate Gaussian. We shall always assume that $f$ has smooth sample paths. In particular, implicitly assuming that $\cZ$ is  a $C^2$ (perhaps stratified)  manifold, we shall assume throughout that $f$ is, with probability one, twice differentiable. For a full theory of such random fields, including many topological properties of their excursion sets which are quoted below without reference, see any of \cite{RFG,ATSF,HRF}. Gaussian random fields, including the topological properties of their sample paths, have been at the core of numerous applications, perhaps the most notable ones being in cosmological and neurophysiological problems. (The history of topology mixed with `modern' random field theory in cosmology probably starts with \cite{Bardeen,Coles1988},  with recent papers including TDA concepts such as persistence including 
\cite{feldbrugge2012,PratyushOld,cole2018,PNAS} and \cite{Nature}. The last two of these, in particular, treat modelling of the CMB (Cosmic Microwave Background radiation) and are motivating applications for the example of Gaussian excursion sets that we are about to discuss. As for the applications of Gaussian processes to neurophysiology, \cite{friston1994statistical}, along with its many thousands of citations, is a good place to enter the literature.)

Because of the importance of Gaussian random fields in stochastic modelling, and the growing interest in them from a topological viewpoint, we shall discuss this example in some detail.

\begin{figure}[h!]
\bc
      \includegraphics[width=6.0in, height=3.4in]{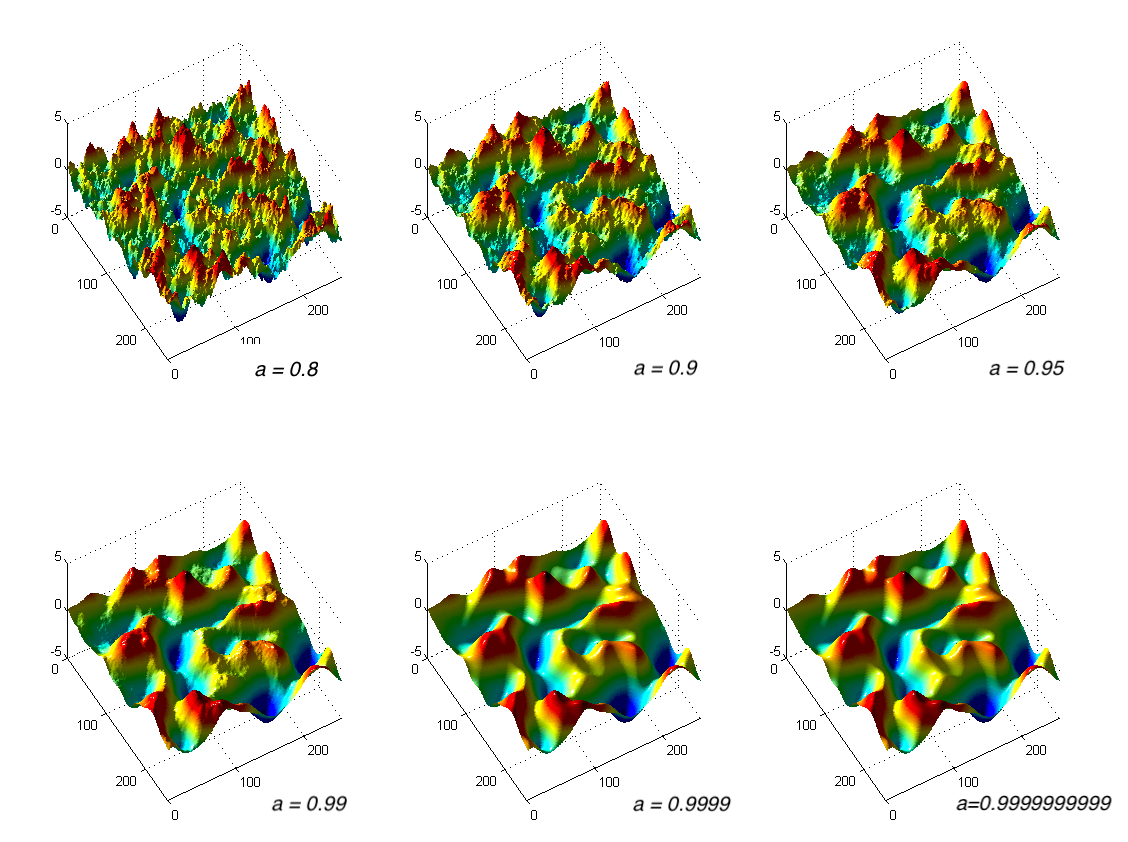}
\caption{\footnotesize
Simulations of a zero mean, two-dimensional, stationary, Gaussian  random field with covariance function \eqref{eq:covariancefunction}. (With thanks to Eliran Subag.) Note that lower values of $a$ lead to higher levels of local variability, while the global structure of the field remains quite stable.}
\label{fig:variantalpha}
\ec
\end{figure}

For a specific class of examples, we take $f$ to be the stationary, isotropic, zero mean, Gaussian random field on $\real^2$ with covariance function 
\beq
\label{eq:covariancefunction}
R(x) \ \equiv \ \E\{f( y)f(y+x)\} \ = \  e ^{-b \|x\|^{2a}},
\eeq
with $b>0$ and $a\in (0,1]$. The parameter $b$ is essentially a spatial scalling parameter, while $a$ controls both long range correlations and the local smoothness of $f$. In particular, if $a=1$, then the random field has $C^\infty$ sample functions (with probability one, a qualifier that we shall drop from now on) while, for $a<1$, $f$ is non-differentiable but continuous, and satisfies a H\"older condition of order  $\alpha$, for all $\alpha\leq a$. Examples of realisations of these processes, over the unit square $[0,1]^2$, for $b=100$ and various $a$ are given in Figure \ref{fig:variantalpha}. (The numbers on the base axes reflect the fact that the simulations are taken over a $256\times 256$ subgrid of $[0,1]^2$.)

\begin{figure}[h!]
    \centering
        \includegraphics[scale=.11]{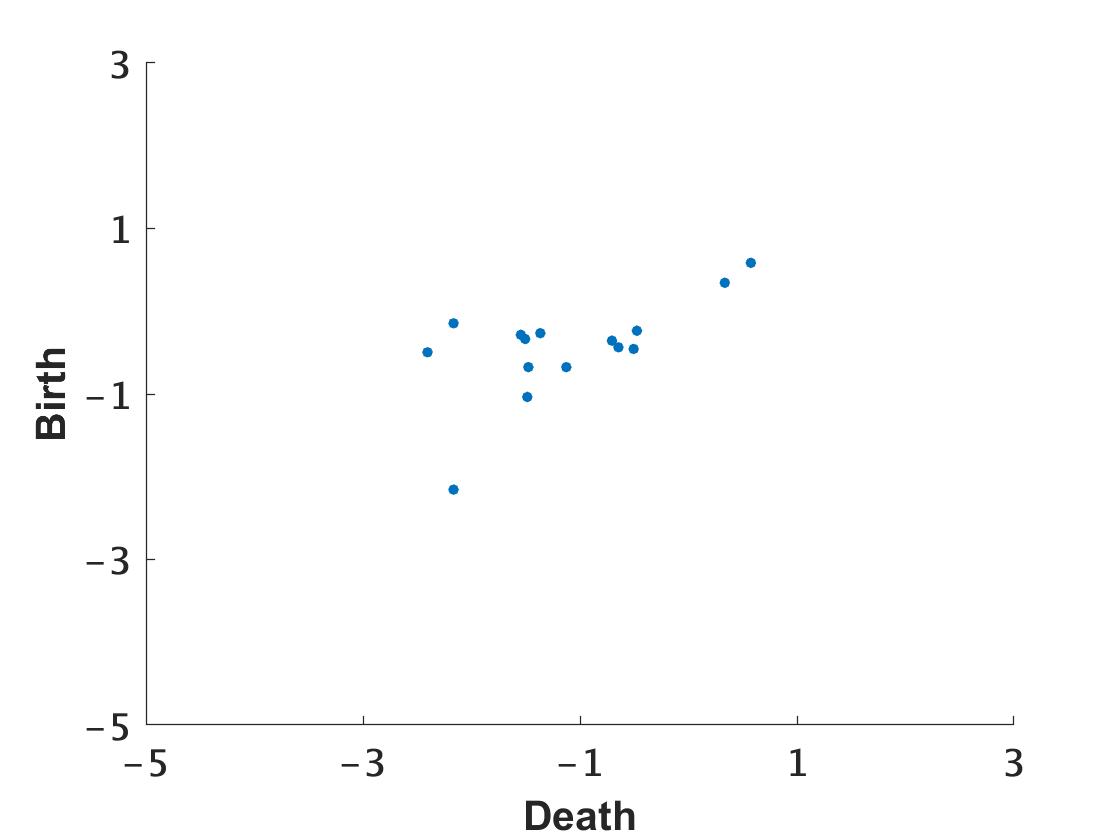}
    \quad
        \includegraphics[scale=.11]{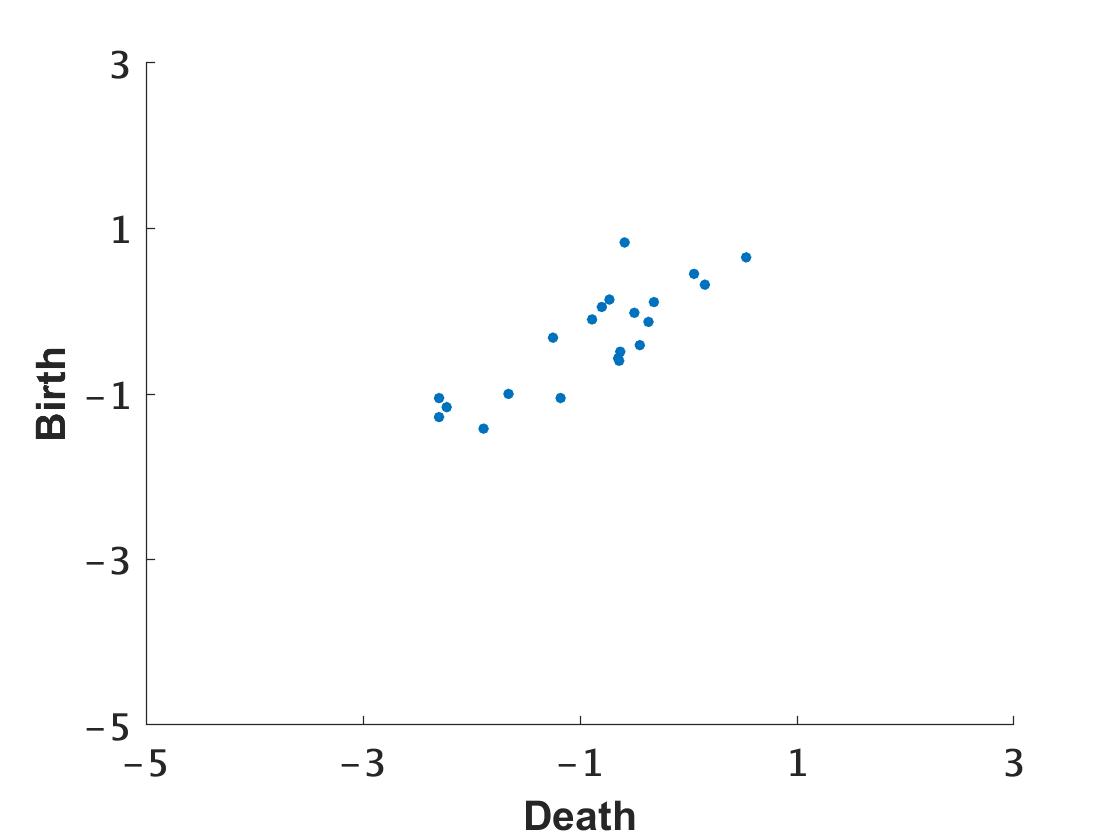}
    \quad
        \includegraphics[scale=.11]{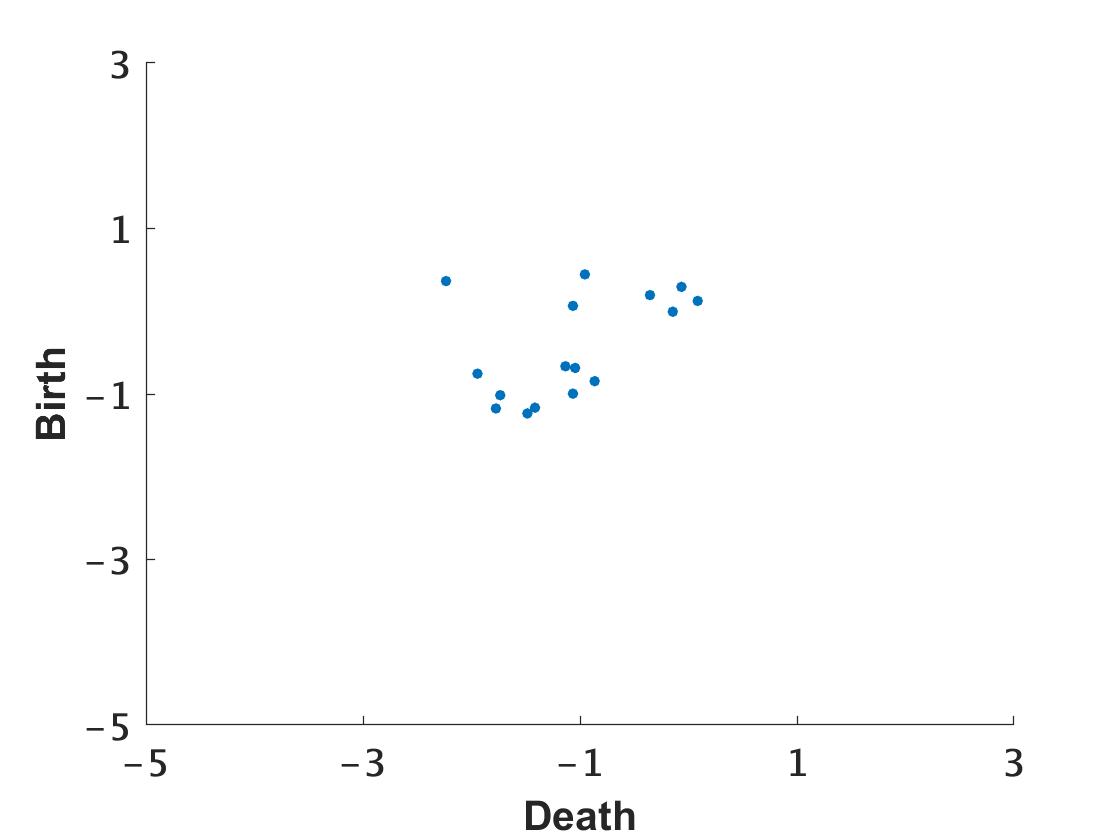}
  \caption{\footnotesize
      $H_1$ persistence diagrams for the excursion sets of three simulations of the Gaussian process of Figure \ref{fig:variantalpha}, with $a=1$.  }
\label{fig:r1}
\end{figure}

Turning now to persistence, 
Figure \ref{fig:r1} gives three examples of  the $H_1$ persistence diagrams of the upper level set filtration  for three simulations of the random field in the bottom right of Figure \ref{fig:variantalpha}. (We shall typically work with $H_1$ diagrams for this example rather than the $H_0$ diagrams; viz.\ we shall work with the topology of holes, rather than connected components, in the excursion sets. However, everything that we shall have to say for this case has an immediate and essentially equivalent parallel in the $H_0$ case.) 
It is hard to say much from these Figures, other than perhaps to surmise that there seems to be a repulsive effects between points in each diagrams, a conjecture that one might put down to either `experience' in looking at examples of point processes, or simple subjectivity.   Certainly there is not much that can be said about the global shape of the point cloud. 

\begin{figure}[h!]
 \centering
       \subfigure[]
    {
      \includegraphics[width=2.9in, height=2.1in]{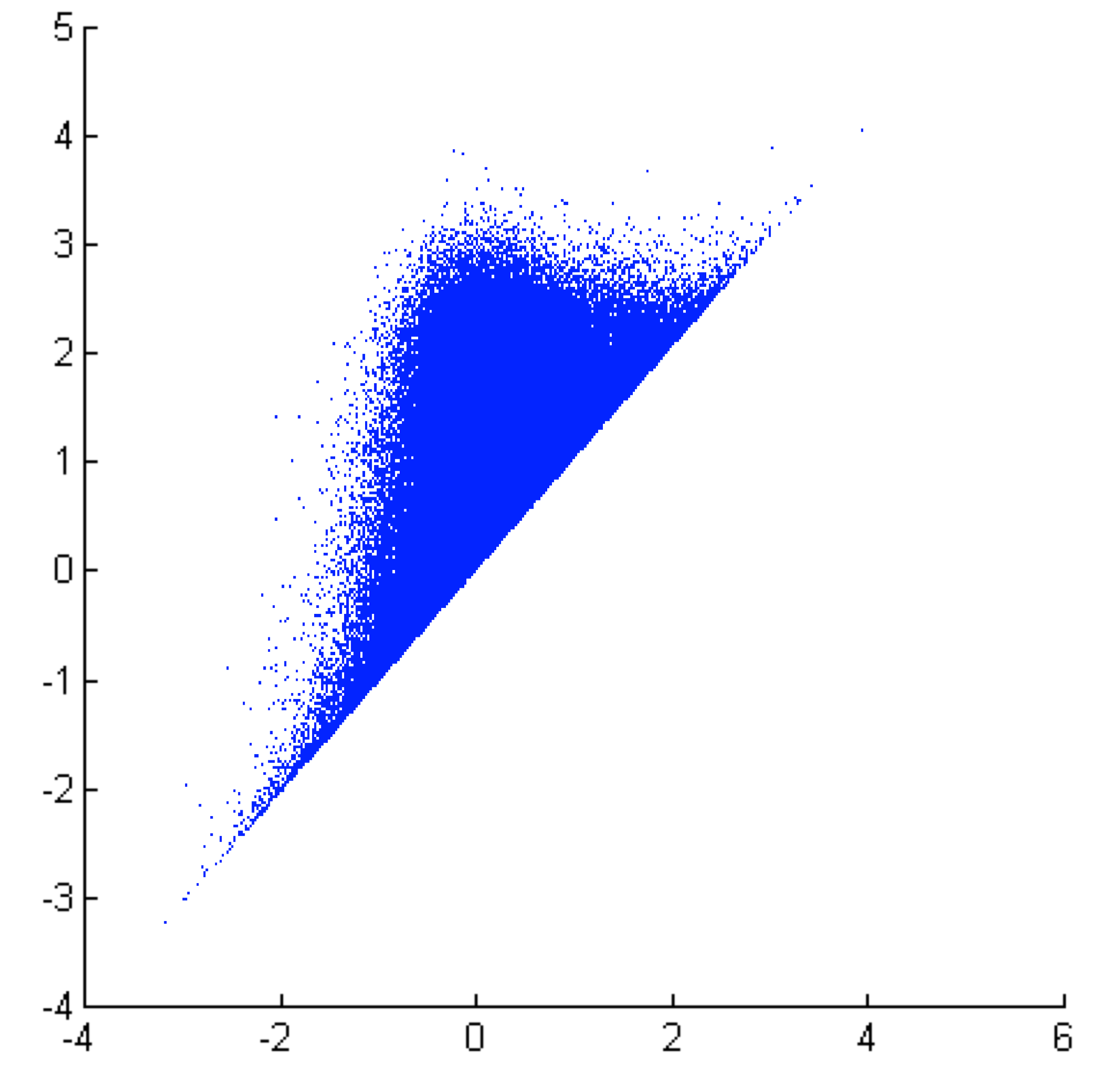}
    }
 \qquad      \subfigure[]
    {
        \includegraphics[width=2.4in, height=1.7in]{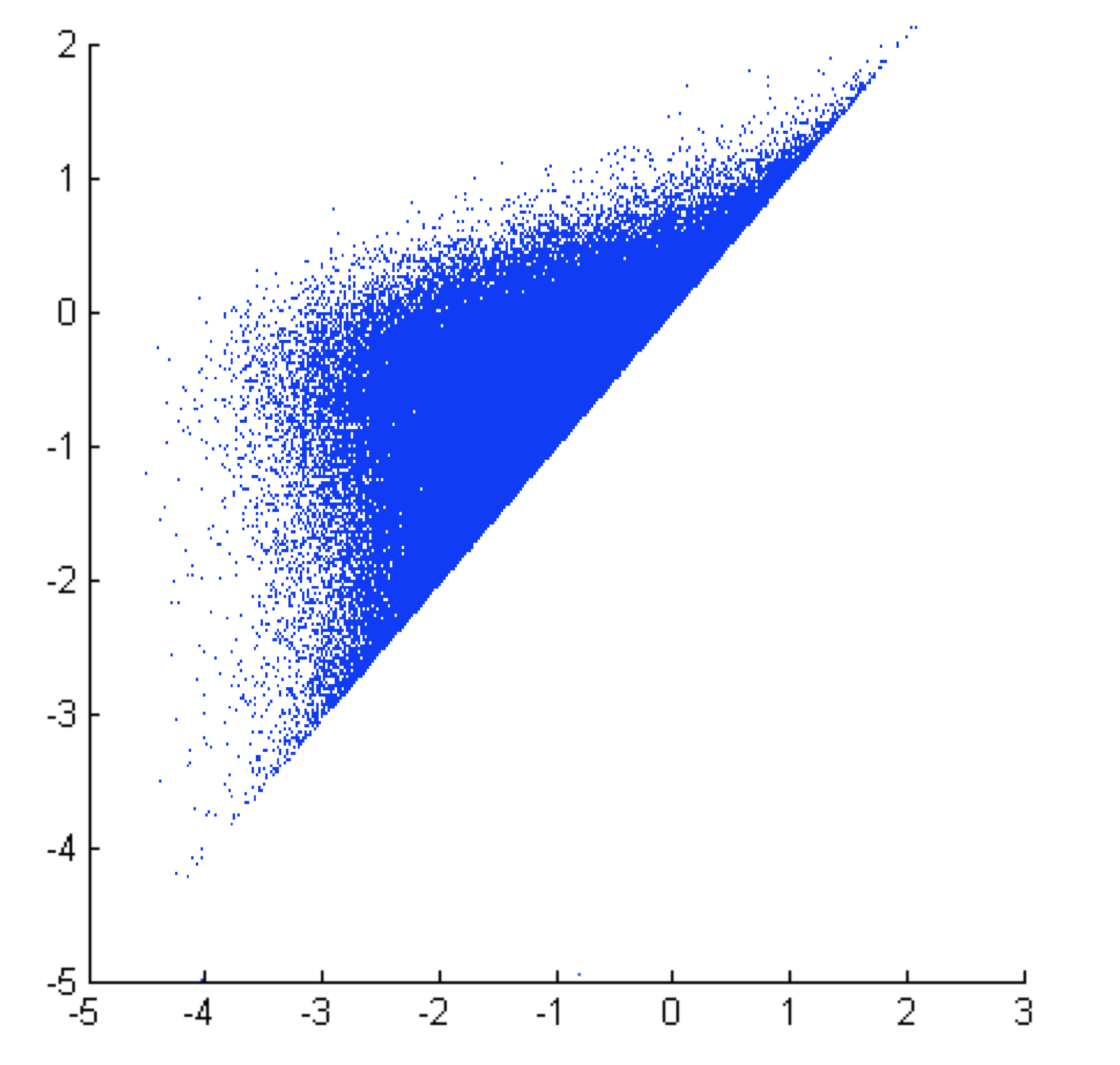}
    }
\caption{\footnotesize    Superimposed persistence diagrams ($H_0$ in (a) and $H_1$ in (b)) for excursion sets of 10,000  Gaussian random simulations. The $H_0$ diagram contains a total of 251,330 points, while there are only 174,006 points in the $H_1$ diagram. That is, there are of the order of 20 points in individual diagrams, as  in Figure \ref{fig:r1}.}
\label{fig:GaussH}  
\end{figure}

On the other hand, superimposing the points from 10,000 such simulations into one picture for both $H_0$ and 
$H_1$ diagrams  yields Figure \ref{fig:GaussH}, from which a global structure of the point clouds is very clear. (As an aside,  note that in the $H_0$ diagram, as in all future examples, the `points at infinity'  in the individual persistence diagrams are not shown.)

Combining the pictorial information in Figures \ref{fig:r1} and \ref{fig:GaussH} it is clear that in looking for a model for random persistence diagrams that will cover this case, we need to find something that allows close points in the diagrams to locally repel one another, while at the same time complying with global constraints.

Before looking for such a model, however, we  note that the `natural' path would be to rigorously derive analytic relationships between  Gaussian random fields -- `parameterised' by their covariance functions  -- and the persistence diagrams of their excursion sets. However, this is beyond the capabilities of current Gaussian theory, and may well be an almost impossible problem. For example, although much is known about the mean Euler characteristics of these excursion sets, as well as some qualitative results about their shape at asymptotically high and low levels, all attempts at computing exact expressions for something as topologically simple as their mean Betti numbers have, so far, met with 
little to no success.

What is known comes from the Morse theory connection between topology and critical points. At least in some cases it is possible to compute precise expressions for the mean number of critical points of Gaussian fields (cf.\ \cite{HRF}, Chapter 6, for a survey) and to use these and the Morse inequalities to say something about mean Betti numbers (cf.\ \cite{feldbrugge2012} for a nice example of this in two dimensions). Another  
observation with rigorous backing comes from the so-called `Slepian models' for Gaussian processes, that show that Gaussian critical points tend not to be close to one another. (This is also well documented in the  connections between the complex zeroes of Gaussian processes and determinantal, planar point processes; e.g.\ \cite{Hough}.) Taking into account that each  of both the birth and death coordinates of each point in a persistence diagram  is also a critical level of $f$, this literature validates the impression gained from Figure \ref{fig:r1} that there is a natural repulsion between the points.
(Exceptions to this will occur at the extreme critical points (i.e.\ high or low) where the classical extremal, Poisson, limit theory of Gaussian processes leads to the Poisson-like thin spread of points at the top of  (a) of Figure \ref{fig:GaussH}   and at the left of (b).)

These models also explain that asymmetries present in both diagrams of Figure \ref{fig:GaussH}, as well as the general shapes of the clouds.

In summary, therefore, there is significant motivation, from this example alone, to develop a model for random persistence diagrams that include both local repulsion and global constraints.
 
\subsection{Sampling from non-concentric circles with different radii}
\label{subsec:2circles-new}

As our second example we take what is basically a toy model, something more familiar in the TDA literature.  It is much simpler, at least in that it requires no previous, specific, probabilistic knowledge in order to appreciate its subtleties, but it still exhibits nuances requiring the same demands of a stochastic model that arose with the random field example.

\begin{figure}[h!]
 \centering
    \subfigure[]
    {
        \includegraphics[width=1.8in, height=1.8in]{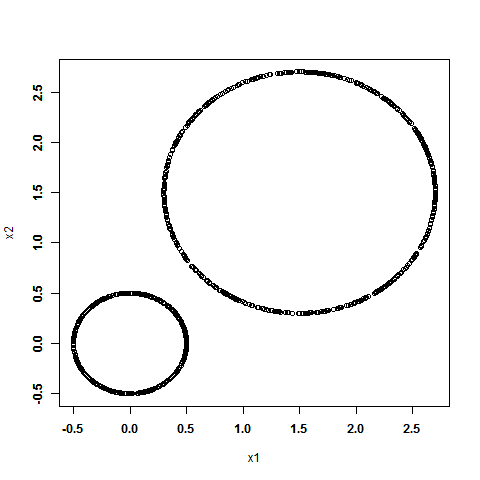}  
          }
      \subfigure[]
       {
      \includegraphics [width=1.8in, height=1.8in]{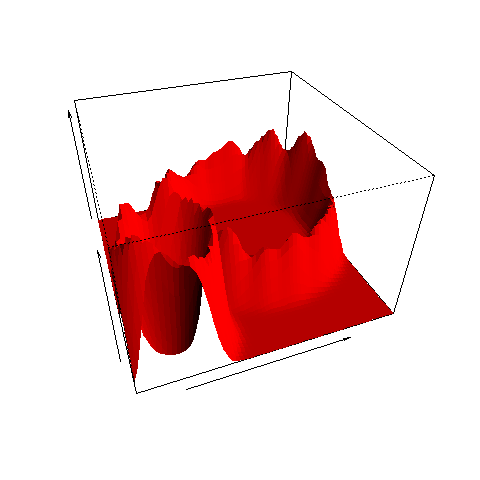} 
       }
       \subfigure[]
       {
\includegraphics[width=1.8in, height=1.8in]{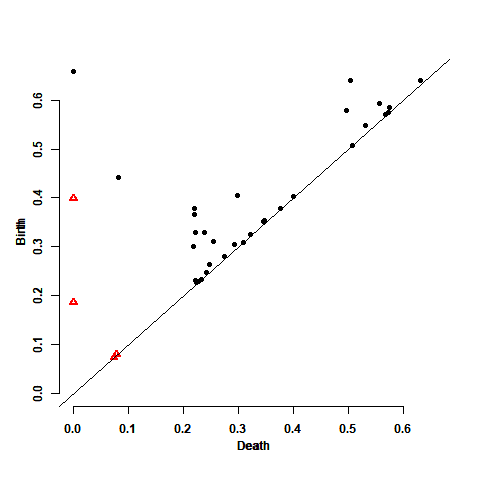}
}
\caption{\footnotesize
 A random sample (a) from two non-concentric circles, 500 points from the  smaller circle and 650 points from the larger circle, with a kernel density estimate (b) and the persistence diagram (c) for its upper level sets. Black dots represent  $H_0$ persistence points, and  red triangles   $H_1$ points.}
\label{fig:farcircle}
\end{figure}

In particular, we sample from the two circles of Figure \ref{fig:farcircle}(a), taking 500 points sampled uniformly at random from  the smaller circle, with  a radius of 0.5, and $650$ points  from the  larger circle, of radius of 1.2.  We then compute a smoothed density function estimate for the sample, using a Gaussian kernel with  bandwidth  $0.1$. The result is shown in (b) of the same figure, and in (c) we  show the persistence diagram of the upper level set filtration of this function. It contains $N=32$ points of $H_0$, these being the the black circles indicating $H_0$ (zero-th homology) persistence. The red triangles correspond to $H_1$. We see two black circles (the two connected components; for this example the `point at infinity' is included, as the point with death coordinate zero)   and two red triangles (holes) somewhat isolated from the other points in the diagram, as is to be expected. 

As in the previous example, it would be nice to be able to say something about the stochastic makeup of both the signal and noise in the persistence diagrams of this example, given the simple sampling that lies behind it. In fact, in this case more can be said, and the reviews \cite{KahleOmerReview} and \cite{KahleS}, along with more recent papers such as \cite{BKS,OmerSayan}, contain many results about the asymptotics of this kind of scenario, but only when the number of points being sampled diverges to infinity.

\begin{figure}[h!]
 \centering
      \includegraphics[width=3.6in, height=2.55in]{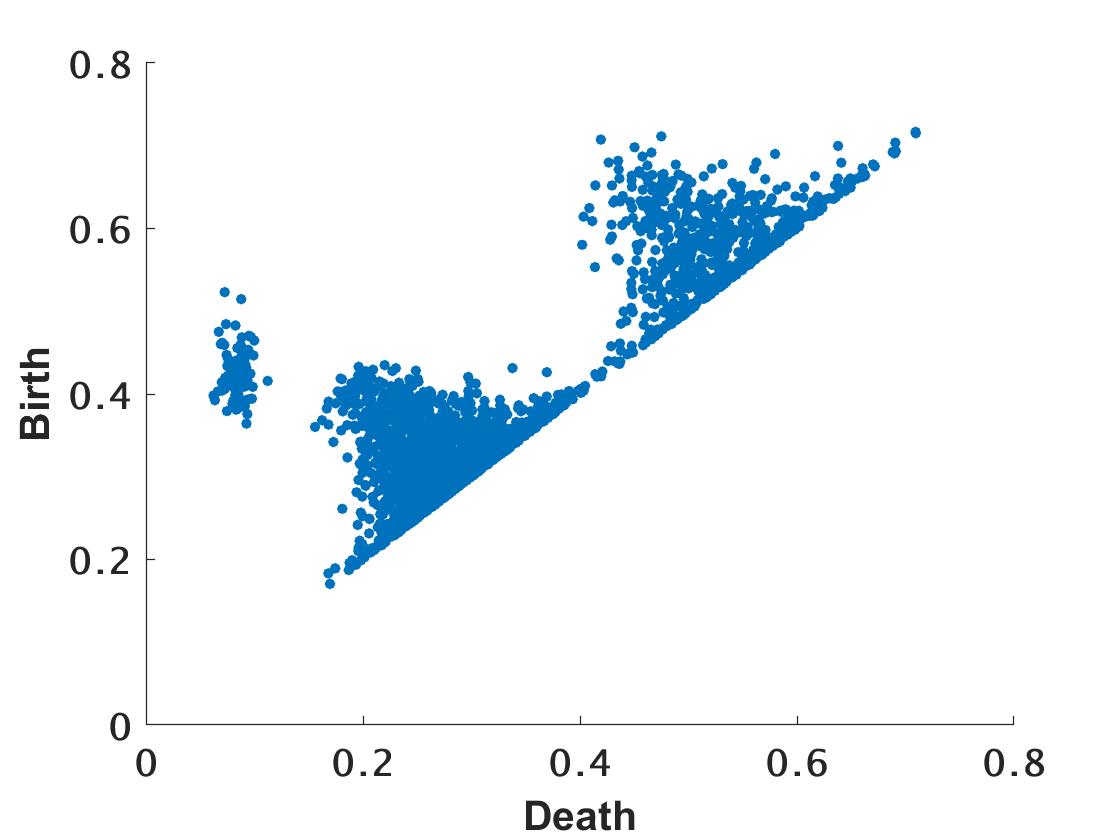}    
\caption{\footnotesize Superimposed persistence diagrams for 100  simulations of  the example of Figure \ref{fig:farcircle}.}
\label{fig:TwoCirclesPH:rob}
\end{figure}

Adopting the same approach as in the Gaussian random field example, Figure \ref{fig:TwoCirclesPH:rob}
shows the points of 100 independent realisations of the $H_0$ persistence diagram of  Figure \ref{fig:farcircle}(c) (but without the points at infinity).  A pattern arises: While the local repulsion between points which can be seen in the individual case is now blurred, the appearance of three clearly delineated collections of points is now clear. The leftmost comes from the connected component which is the smaller circle. (The larger circle generated the missing point at infinity.) The remaining two groups represent the topology of noise; viz.\ the sampling error which generated the irregularities in the heights of the two cones in (b). 

Thus, once again, we are in a situation in which random persistence diagrams show both local repulsion between points, but restrained  in a quiet remarkable fashion by  global constraints.

%
%


\section{Modelling point processes and persistence diagrams}
\label{sec:modelling-new}

\subsection{An infeasible Gibbs distribution}
The model that we plan to adopt for persistence diagrams will, as in \cite{PNAS} be based on a Gibbs distribution. That is, given a finite collection $\tX_N=\{X_1,\dots,X_N\}$ of continuous random variables, each taking values in  $\real^D$, for some $D\geq 1$, and with joint probability density $\ph_\Theta(\tx_N)$,
dependent on a multi-dimensional parameter $\Theta$, we assume  that
$\ph_\Theta$ can be naturally written in the form
\beq
\ph_\Theta(\tx_N) = \frac{1}{Z_\Theta} \exp ( - H_\Theta (\tx_N)),
\label{eq:Gibbs}
\eeq
where the `Hamiltonian' $H_\Theta: \real^N\to\real$
describes the `energy' of configurations $\tx_N$. 
The normalisation, or partition function, $Z_\Theta$,  a function of $\Theta$, is infamously hard to evaluate, and so we shall soon move from this model to a slightly different one. 

In view of the motivating examples of the previous section, we start with a Hamiltonian component that controls global shape, by assuming that there is an underlying density $f^G_{\Theta_G}(x)$ which describes the univariate (common, marginal) density of the $X_i$. Under $f^G_{\Theta_G}$, if the $X_i$ were independent random variables, then we could write their joint density \eqref{eq:Gibbs} as 
\beq
\notag
\ph^G_{\alpha,\Theta_G}(\tx_N) \ =\     \frac{1}{Z^G_{\alpha,\Theta_G}}  \prod_{i=1}^N \left(f^G_{\Theta_G}(x_i)\right)^\alpha  &=& \ 
\frac{1}{Z^G_{\alpha,\Theta_G}} \exp\left(\alpha \sum_{i=1}^N \ln \left(f^G_{\Theta_G}(x_i)\right)\right)\\  &\definedas& 
\frac{1}{Z^G_{\alpha,\Theta_G}} \exp\left(\alpha  H^G_{\Theta_G}(\tx_N)\right),
\label{eq:Gibbsalpha}
\eeq
with $\alpha=Z^G_{\Theta_G}=1$. Taking values of $\alpha$ other than $\alpha=1$ widens the family of distributions, while retaining independence, and the need for   additional freedom of this form will become clear in a moment. Note, however, that, whatever the value $\alpha$, the $X_i$ are always independent and identically  distributed under \eqref{eq:Gibbsalpha}. Furthermore, their spatial positioning is determined by  the support of $f^G_{\Theta_G}$.

We now turn the local interactions for which we introduce the joint densities
\beq
\ph^L_{\Theta_L}(\tx_N) \ =\     \frac{1}{Z^L_{\Theta_L}}  \exp\left(-H^L_{\Theta_L} (\tx_N)\right) \
\definedas \  \frac{1}{Z^L_{\Theta_L}} 
 \exp\left(-\sum_{i=1}^N   \L^K_{\Theta_L} \left(x_i\,\big| \, \cN(x_i) \right)\right)
  \label{eq:GibbsH2}
\eeq
where, for a $K\geq 1$,  and $x,y\in\real^D$,  
\beq
  \L^K_{\Theta_L} \left(y\big| \, \cN(x) \right) \ \definedas \ \sum_{k=1}^K \theta_L(k) \sum_{z\in \cN_k(x)}  \|z-y\| ,
  \label{eq:LKL}
  \eeq
  and $\cN_k(x)$ is the set of the $k$ nearest neighbours to $x$.
  
  Putting all this together, and combining the parameters $\alpha$, $\Theta_G$ and $\Theta_L$  into a single vector $\Theta$, and  we can define a Gibbs distribution with Hamiltonian
  \beq
 H_\Theta (\tx_N)) \ = H^L_{\Theta_L}(\tx_N) - \alpha  H^G_{\Theta_G}(\tx_N),  
    \label{eq:TheGibbs}
  \eeq
  and corresponding partition function.
  
 The probabilistic import of the combined Hamiltonian should be clear, as well as the reason for introducing the parameter $\alpha$.
 $ H^L_{\Theta_L}$ controls local interactions. Acting alone  -- viz.\ without the presence of $H^G_{\Theta_G}$ --  negative values of the $\theta_L(k)$ encourage repulsion among points, and positive values encourage attraction. On the other hand,  taking $\alpha>0$,  $H^G_{\Theta_G}$ acts to keep the points distributed within the support of $f^G_{\Theta_G}$. The larger $\alpha$ is, the stronger is this global effect.
 
 Unfortunately, although the above model fills the requirements outlined in the previous section, it has the drawback of all such models in that the difficulties in computing the `corresponding partition function', even numerically, are formidable. In addition, it requires that we have, {\it a priori}, a specific choice   for the family of densities   $f^G_{\Theta_G}$. To overcome these issues, for both model estimation and subsequent simulation, we take the two-pronged approach of replacing $f^G$ by a non-parametric estimate of it, based on the data $\tx_N$, and replacing the true Gibbs model by a pseudo-likelihood version of it.

\subsection{A feasible, semi-parametric,  pseudo-likelihood approach}
\subsubsection{Choosing the shape prior}
\label{subsubsec:choosingprior}
The first step in applying the model of the previous subsection is to decide on a `global shape' density function $f^G$ (where we now drop the explicit dependence on the parameter set $\Theta_G$).

A number of paths are possible. For example, one could adopt, 
 {\it a priori},  a family of candidates, to then be estimated, typically parametrically, from the data. While attractive, this adds computational complexity to a problem that is already (as we shall see below) computationally intensive, and so impractical from a  computational viewpoint. More importantly, given how little is known,  theoretically, about the distributions of persistence diagrams, it is hard to know what would form a good class of candidates. 

If one is working with a specific problem in which the point set arises from a family of experiments that have been carried out often in the past, a number of options are available. For example, if one knew that we were dealing with a persistence diagram that arose from the excursion set filtration of a Gaussian process, without knowing anything additional about the process, then the fact that persistence diagrams in this setting typically have a shape akin to those of Figure \ref{fig:GaussH} might lead one to choose a parametrised class of densities that exhibit such behaviour, and then estimate the parameters for the particular problem. This could be done prior to the estimation of the full model that we are about to develop, or as part of a likelihood maximisation for the full model. 
Many additional methods, all based on some sort of averaging of point processes and density estimation, are also available. 

In the setting of persistence diagrams, one could, for example, work with average  persistence landscapes, as in \cite{Bubenik15}, or  the Fr\'echet mean of a collection of diagrams, as in 
 \cite{TurnerFrechet1}. Any of these can exploited to obtain an estimate of $f^G$.
 
 However, for applications to TDA, we wish to remain in the setting of \cite{PNAS}, in which no previous information is available, and in which only one persistence diagram is available. As explained there, this is more often than not the situation in applications, particularly if they are `Big Data' or Cosmology applications, in which there is only one data set.
 
 Given this, there is little we can do beyond estimating $f^G$ by a simple, non-parametric, density estimator, and so, for the point set $\tx_N$, we choose the kernel density estimate (KDE)
 \beq
\label{eq:kernel}
\hat f^G(x) \ = \  \frac{1}{N(\text{det}\Sigma)^{1/2}  (2\pi)^{D/2}} \sum_{i=1}^{N}
 \exp\left(-\smallhalf\Sigma^{-1/2}(x-x_i)(x-x_i)'\Sigma^{-1/2}\right),\qquad x\in\mathbb R^D,
\eeq
 where our vectors are column vectors and $\Sigma$ is a symmetric, positive definite, scaling matrix. In many cases one can take $\Sigma$ to $\eta^2$ times the identity matrix, where $\eta >0$ is a bandwidth parameter (e.g.\ \cite{Wand}).
Standard asymptotic  theory indicates that, in this case, $\eta$, or, in the general case $\max_k \Sigma(k,k)$, should be chosen to be of order $O(N^{-1/(4+D)})$ (e.g.\ \cite{SilvermanDensity,WassermanAll}). However, in most of the examples that we shall treat later $N$ is in the range 15--25, and so we shall need something other than asymptotic theory to chose $\eta$ or $\Sigma$.
 
 Note that even if the points in $\tx_N$ lie in a (perhaps bounded, or even lower dimensional) subset of $\real^D$, $\hat f^G$ is defined over all of $\real^D$.

 \subsubsection{Moving to pseudo-likelihoods}
 \label{subsubsec:moving}
 
 We now turn to role of the Hamiltonian $H^L_{\Theta_L}$ of \eqref{eq:GibbsH2}
 in the context of  parameter estimation for the density $\varphi$  of \eqref{eq:Gibbs}. Unfortunately, if we retain it in its current form,  
parameter  estimation  by a method such as direct maximum likelihood is precluded by the fact that we neither have an analytic form for  $Z_{\Theta}$, nor is there any way to compute it numerically in any reasonable time frame.

The standard way around this problem, which we adopted in \cite{PNAS}, is the pseudolikelihood approach; e.g.\  \cite{Besag,chalmond}, which originated  in the context of  point cloud data with spatial dependence. In particular, this method exploits the inherent spatial Markovianess of a Gibbs distribution to replace the overall probability of, in our case, $\varphi_{\Theta}$, by the (semi-parametric, due to the inclusion of the non-parametric $\hat f^G$) pseudo-likelihood
\beq
\label{eq:pseudo1}
\tilde{L}^ K_{\alpha,\Theta}(\tx_N)  \definedas
\prod _{x\in\tx_N}  \hat\varphi_{\alpha,\Theta} (x),
\eeq
 where
 \beq
  \hat\varphi_{\alpha,\Theta} (x) 
  &=&  \notag
 \frac{
 (\hat f^G(x))^{\alpha}\times \exp \left(-\cL^K_\Theta  \left(x\,\big| \, \cN(x) \right)  \right)
  }{
\int _{\real^D}\hat f^G(z))^{\alpha} \times \exp \left(-\cL^K_\Theta  \left(z\,\big| \, \cN(x) \right)\right)\,dz}\\
 &=&  
 \frac{
 \exp \left(-\cL^K_\Theta  \left(x\,\big| \, \cN(x) \right) +\alpha\ln (\hat f^G(x))       \right)
  }{\int _{\real^D} \exp \left(-\cL^K_\Theta  \left(z\,\big| \, \cN(x) \right) +\alpha\ln (\hat f^G(z))   \right)\,dz}.
\label{eq:conditionalham1}
\eeq
Note that the main simplification of moving from the original Gibbs distribution to the pseudo-likelihood above is that the normalising constant/partition function now involves the product of $N$ integrals over $\real^D$, rather than a single, but numerically much more  demanding,  integral over $\real^{ND}$.

The remainder of the paper takes the model \eqref{eq:pseudo1} as given, and examines parameter estimation for it, how to use it to replicate point processes --  in particular, persistence diagrams  --  and how to use these tools to develop  statistically sound hypothesis tests in persistence based TDA. Evaluation of how well the procedures  work will be undertaken via  numerical experiments for specific examples.

\section{Modelling, estimating, and replicating persistence diagrams}
\label{sec:MERP}

\subsection{Persistence diagrams}
Although the development of the previous section was quite general, for the rest of the paper we shall concentrate on its application in the setting of persistence diagrams. 

The first consequence of this is that the point process of interest is now a subset of $\real^2$, restricted to the half plane above the $45^\circ$ line passing through the origin. 

In fact, mainly for coding convenience, given a persistence diagram with points $\{(d_i,b_i)\}_{i=1}^N$ we define a new point set
$\tx_N =\{x_i\}_{i=1}^N=\{(x^{(1)}_i,x^{(2)}_i\}_{i=1}^N$, defined by 
\beqq
x^{(1)}_i\ =\ d_i,\qquad x^{(2)}_i=b_i-d_i.
\eeqq
This (trivially invertible) transformation has the effect of moving the points in the original persistence diagram downwards, so that the diagonal line projects onto the horizontal axis, but still leaves a visually informative diagram, which we called the projected persistence diagram, or PPD, in \cite{PNAS}.   The points of the PPD lie in  the upper half plane $\real\times\real_+$, and it is for these points that we want to apply the model of the previous section.

Since the PPD is restricted to a half plane, we shall assume that the density $f^G$ of \eqref{eq:Gibbsalpha} is supported there. Similarly, the empirical density $\hat f^G$ of \eqref{eq:kernel} needs to be replaced by
\beq
\label{eq:kernel-new}
\bar f^G(x) \ \definedas \  \frac{\hat f^G(x) {\bf 1}_{\real\times\real_+}(x)}{\int_{\real\times\real_+} \hat f^G(x)\,dx}.
\eeq
Replacing $\hat f^G$ in the definition \eqref{eq:conditionalham1} of $\hat \varphi_{\alpha,\Theta}$ by $\bar f^G$ leads to the pseudolikelihood 
\beq
  \bar\varphi_{\alpha,\Theta} (x) 
  &=&  
 \frac{
 \exp \left(-\cL^K_\Theta  \left(x\,\big| \, \cN(x) \right) +\alpha\ln (\bar f^G(x))       \right)
  }{\int _{\real^D} \exp \left(-\cL^K_\Theta  \left(z\,\big| \, \cN(x) \right) +\alpha\ln (\bar f^G(z))   \right)\,dz}.
\label{eq:finalpseudo}
\eeq
with which we shall work from now on.

\subsection{Parameter estimation}

There are four quite different classes of parameters to estimate in the likelihood \eqref{eq:pseudo1}: the bandwidth parameters $\Sigma$ or 
$\eta$, the maximal neighbourhood size $K$, the weighting $\alpha$ of the global shape, and the local interaction 
parameters $\theta_1,\dots,\theta_K$. There are also a number of difficult numerical issues associated with their estimation, which we now describe.

As for the bandwidth parameters of \eqref{eq:kernel}  for estimating $f^G$, in practice we found two different scenarios. When the persistence diagram had a locally `regular' (in a very imprecise sense) shape, then it sufficed to take $\Sigma$ to be diagonal, so that only $\eta$ required estimation. As noted above, asymptotic theory gives that $\eta$ is optimally chosen to be approximately $\max (\sigma_1,\sigma_2)N^{-1/6}$, where $\sigma_j^2$ is the empirical variance of the $x^{(j)}_i$.  However, in practice (not surprisingly, since our samples -- persistence diagrams -- did not have `asymptotically many' points) we used this figure as a guide, followed by ad hoc decisions. 

On the other hand, for situations in which there were important subtleties in the persistence diagram, such as in the case of the Gaussian random field excursion filtration of Section   \ref{subsec:Gauss-new} (cf.\ Figure \ref{fig:GaussH}), estimation of the full, non-isotropic, bandwidth matrix $\Sigma$ was preferable. Since asymptotic methods require particularly large sample sizes in this case (\cite{SilvermanDensity,WassermanAll}) we found the data driven method of \cite{Duong}, as implemented in the R package {\it Hpi}, worked well.

%
%
%
%

As for choosing the parameter $K$, the number of nearest neighbourhood regions for the local part of the Hamiltonian, prior experience, not only  from  \cite{PNAS} but from decades of modelling in Statistical Mechanics, implies that it suffices to consider $K\leq 3$. Choosing among these values, and deciding which, if any, of the neighbourhoods should be excluded from the Hamiltonian (i.e.\ by setting their $\theta_L$ value identically zero) is then easy to do via  standard statistical procedures such as those based on the Akaike or Bayes Information Criteria, AIC and BIC, etc.  (cf.\ \cite{burnham}.)  Of course, this procedure implies that, given an estimate $\bar f^G$ as above,  for each such model -- i.e.\ for each choice of non-zero $\theta_L$ --  we need to maximise the pseudo-likelihood \eqref{eq:pseudo1} over $\alpha$ and the non-zero $\theta_L$.

This maximisation was carried out in a two stage procedure: We searched over $\alpha$ by the bisection method, and maximised the pseudo-likelihood for the $\theta_L$'s using a standard non-linear optimisation procedure. (To be more specific, the Matlab routine {\it fminunc}.)   During this stage, the integral in the denominator of \eqref{eq:conditionalham1} requires computation, which we carried out numerically by the trapezoid method. (To be more specific, again, this was done on a $101\times 101$ grid. For the $x^{(1)}$ variable, the range extended to at least 4 standard deviations, in each direction, from the mean, where means and standard deviations were calculated from the original data. For the positive variable $x^{(2)}$, the range went from 0 to at least the mean plus 4 standard deviations.)

As an (important) aside,  note that,  after considerable experimentation, we found that it sufficed to restrict $\alpha$ to the interval  $[0,3]$.

\subsection{Replicating persistence diagrams}
\label{sec:replications}

Given a pseudolikelihood as in the previous section (with known or  estimated parameters), generating simulated replications of the associated point set  via Markov Chain Monte Carlo (MCMC)  is standard, and so we include now only the briefest of descriptions. (As in \cite{PNAS} we take a  Metropolis-Hastings MCMC approach, and refer the  reader to \cite{RobertCasella,Handbook} for technical background.  In particular, see \cite{RobertCasella}  Sec.\ 10.3.3, in which the  approach we take is called `Metropolis-within-Gibbs' and its  properties are discussed.)

Thus, all we need to do is to show how, given an initial version $\tx^0$ of our point set/persistence diagram, we make the step from version $\tx^n$ to the next stage, $\tx^{n+1}$. This, as usual, is done one point at a time.

Firstly, we choose a new point in $x^*\in\real\times\real_+$ according to the density $\bar f^G$ of \eqref{eq:kernel-new}. Note that this choice  depends only on the initial state $\tx^0$, and so can be done in advance, before entering the MCMC procedure. (Of course, one could replace this with using $\tx^n_N$ at the $n$-th step rather than using only $\tx^0$, but we found this to be numerically less stable and produce poorer simulations. Also, as mentioned in Section \ref{subsubsec:choosingprior}, if there are better estimates of $f^G$ based on prior information, then these should be used in place of  $\bar f^G$.)  To choose $x^*$ we used a standard inverse transform method \cite{RobertCasella,Handbook}. Details are given in Appendix \ref{sec:appendixglobalshape}.

Next, for two points $x,x^*\in\real\times\real_+$  define an `acceptance probability', according to which $x\in\tx_N$ is replaced by $x^*$, leading to the updated set $\tx_N^*$, as
\beqq
\rho \left(x, x^*\right)
& =&  \min \left\{1,\ \frac{(\bar f^G(x^*))^\alpha e^{-\cL^K_\Theta(x^*|\cN_x)}  \cdot \bar f^G (x)}{
 (\bar f_G(x))^\alpha e^{-\cL^K_\Theta(x|\cN_x) }  \cdot \bar f^G \left(x^*\right)}\right\} \\
& =&  \min \left\{1,\  e^{-\cL^K_\Theta(x^*|\cN_x)  +\cL^K_\Theta(x|\cN_x)}\left({ \bar f^G (x)}/{\bar f^G (x^*)}\right)^{1-\alpha}\right\} .
\eeqq
Recall from \eqref{eq:pseudo1} and \eqref{eq:conditionalham1} that the current state of the entire system, $\tx_N$, appears in the $\cL$ terms  via the  sums over  neighbourhoods in the local Hamiltonian.

The one-step replacement  algorithm can now be described, as in  Algorithm \ref{MCMC:algorithm}.

\begin{algorithm}
\caption{MCMC step updating diagram for $\tx_N$}
\label{MCMC:algorithm}
\begin{algorithmic}[1]
\State  $k =0$
 \State $k \gets k+1$
\State Choose $x^*$ according to $\bar f^G$
\State Compute $\rho(x_k,x^*)$
\State Choose $U$ a standard uniform variable on $[0,1]$

\If{$U<\rho(x_k,x^*)$} set $x_k =x^*$
\EndIf
\If{$k<N$}  go to Step 2
\EndIf
\end{algorithmic}
\end{algorithm}

To obtain  approximately independent copies of the point set/persistence diagram, the procedure dependents on three  parameters, $n_b$, $n_r$ and $n_R$.  Starting with the original diagram, run the algorithm for a burn in period. Then, starting with the final result from the burn in, run the algorithm a further $n_b$ times, choosing the last output of this block of $n_b$ iterations as the first simulated diagram. Repeat $n_r$ times, each time starting with the most recently simulated diagram; viz.\ the output of the previous block. Finally, replicate the entire procedure $n_R$ times, so that we have a total of $n_r\times n_R$ simulated diagrams.
The optimal choice of $n_b$, $n_r$ and $n_R$ and typically depends on the specific problem, and is discussed  in the examples below and in \cite{PNAS} SI Appendix (Sec.\ 2.1).

Recall that the procedure above is for {\it projected} persistence diagrams, but conversion from these back to the usual diagrams is trivial, via the mapping {$x_i\to (x_i^{(1)}, x_i^{(1)}+x_i^{(2)})=(d_i,b_i)$ }.


\subsection{RST and statistical inference for persistence diagrams}
\label{sec:inference-new}

The basic idea of, and need for, RST (Replicating Statistical Topology), as introduced in \cite{PNAS} is quite simple: In many applications of TDA only one original data set is available, and so only one persistence diagram. Replicating the experiment may be unfeasible, and as a result making statistically meaningful statements about an observed persistence diagram is essentially impossible.

The idea behind RST is that, given a single persistence diagram, it can be modelled, estimated, and simulated as described in the preceding sections, and on each simulated diagram any number of meaningful statistics calculated. The MCMC replications of the diagrams thus give a sequence of (almost) independent values of the statistic, enabling standard statistical analysis. Precisely which statistics one chooses will depend on the hypotheses one wants to test, and in the examples of the following section we shall look at a number of possibilities.

The main difference between the present paper and \cite{PNAS} is in the choice of model. The model there was of a similar nature, in that there was a `local' Hamiltonian that described nearest neighbour interactions much as in the model of the current Section  \ref{sec:modelling-new}, but there was no attempt to capture the global structure via a term akin to the density $\bar f^G$. To nevertheless place some restrictions on the overall shape, there was a term in the Hamiltonian conserving the centre and the spread of the simulated diagrams. (For completeness, the earlier model is described in Appendix \ref{sec:appendixearliermodel}.)  While that model  was a good first step, the current model does a much better job of replicating persistence diagrams, as the examples  of the next section will all show. In particular, the stationary distributions of the MCMC procedure are much closer to the true distributions of the persistence diagrams, where the latter are estimated by simulating the underlying experiment to obtain replications. This procedure, of course, is easy to carry out for a simulation study, but, as just noted above, typically impossible or prohibitively expensive in  applications.

An alternative, and earlier, approach to replicating persistence diagrams was based on resampling methods such as the bootstrap and jacknife, e.g.\  \cite{chazal,fasy,RobTurner}. Resampling can be done either at the level of the basic experimental data, in which case a new persistence diagram needs to be calculated for each resampled data set, or at the level of the persistence diagram itself. Since we have already argued in \cite{PNAS} that RST compares favourably with these procedures, and since the version of RST presented here improves on the earlier one, we shall say no more about these methods in the current paper, other than to make two points in their favour. Firstly, 
they can be computationally much faster than the method that we propose, for which the parameter estimation is often delicate and the MCMC can take considerable computer time. (Typical time differences are about a factor of 10, from seconds to minutes.)  However, this computational advantage is  typically only for the case when the subsampling is at the level of the persistence diagram. When it is at the level of the experimental data, each new data set requires the computation of persistence diagrams, and this is usually time consuming.  The second advantage of the resampling approach is that, at least at this point of time, it is easier to fit it into an existing framework of  asymptotic statistical theory, which provides rigorous information about its behaviour, as various parameters tend to infinity or zero. Whether or not one feels that this is important for applications is, to a large extent, a question of taste, best left to the individual reader.

\section{Examples}
\label{sec:examples-new}

Given the lack of any theory (other than asymptotic, and even then in very special cases) about the distributions of persistence diagrams it is, {\it a fortiori},  impossible to bring theorems about the efficiency of the approach of the previous two sections for modelling them. Furthermore, since in any application 
``the proof of the pudding is in the eating\footnote{Originally, ``al fre\'ir de los huevos lo ver\'a" (you will see when the eggs are fried). See \cite{pudding}.}", establishing the value of the preceding sections can really only done by seeing how they work in practice.

Thus, in the following subsections, we shall investigate a number of examples, in each one concentrating on these aspects of model that are most relevant for the example in question. For example, in the two motivating examples we have already seen -- persistence diagrams arising from Gaussian excursion sets and sampling from non-concentric circles -- the emphasis will be on reproducing, by simulation, high  quality samples with  (a reasonably close approximation to) 
the true distribution (as estimated by simulation of the underlying phenomenon). In other examples we shall combine these methods with statistical inference to see how well our proposed procedure of Replicating Statistical Topology works in a hypothesis testing scenarios. {\it En passant}, we shall show how and why the current model improves on that in \cite{PNAS}.

\subsection{Gaussian excursion sets}
\label{subsec:GaussExample}

Returning now to the setting of Section  \ref{subsec:Gauss-new}, our aim for this example is to see how well we can approximate the true distribution of the $H_1$ persistence diagrams of Gaussian excursion sets.

Specifically, we take the Gaussian random field on $[0,1]^2$ with covariance \eqref{eq:covariancefunction} (with $b=100$ and $a =1$)  and  our  main aim is to see if we can capture the shape of Figure \ref{fig:GaussH}(b). To do this, we took the  600 simulations that led to this figure,  along with their  excursion set persistence diagrams. We shall call these the `original' diagrams. Three examples were shown in Figure \ref{fig:r1}. For each of the original diagrams we fitted a model with pseudolikelihood \eqref{eq:finalpseudo}, and then ran a MCMC simulation as described in Section \ref{sec:replications}, running it for 1,000 steps, so that in the final analysis we had $600\times 1,000=600,000$  diagrams from the MCMC procedure. We shall call these the `simulated' diagrams.

In fitting a model, we found an interesting phenomenon, summarised in Table \ref{table:Models}: In most of the 600 cases, a simple model with $K=1$ (i.e.\ only one neighbourhood term in the Hamiltonian) was chosen as optimal by AIC and BIC criteria. A natural (at least  `natural', {\it a fortiori}) explanation for this comes from the numbers
in Table \ref{table:correlations}.

\begin{table}[h!]
\begin{center}
\begin{tabular}{|l|rrrrrrr|}
 \hline  Active $\theta$ parameters &$\theta_1$ & $\theta_2$ & $\theta_3$ &  $(\theta_1,\theta_2)$  &  $(\theta_1,\theta_3)$   & $(\theta_2,\theta_3) $& $(\theta_1,\theta_2,\theta_3)$   \\  \hline 
Count &430 &0 & 0&71 &79 &0 &20 \\
Percentage & 72.7&0 &0 &11.8 &13.2 & 0& 3.3 
\\ \hline
\end{tabular}
\caption{\footnotesize{Numbers and percentages of models of different types for the Gaussian excursion set example.}}
\label{table:Models}
\end{center}
\end{table}

What Table \ref{table:correlations} shows are very strong negative correlations between $\alpha$ and $\theta_2$ and $\theta_3$, and between $\theta_1$ and $\theta_2$ and $\theta_3$. (All  were statistically significant with $p$-values less than $10^{-6}$.) The  negative correlations 
between the estimates of $\alpha$ and the $\theta_i$ are natural in terms of our original discussion of the model, and capture the `competition' between its global and local aspects. The second set indicates a similar competition, and indicates  that by playing with one of the parameters it is possible to reduce the effect of another. Precisely quantifying this effect is what the AIC and BIC indices do, and is presumably behind the preponderance of the model with $K=1$ in Table  \ref{table:Models}.  (As an aside, we note that in this model the correlation between $\alpha$ and $\theta_1$ was -0.1423, with a $p$-value of $0.0007$.)

\begin{table}[h!]
\begin{center}
\begin{tabular}{|l|rrr|}
\hline
 & $\theta_1$ & $\theta_2$ & $\theta_3$ \\  \hline 
$\alpha$ & -0.0064 \ \    & -0.3514*   & -0.4108*  \\
$\theta_1$    &  & -0.2409*   & -0.1822*    \\
$\theta_2$  &   &   & -0.0776\ \ 
\\ \hline
\end{tabular}
\caption{\footnotesize{Correlations between parameter estimates for the Gaussian excursion set  example for the model with $K=3$; i.e.\ with  three $\theta$ parameters. Statistically significant correlations are marked with asterisks.}}
\label{table:correlations}
\end{center}
\end{table}

From the point of view of general Statistical Mechanics, the discipline underlying Gibbs distributions, the preference for the $K=1$ model is consistent with experience. Nearest neighbour models (i.e.\ Hamiltonians with only a single interaction term) are almost ubiquitous as models of physical phenomena, and higher level interactions are only rarely needed.

\begin{figure}[h!]
        \subfigure[]
        {
        \includegraphics[scale=.08]{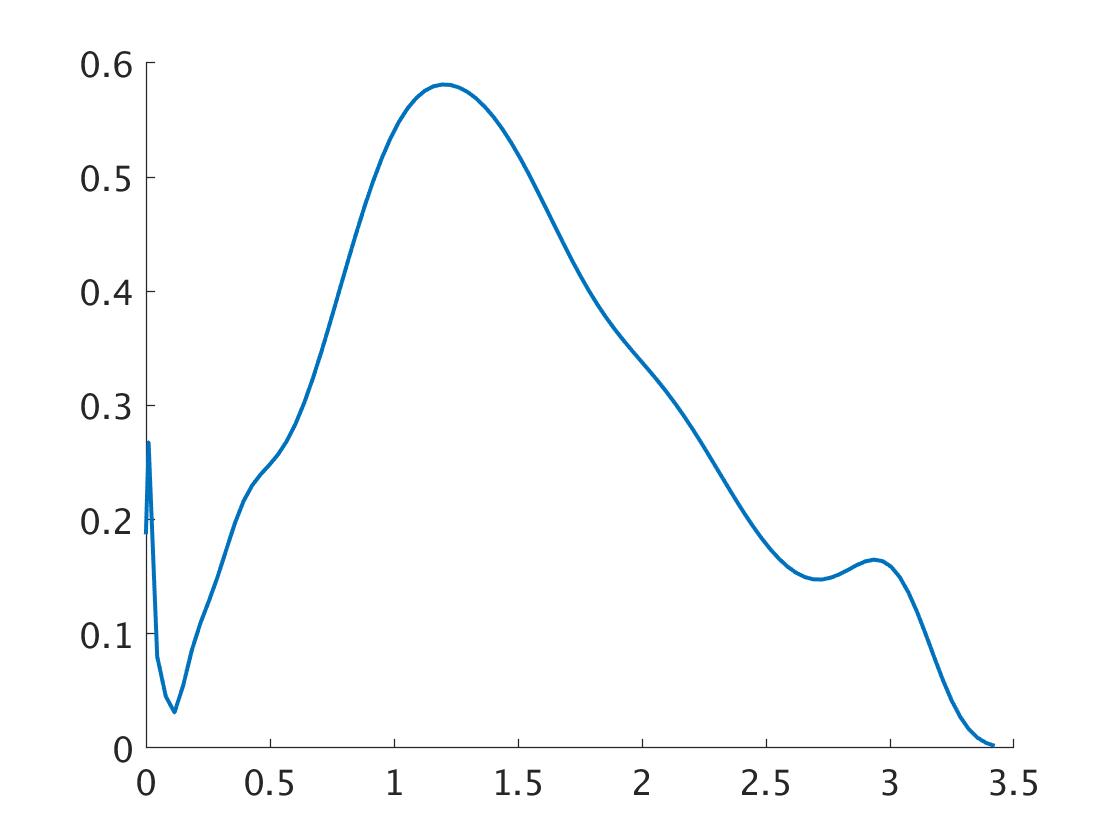}
    } 
    \subfigure[]
    {
        \includegraphics[scale=.08]{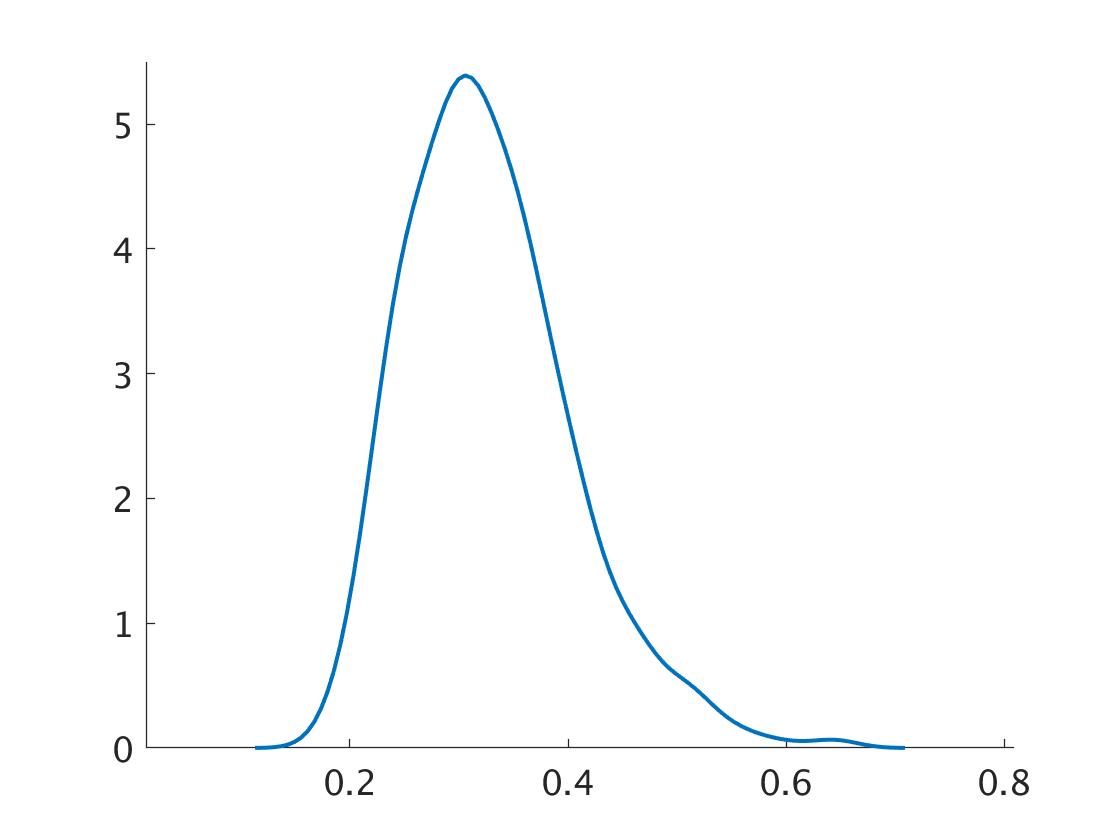}
    }
    \\
   \subfigure[]
    {
        \includegraphics[scale=.08]{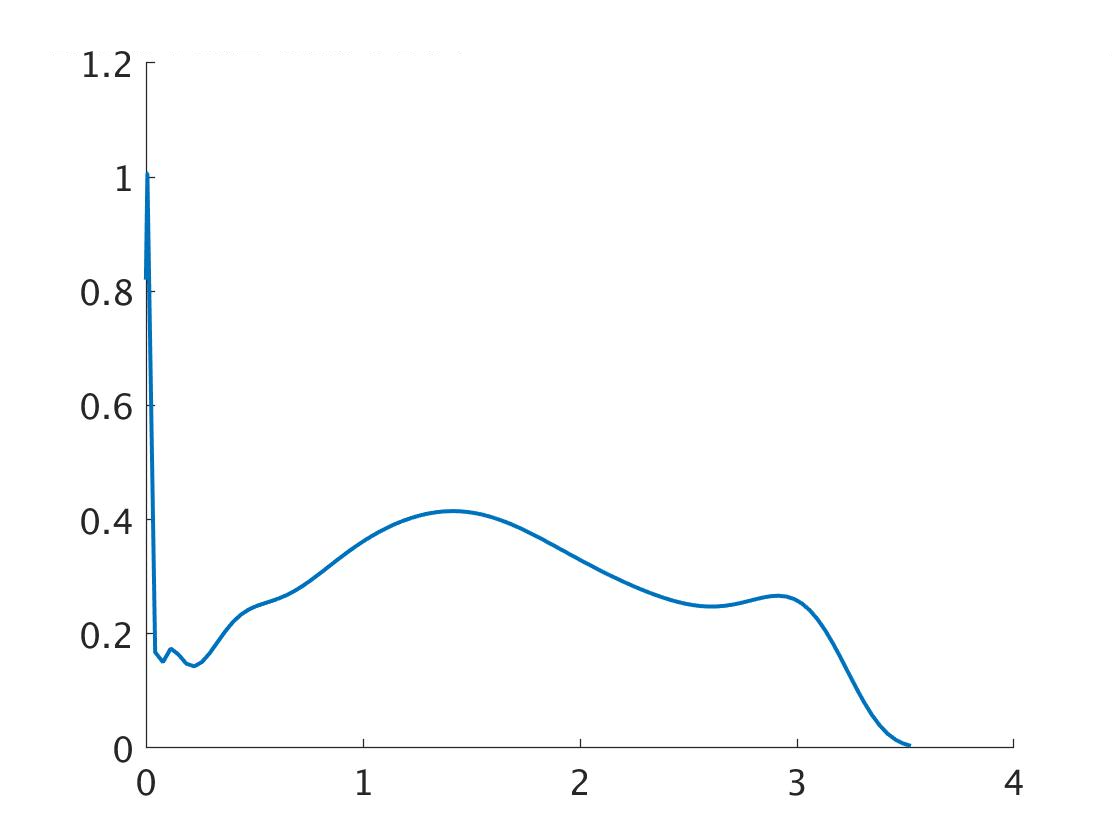}
    }
    \subfigure[]
    {
        \includegraphics[scale=.08]{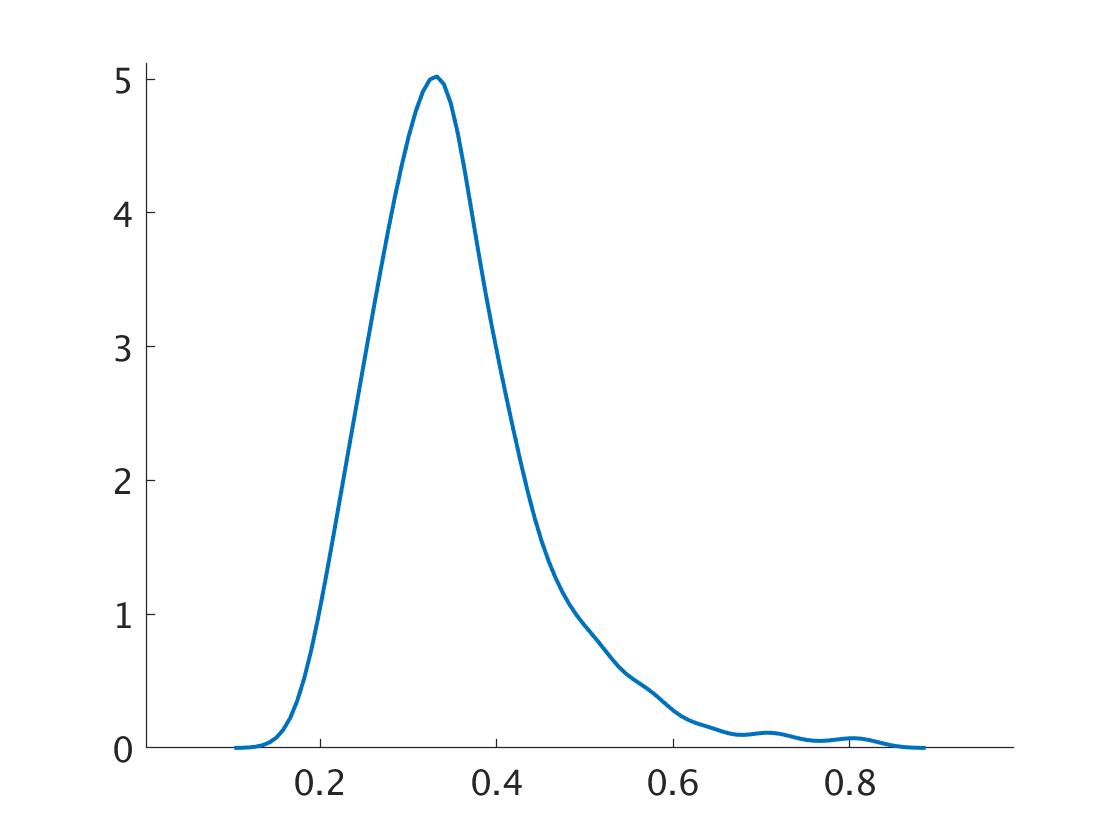}
    }
     \subfigure[]
    {
        \includegraphics[scale=.08]{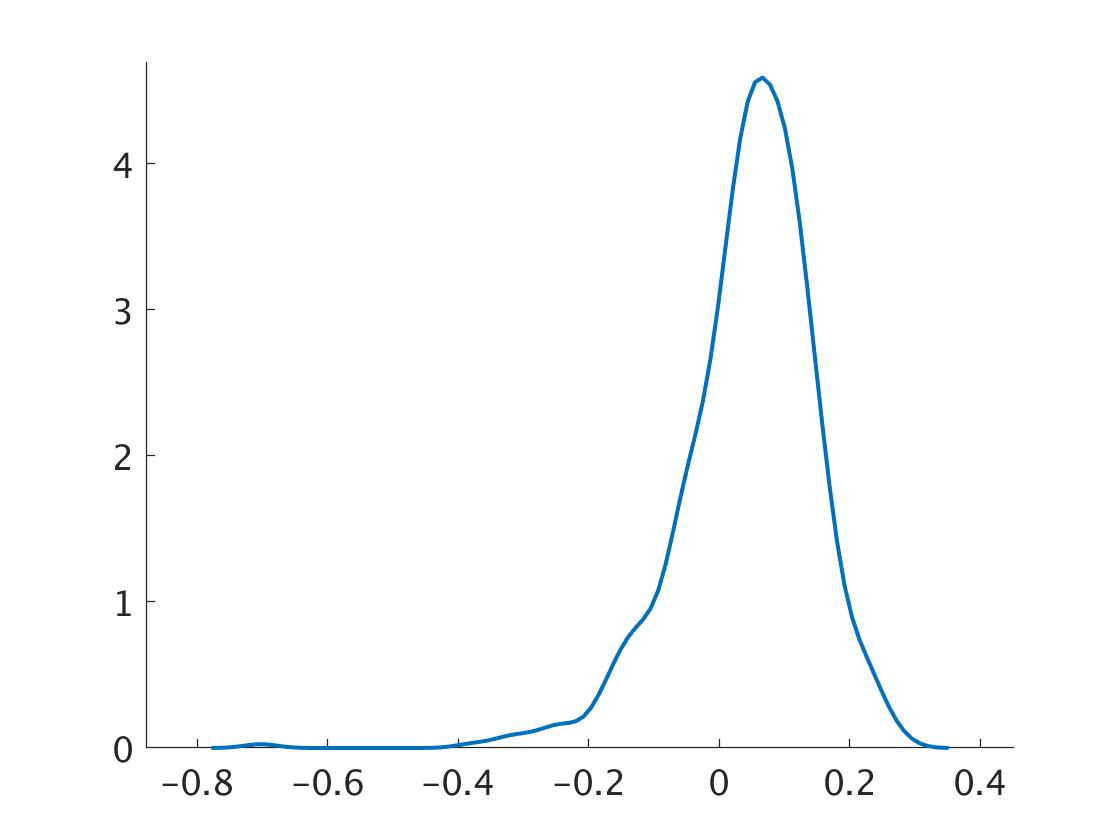}
    }
    \subfigure[]
    {
        \includegraphics[scale=.08]{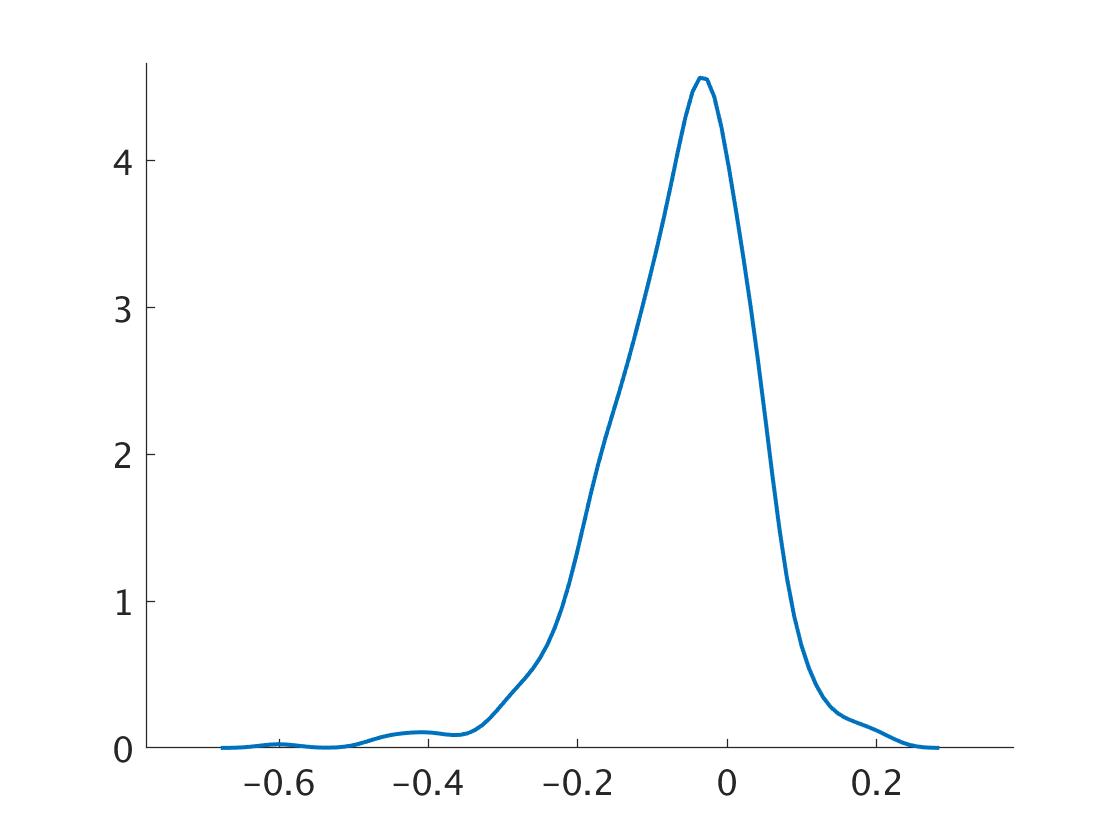}
    }    
  \caption{\footnotesize
     Empirical distributions of $\alpha$ (a) and $\theta_1$ (b) in the model with $K=1$ for the Gaussian excursions persistence diagrams. For the model with $K=3$, the densities are for  $\alpha$ (c),  $\theta_1$ (d), $\theta_2$ (e) and $\theta_3$ (f).    The values of $\theta$ have all been normalized by the number of points in the diagram.  }
\label{fig:r33}
\end{figure}

Further information on the stability of the  parameter estimation is provided in Figure \ref{fig:r33}, which shows
 the empirical densities of the parameter estimates  for all 600, when we fix the model: Results for the  model with  $K=1$ are shown in Panels (a) and (b), and  for the case  $K=3$ in Panels  (c)-(f).  In the case of the $\theta_j$, we normalise its value for each diagram by   the number of points, $n$, in the diagram. For reasons that are currently unclear to us, this normalisation not only produces a slightly smoother empirical density, but  also makes comparison between different examples more natural. In particular, whereas in the current example the number of points in each diagram ranges from 11 to 23, with  about 90\% 
 in the range 14--19,  in later examples (sampling from non-concentric or concentric circles, and from a 2-sphere) the numbers are  often much higher. Estimates of the $\theta_j$ grow with this number, whereas the $\theta_j/n$ remain in the range  $(-5,5)$.  Overall, with this normalisation, the empirical densities for the 
$\theta_j/n$ estimates, along with the corresponding ones for the later examples, are an invitation to believe in an asymptotic normal limit for the estimates. Such limits, while commonly observed in settings like the current one, are infamously hard to prove. 

Another interesting phenomenon, shown in Figure \ref{fig:r33} (a) and (c) is the fact that while the typical choice of the parameter $\alpha$ is around 1.5, the estimation procedure is not totally averse to setting $\alpha=0$ (20 cases out of 600).That is,  the model is prepared to do without the global restrictions imposed by the density $\bar f^G$, but  does so only rarely. 

\begin{figure}[h!]
    \centering
        \subfigure[]
        {
        \includegraphics[scale=.11]{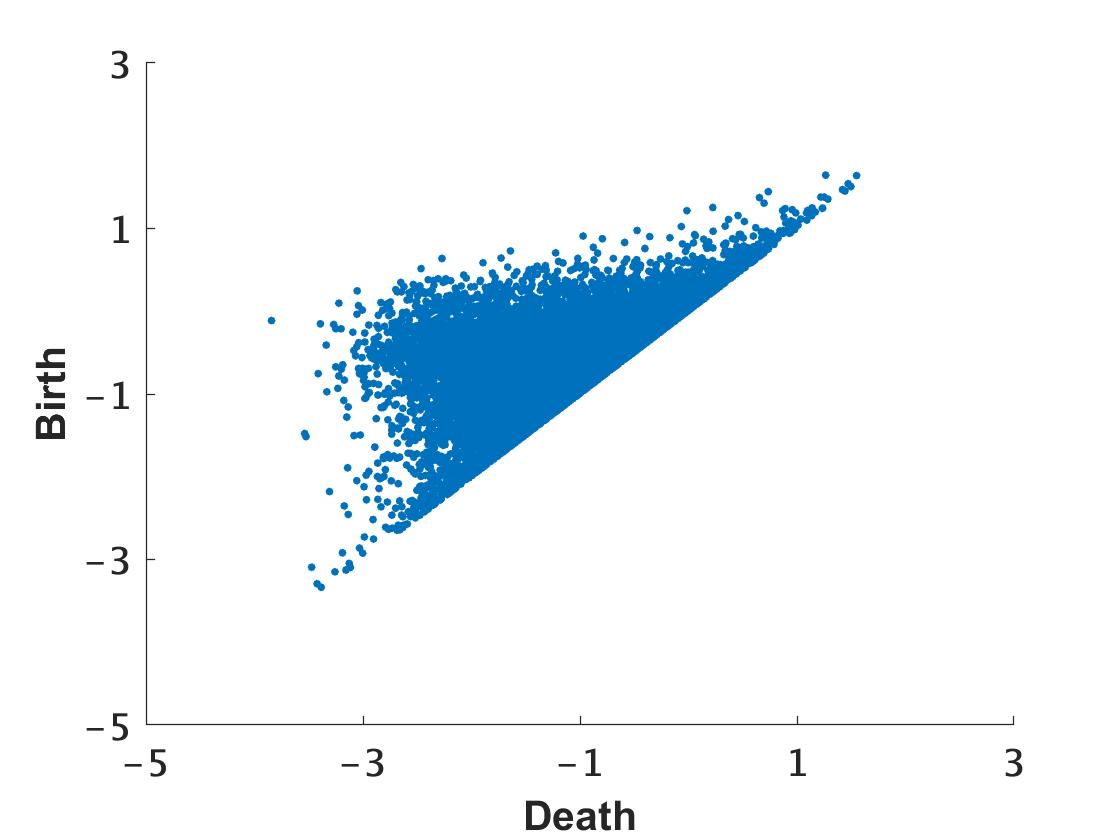}
    }
    \subfigure[]
    {
        \includegraphics[scale=.11]{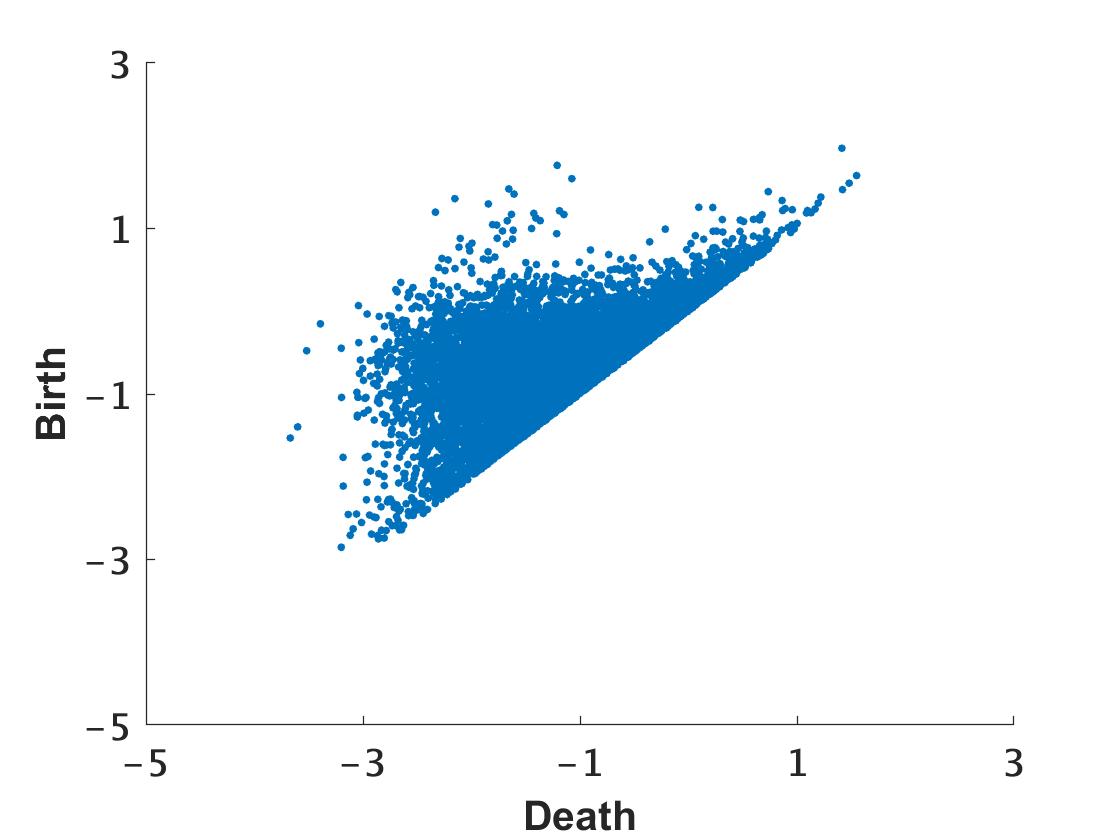}
    }
   \subfigure[]
    {
        \includegraphics[scale=.11]{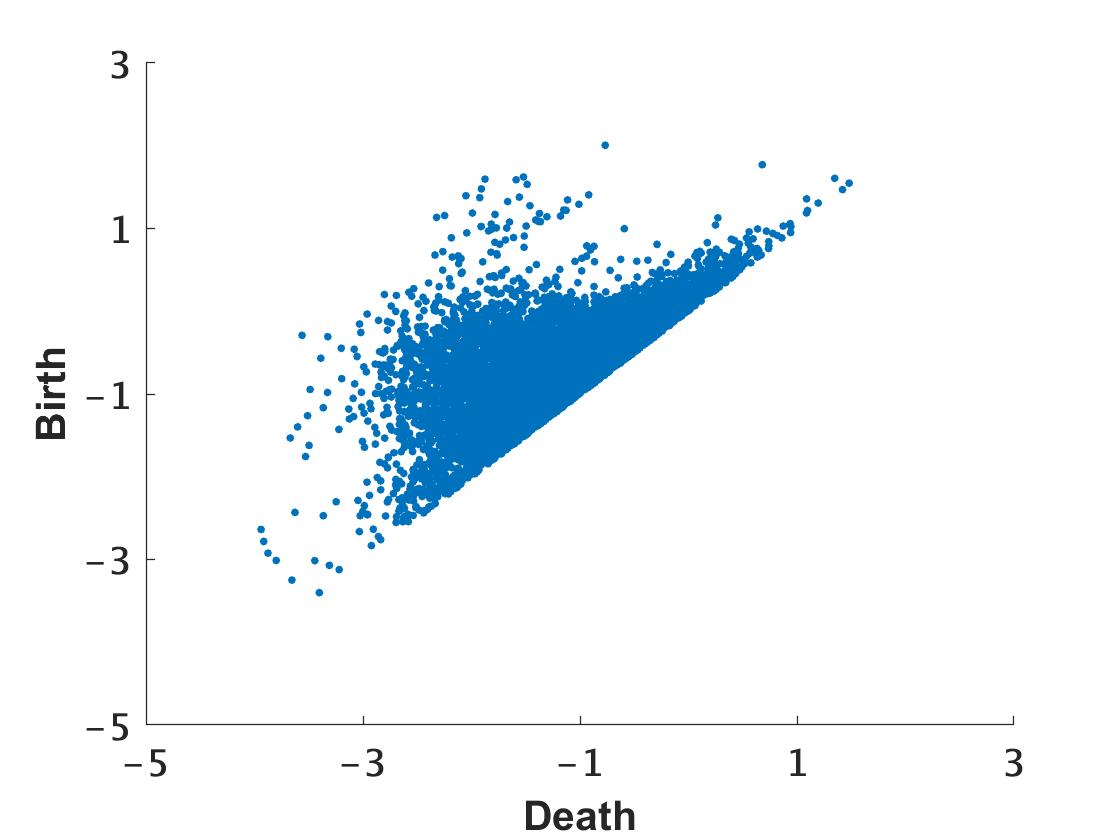}
    }
    \\
    \subfigure[]
    {
        \includegraphics[scale=.11]{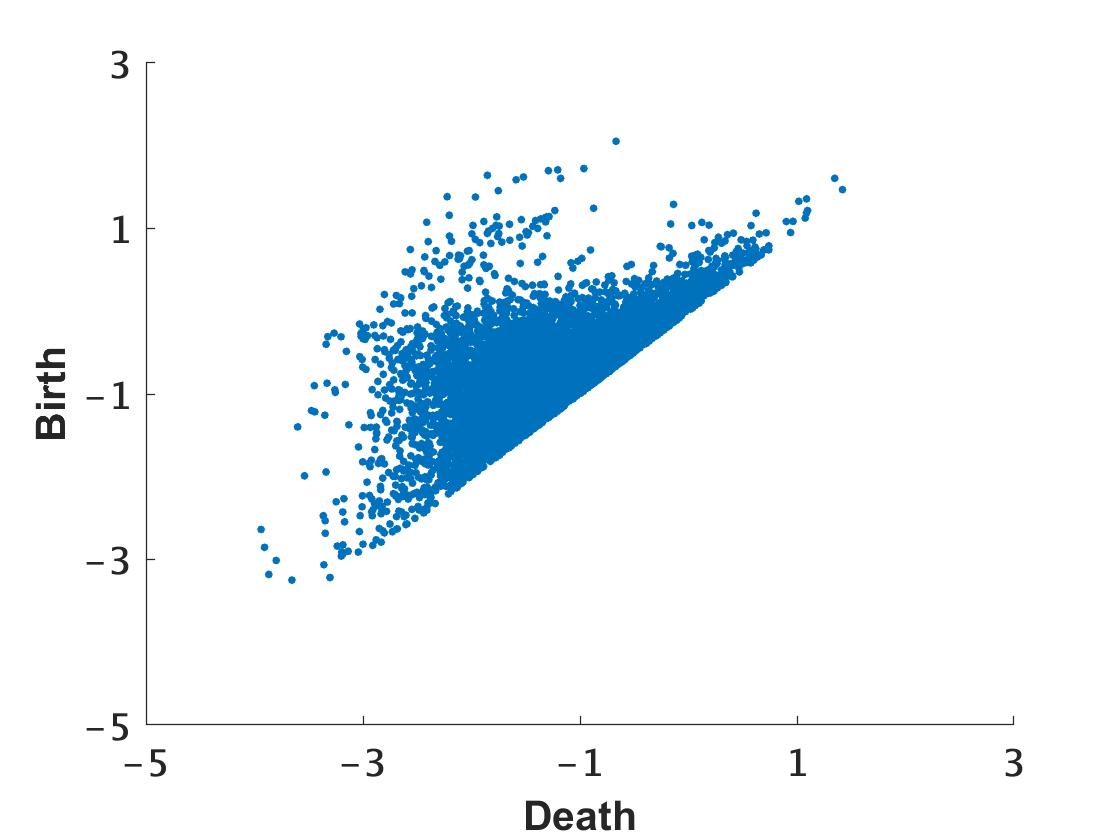}
    }
    \subfigure[]
    {
        \includegraphics[scale=.11]{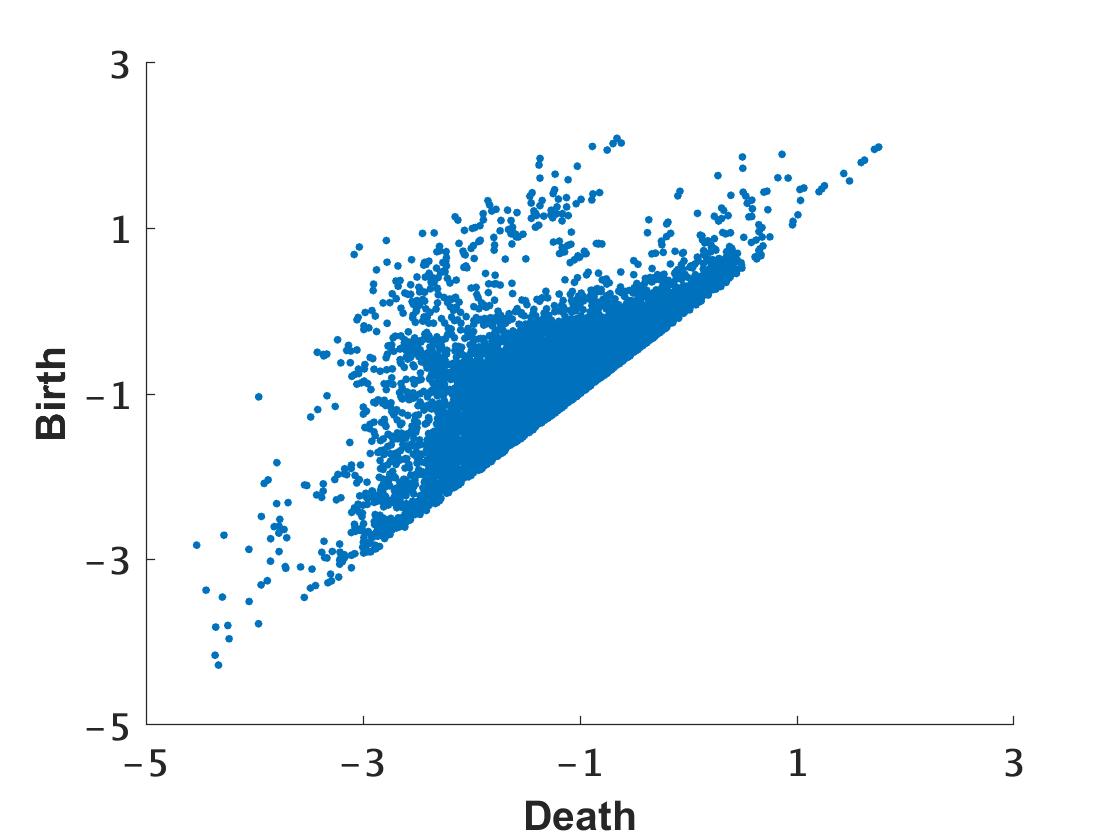}
    }
    \subfigure[]
    {
        \includegraphics[scale=.11]{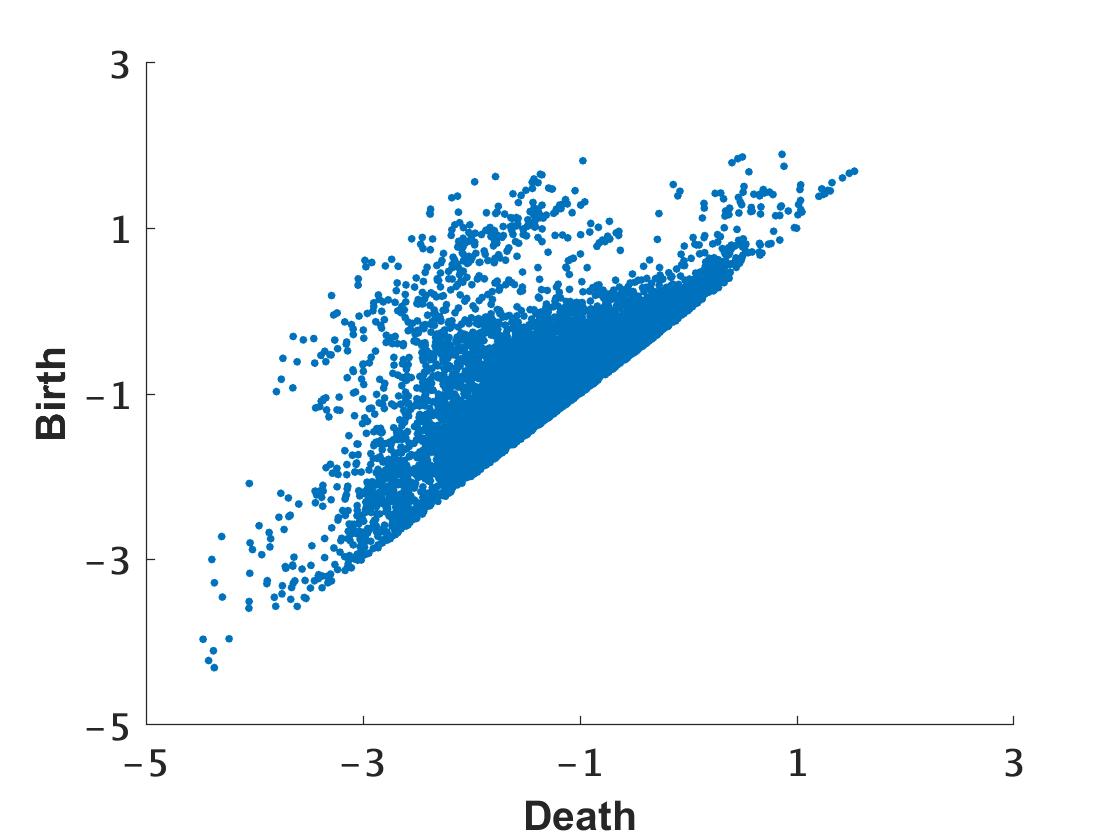}
    }    
  \caption{\footnotesize
     (a) is a superposition of 600 `original' Gaussian excursion set persistence diagrams.  (b)--(f) are similar superpositions from 600 MCMC simulations, at steps 10, 50, 100, 500, 1,000.}
\label{fig:r2}
\end{figure}

We now turn to what was our main motivation for this example: seeing if  our model can reproduce the point clouds of Figure     \ref{fig:GaussH}. 
To this end, we took each of the simulations of `original' diagrams just described, and, on each one, with the parameters estimated as above, ran an MCMC simulation as in  Section \ref{sec:replications}. Figure \ref{fig:r2} summarises the results, superimposing the diagrams of  600  such simulations, at 10, 50, 100, 500 and 1,000 steps into the simulation, along with the diagrams for the original diagrams.

Two facts are immediately noticeable. The first is that the general shape in Figure    \ref{fig:GaussH}    is, in fact, well preserved in the early stages of the simulations.  The second is that, as the MCMC progresses, a cloud of points separates from the main collection, spreading out in the directions of lower and higher death times. While visually striking, this phenomenon is not as remarkable as it first seems. For example, define a `vertical outlier' in a diagram to be a point for which its distance from the diagonal is at least $Q_3 + 5(Q_3-Q_1)$, where $Q_1$ and $Q_3$ are the first and third quartile of all distances in the diagram. Then Table \ref{table:outliers} shows the distribution of  the 18 diagrams out of 600 which contain vertical outliers,  after 500 steps of the MCMC algorithm (this was the  case with the largest number of outliers among those in Figure \ref{fig:r2}).  From the point of TDA, in which many statistical tests are based on vertical outliers (which measure the bottleneck distance of a diagram from the `zero' diagram, with all points on the diagonal) the implication is that these outliers are going to have minimal impact on most statistical tests and, even when there is an impact, it will tend to lead to conservative tests.

\begin{table}[h!]
\begin{center}
\begin{tabular}{|l|cccccccccccc|}
\hline
Number of outliers & 0 & 1  & 2 & 3 &4  & 5 & 6 &7  &8  & 9 & 10  & $>10$   \\ \hline
 Number of diagrams & 582 &  2 &10  & 4 & 0  & 0  & 1 &0  &0  & 0 & 1 & 0   
\\ \hline
\end{tabular}
\end{center}
\caption{{\footnotesize  Numbers of diagrams with a specific number  of  outliers among the 600 simulations of Figure \ref{fig:r2} (e) (i.e.\ at 500 MCMC steps).}}  
\label{table:outliers}
\end{table}

The same general comments can be made about the `horizontal spread' in the later diagrams of Figure \ref{fig:r2}. In any case, in applying these procedures in a RST/TDA setting, it is clear that taking 100 steps into the MCMC procedure in order to produce random perturbations of a given diagram, with similar statistical characteristics, is safe. 

\subsection{Non-concentric  circles}
\label{sec:twononc}

We now return to the example of Section \ref{subsec:2circles-new}, coming from sampling from two non-concentric circles. In particular, we were motivated there by the collection of $H_0$ persistence diagrams in Figure   \ref{fig:TwoCirclesPH:rob}, which showed three almost disjoint clouds of points; one indicating one of the circles (the other circle being represented by the removed `points at infinity')   and two being essentially noise at two scales, corresponding to the different sized circles.

We simulated 100 samples, thus providing 100 `original' persistence diagrams, one of which is in Figure \ref{fig:farcircle}(c). For each of these we estimated a model, as in the previous section, taking $K=3$ for the maximum number of possible neighbourhoods contributing to the Hamiltonian. In doing so, as in the Gaussian excursion example of the previous section, we found very high correlations between the parameter estimates, with $\rho(\theta_1,\theta_2)= -0.6772$ and $\rho(\theta_1,\theta_3)= -0.6948$, both  
statistically significant with $p$-values of order $10^{-14}$ and $10^{-15}$, respectively.  Thus it was not surprising that the optimal model (via considerations of AIC and BIC) overall was often that with only $\theta_1$ (i.e.\ $K=1$) and so we adopt this for the remainder of this example.

\begin{figure}[h!]
    \centering
       \subfigure[]
    {
      \includegraphics[width=2.4in, height=1.7in]{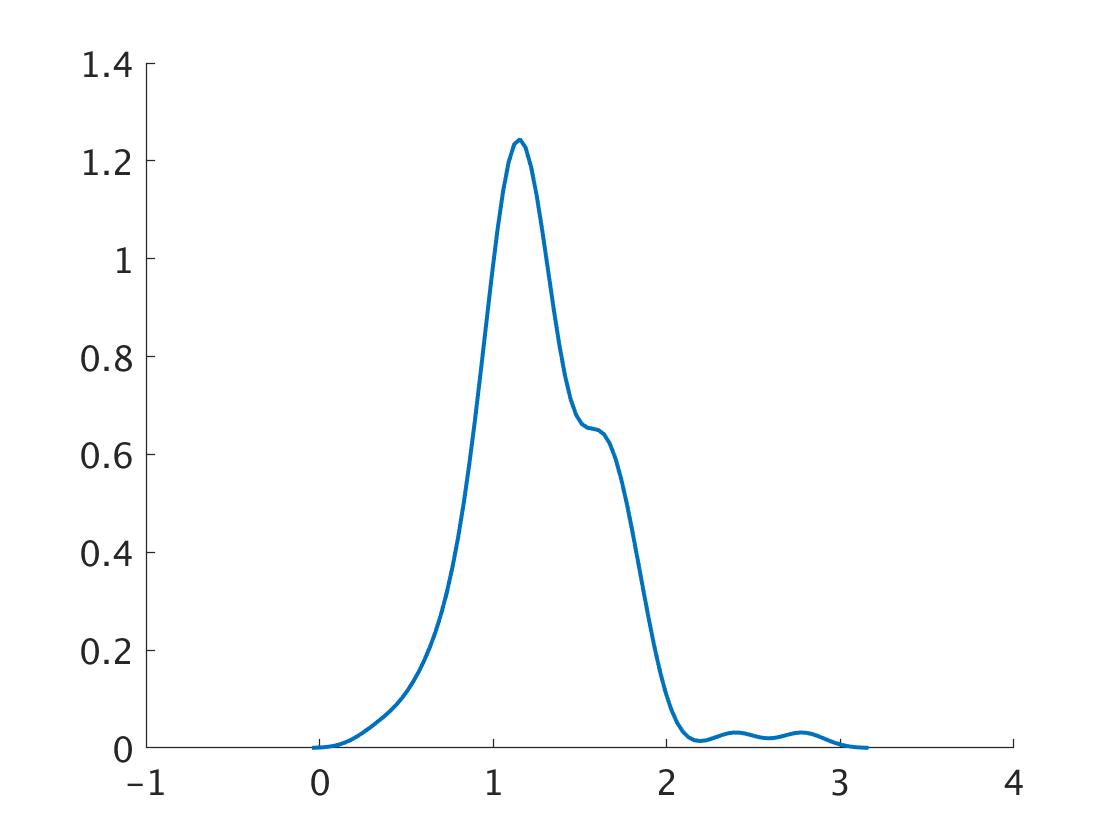}
    }
       \subfigure[]
    {
        \includegraphics[width=2.4in, height=1.7in]{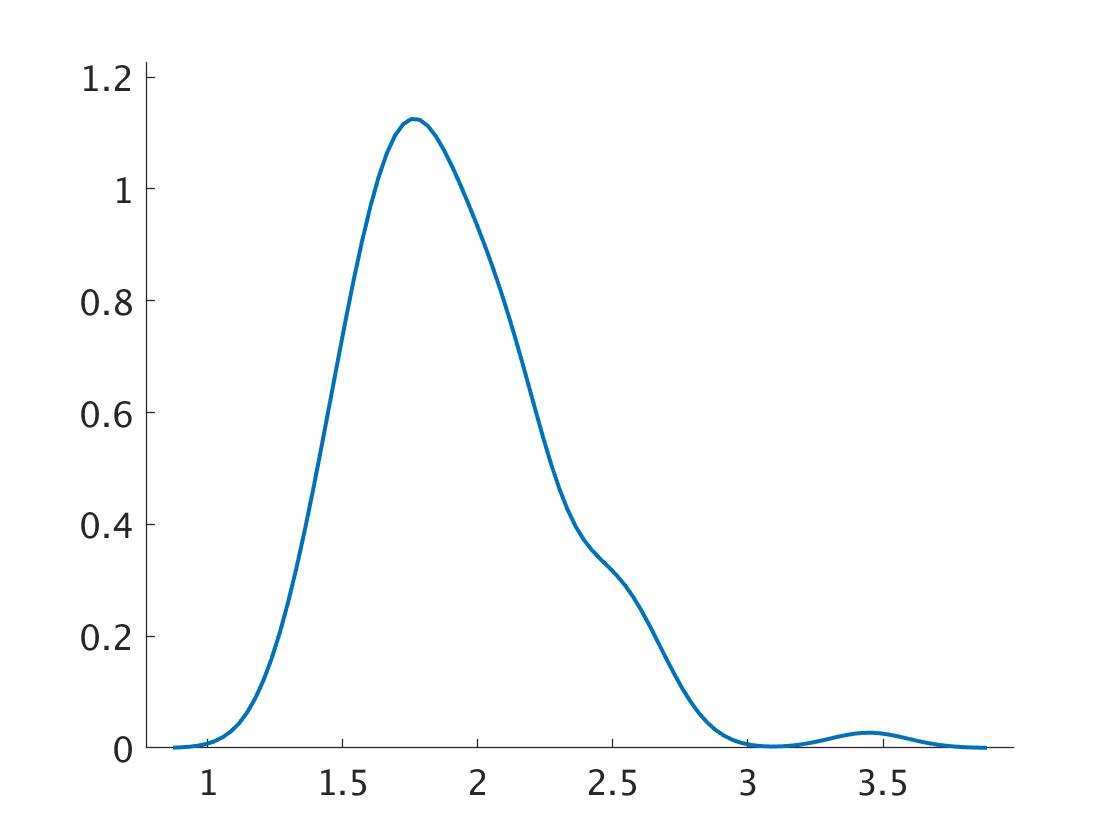}
    }
\caption{\footnotesize
 Smoothed empirical densities for the parameter estimates for the
$H_0$ persistence diagram, based on 100  simulations of non-concentric  circles. (a) $\alpha$, (b) $\theta_1/n$, where $n$ is the number of points in each persistence diagram (with the point at infinity, as usual, removed).}.
\label{fig:2farest}
\end{figure}

Figure \ref{fig:2farest} shows the empirical densities of the estimates of $\alpha$ and $\theta_1/n$, where, as before, $n$ is  the number of points in each diagram (without the `point at infinity').  As opposed to the Gaussian excursion case, very small values of $\alpha$ did not appear among the estimates.

\begin{figure}[h!]
    \centering
    \subfigure[]
    {
        \includegraphics[scale=.11]{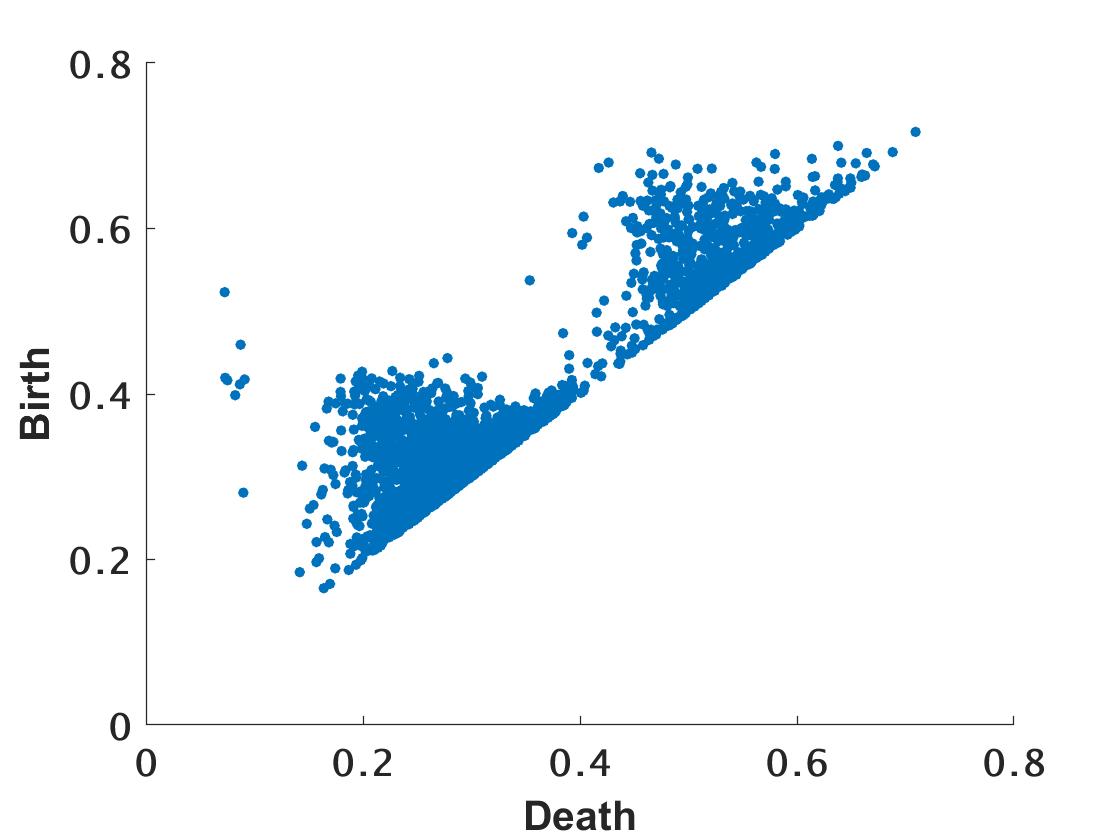}
    }
    \qquad
    \subfigure[]
    {
        \includegraphics[scale=.11]{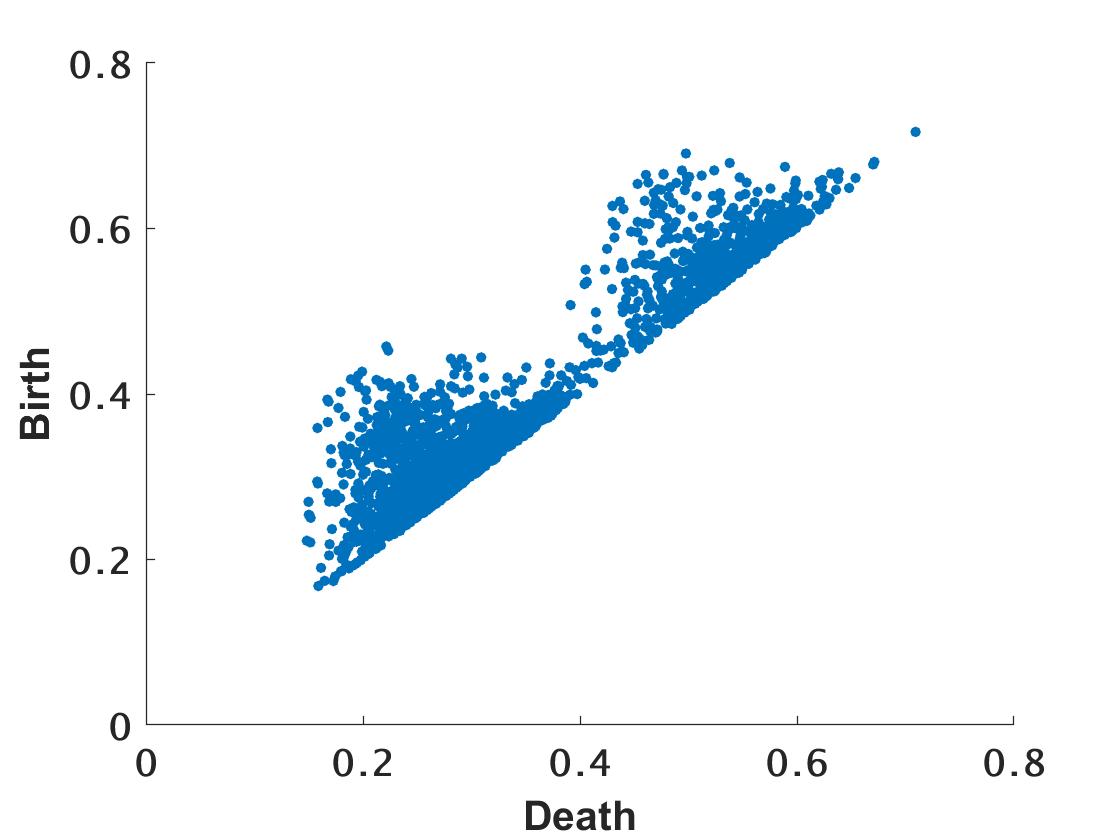}
    }
     \subfigure[]
    {
        \includegraphics[scale=.11]{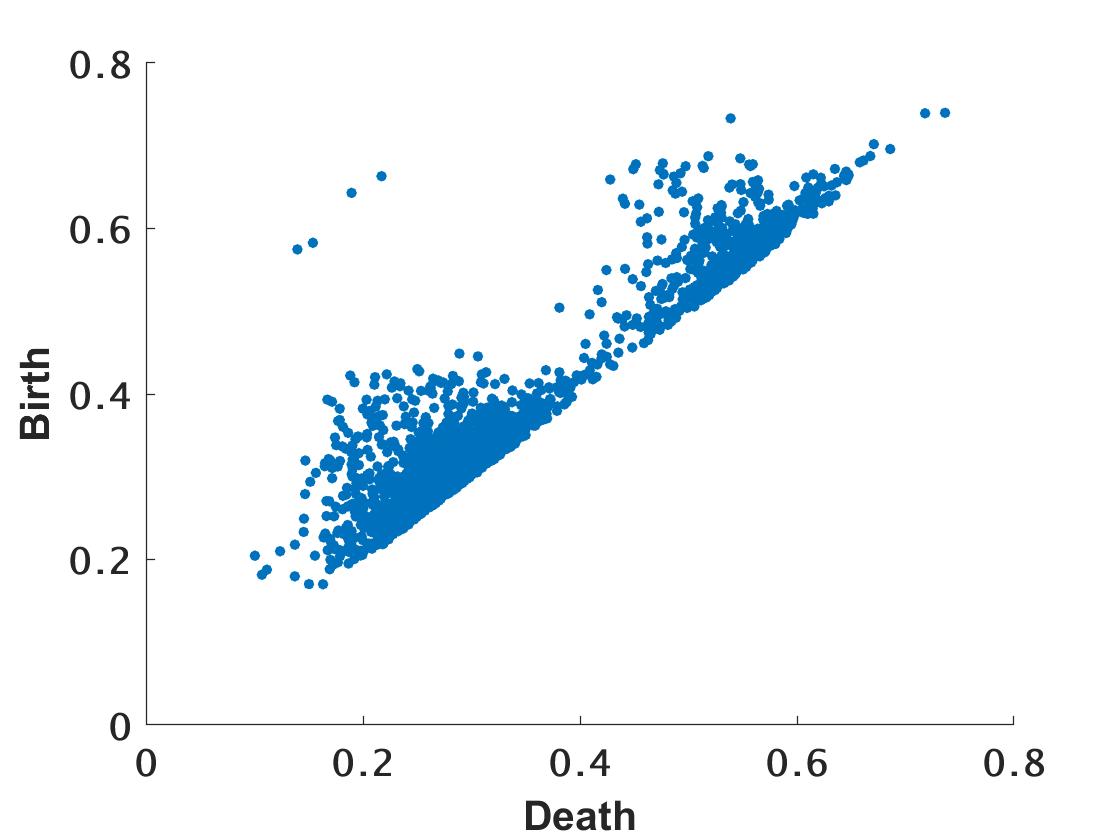}
    }  
  \caption{\footnotesize
   100 superimposed simulated  persistence diagrams, at  (a) 10, (b) 100, and (c)  500 steps into the MCMC routine. }
\label{fig:r3}
\end{figure}

In order to see if our model can reproduce the points clouds of Figure     \ref{fig:TwoCirclesPH:rob}, we took each of the simulations of `original' diagrams just described, and, on each one, with the parameters estimated as above, ran an MCMC simulation as in  Section \ref{sec:replications}. Figure \ref{fig:r3} summarises the results, superimposing the diagrams of  100  such simulations, at 10, 100, and 500 steps into the simulation.  Again, two facts are immediately noticeable. The first, as before, is that the general shape in Figure  \ref{fig:TwoCirclesPH:rob} is, in fact, well preserved in the simulations. The other, somewhat unexpected, is the quite rapid disappearance of the small cloud at the left in Figure     \ref{fig:TwoCirclesPH:rob}; viz.\ the cloud representing `signal' rather than `noise' in the original diagrams. Given that the main TDA motivation for replicating these diagrams is to detect signal points as outliers, this is a rather advantageous phenomenon for statistical testing via persistence diagrams.

\begin{figure}[h!]
    \centering
    \subfigure[]
    {
        \includegraphics[scale=.11]{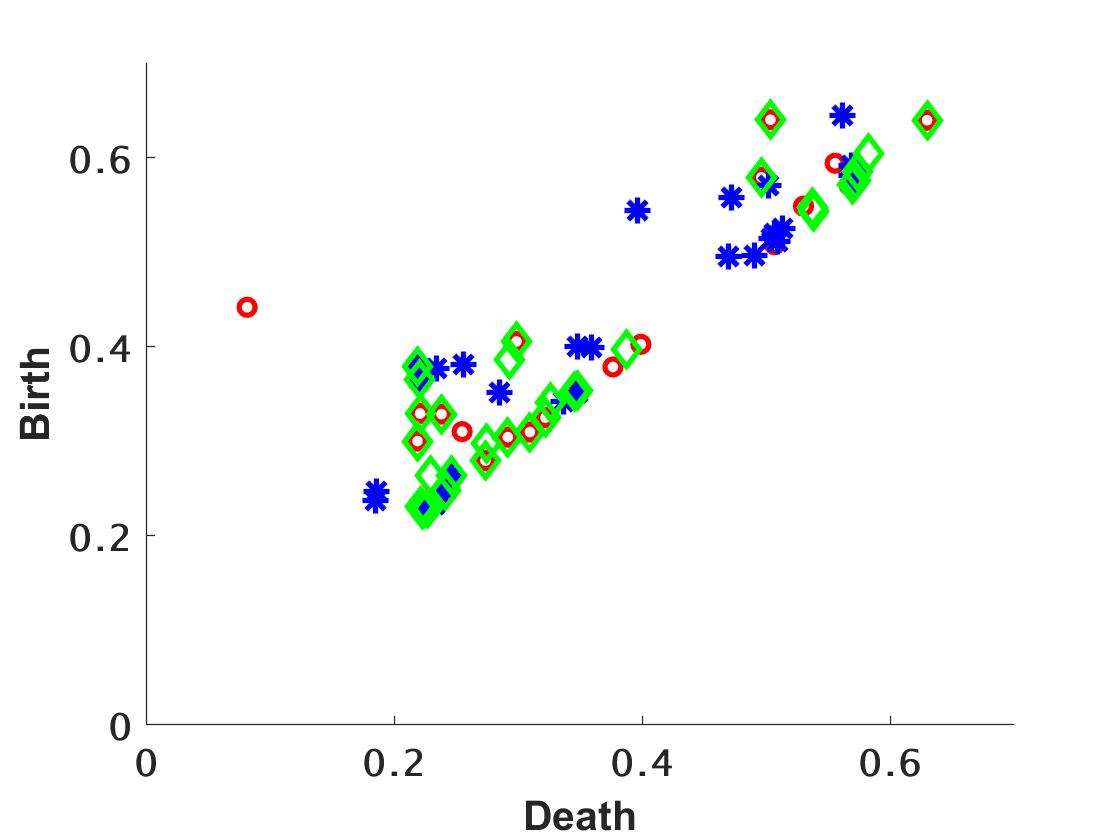}
    }
  \subfigure[]
    {
        \includegraphics[scale=.11]{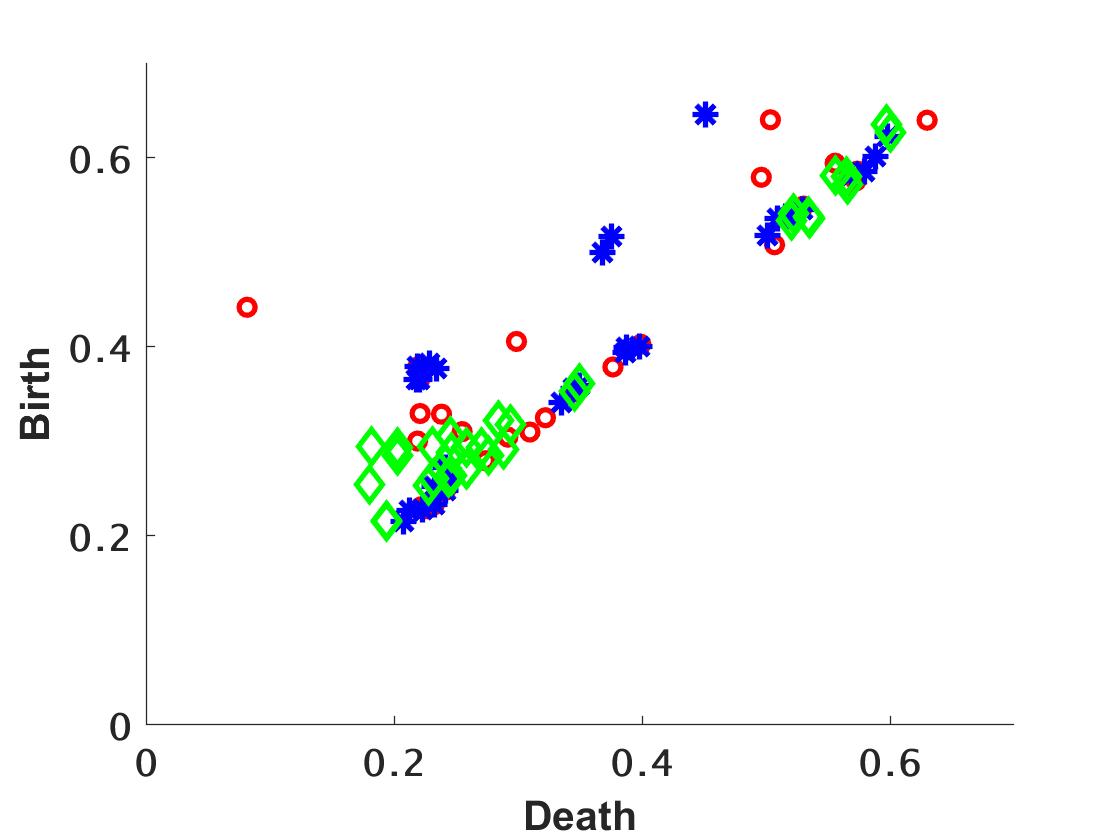}
    }
    \subfigure[]
    {
        \includegraphics[scale=.11]{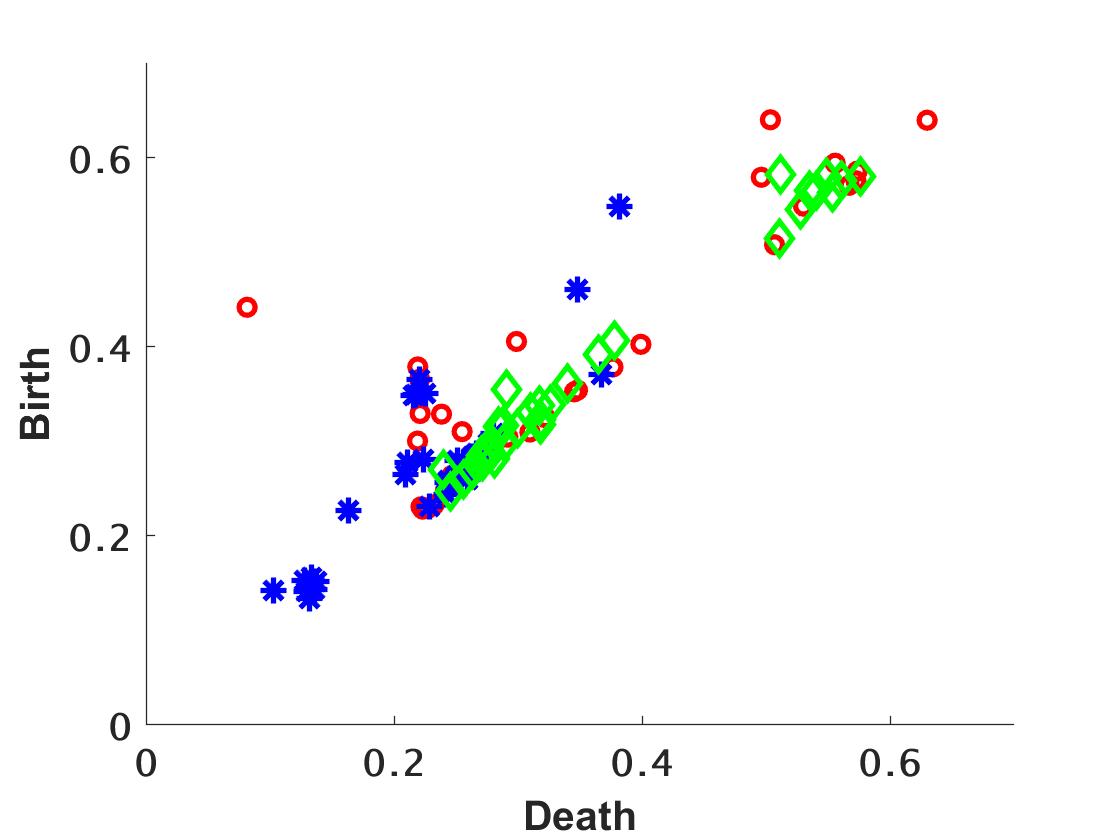}
    }

\caption{\footnotesize
    Comparison, for the  non-concentric circles example, of MCMC persistence diagram simulations based on the model from \cite{PNAS}  (blue stars),  with that based on the  model of the current paper (green diamonds), and with the original persistence diagram (red circles).  The numbers of MCMC steps are (a) 10,  (b) 100, and (c) 1,000. 
     }
     \label{2farRST}
\end{figure}

The final use that we shall make of this example is related to the  model developed in the earlier paper \cite{PNAS}. As noted in the Introduction, this model was not as efficient for modelling complex global structure in diagrams. To show that this is indeed the case here, consider the 
  persistence diagram of  Figure  \ref{fig:farcircle}(c), with  $n=31$ points. Fitting a model with $K=1$, yields estimate   $\alpha=1.6974$ and  $\theta _{1} = 55.0144$, or 
 $\theta _1/n =1.7747$; viz.\ when referring to the densities in Figure \ref{fig:2farest}, this is a reasonably typical case.
We also fitted the model from \cite{PNAS}, with the  Hamiltonian \eqref{eq:HamiltonianOld}, and then ran the MCMC model as usual, for each model. The results are presented in Figure 
 \ref{2farRST},  which  shows the results 10, 100 and 1,000 steps into the simulation. The red circles are the points of the original  persistence diagram, the green diamonds are the points of the simulated  diagram based on the model of this paper, and 
 the blue stars are  the simulated persistence diagram based on the  model from \cite{PNAS}.  The inability of the older model to capture the two clusters is clear even from this simple example, a fact also confirmed by larger scale simulations.  

\subsection{Concentric  circles}
\label{sec:concentric}

\begin{figure}[h]
\bc
 \subfigure[]
    {
    \includegraphics[scale=.22]{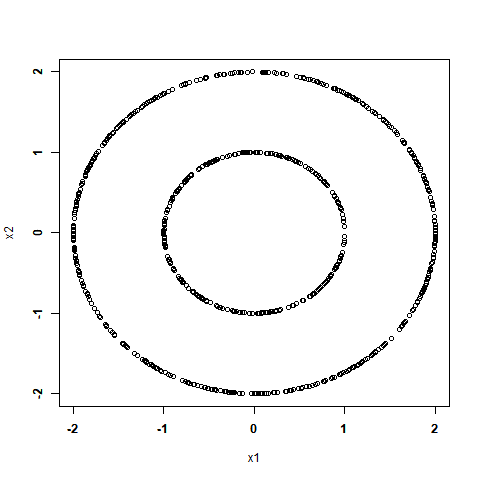}
    } 
\quad
 \subfigure[]
    {
     \includegraphics[scale=.12]{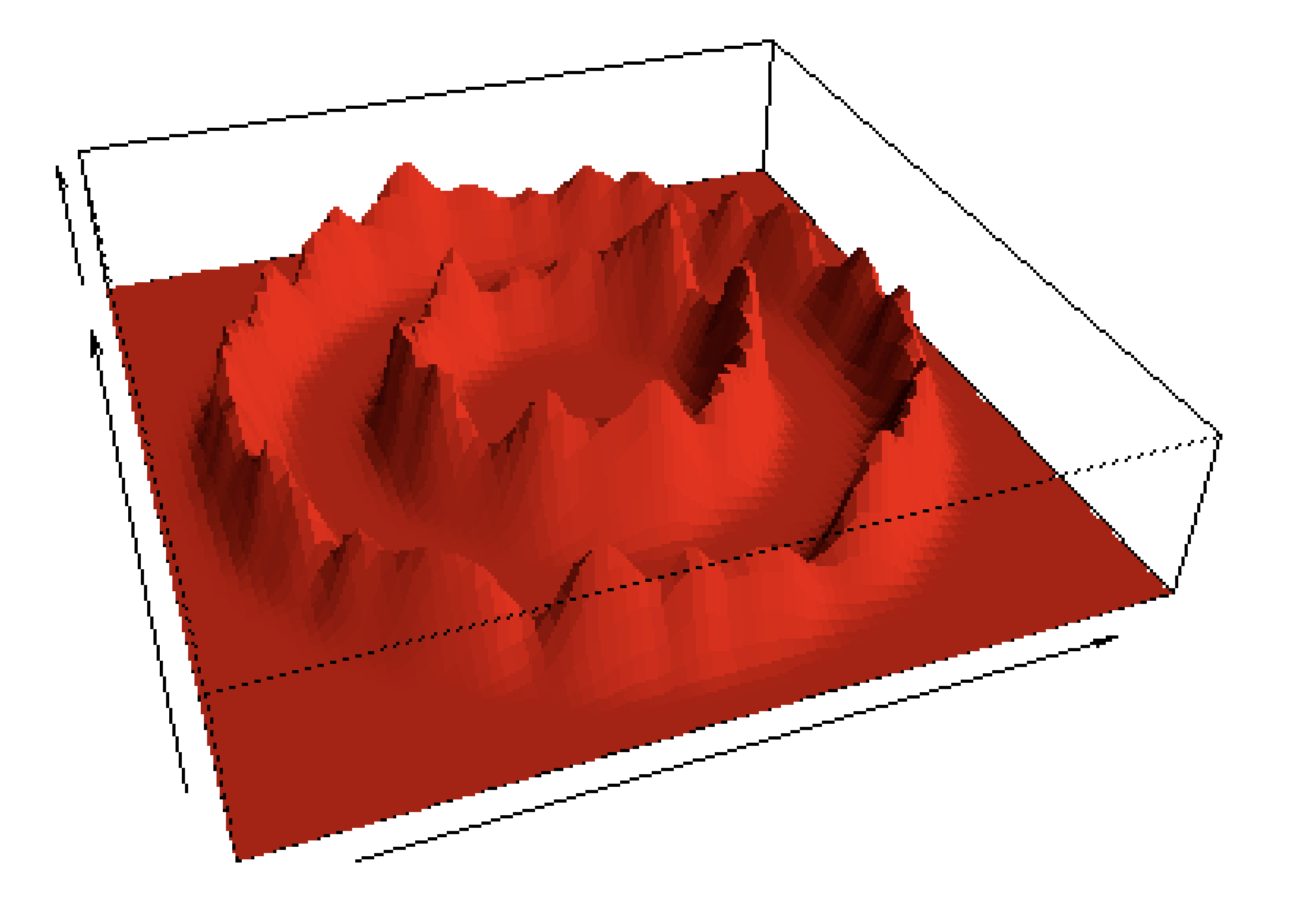} 
}
\quad 
 \subfigure[]
    {
    \includegraphics[scale=.22]{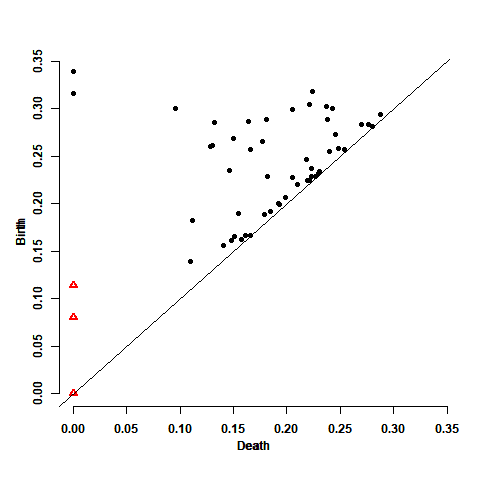}
}
\ec
\caption{\footnotesize
 A random sample from two circles: 500 points from the larger circle and 300 from the smaller one, with a kernel density estimate and the persistence diagram for its upper level sets. Black circles are $H_0$ persistence points, red triangles are  $H_1$ points. Birth times are on the vertical axis.}
\label{fig:twocircle}
\end{figure}

The next example was also investigated in \cite{PNAS}, and is based on a random sample of $n=800$ points from two concentric  circles of diameters 2 and 4,  with 500 points  chosen from the larger circle, and 300  from the smaller one.  Figure  \ref{fig:twocircle}(a) shows a typical  sample, followed in (b) by the corresponding kernel density estimate, with bandwidth $\eta =0.1$, and (c) shows the  persistence diagrams of the upper level set filtration for both   $H_0$  (circles) and  $H_1$ (triangles). The `point at infinity'  is still there,  and the two components and two circles, all of which have death height 0, are clearly represented. 

\begin{figure}[h!]
    \centering
        \subfigure[]
    {
        \includegraphics[scale=.08]{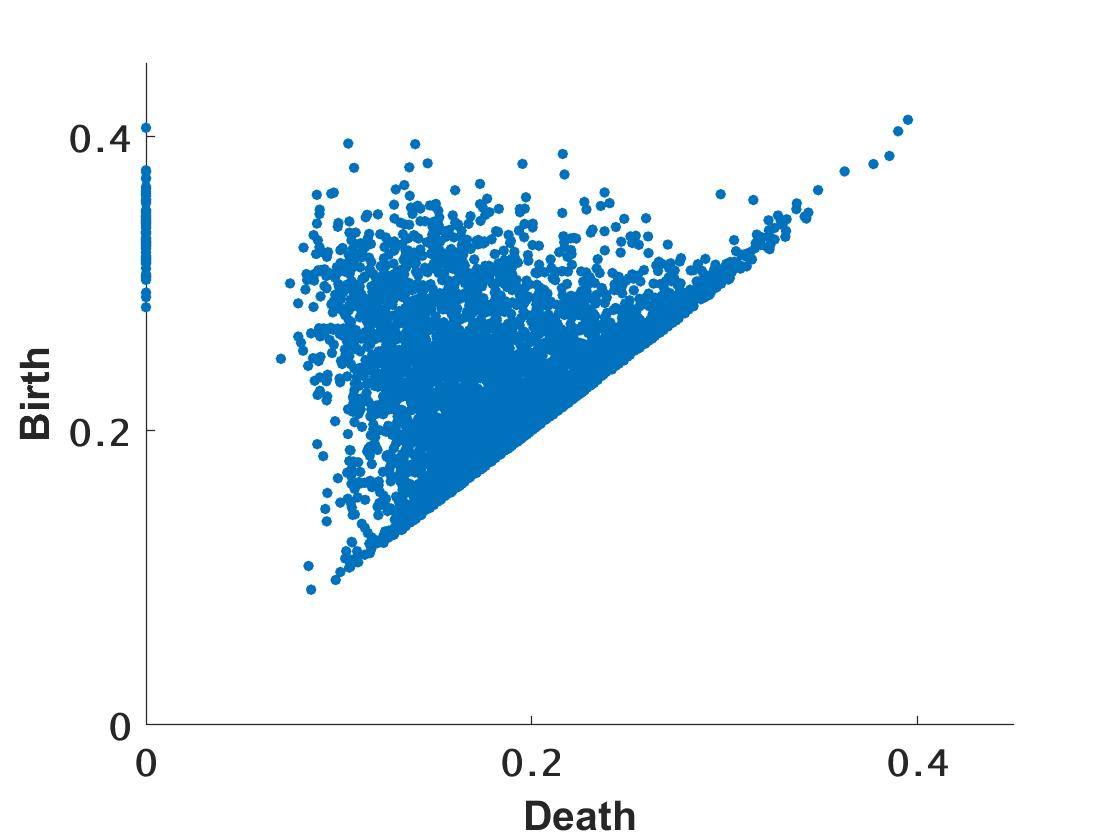}
    }
      \subfigure[]
    {
        \includegraphics[scale=.08]{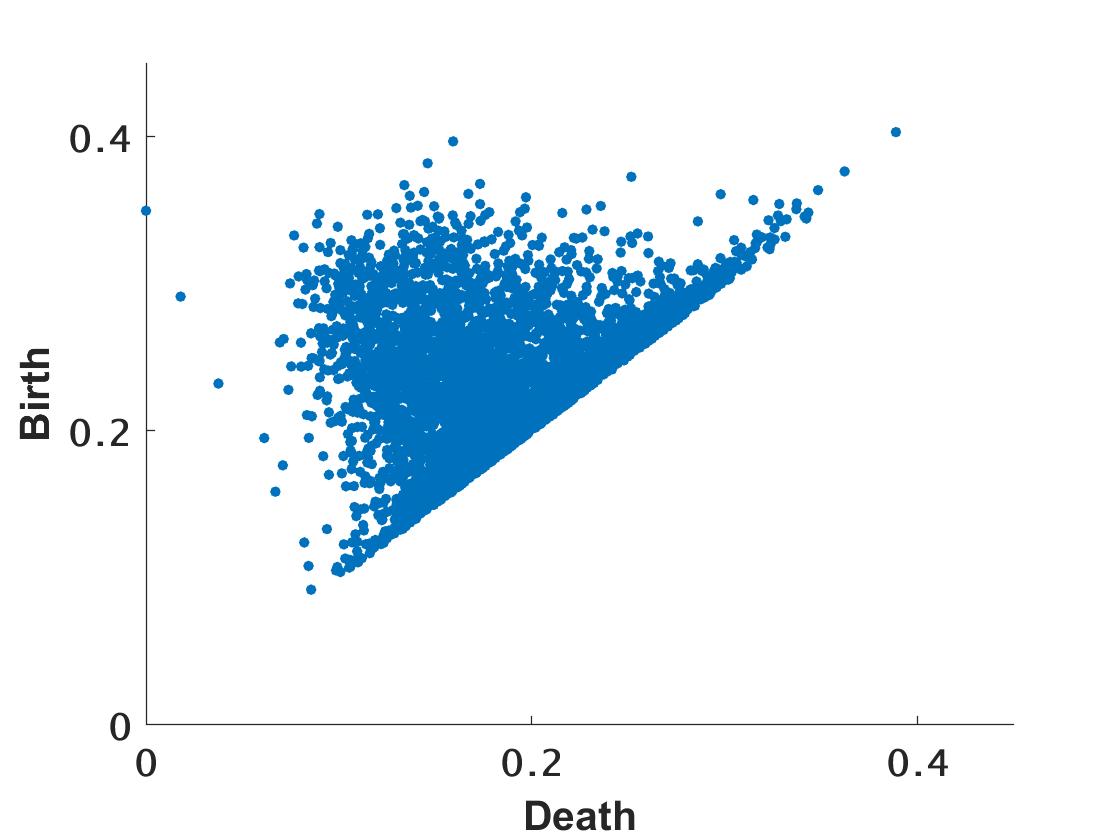}
    }
      \subfigure[]
    {
        \includegraphics[scale=.08]{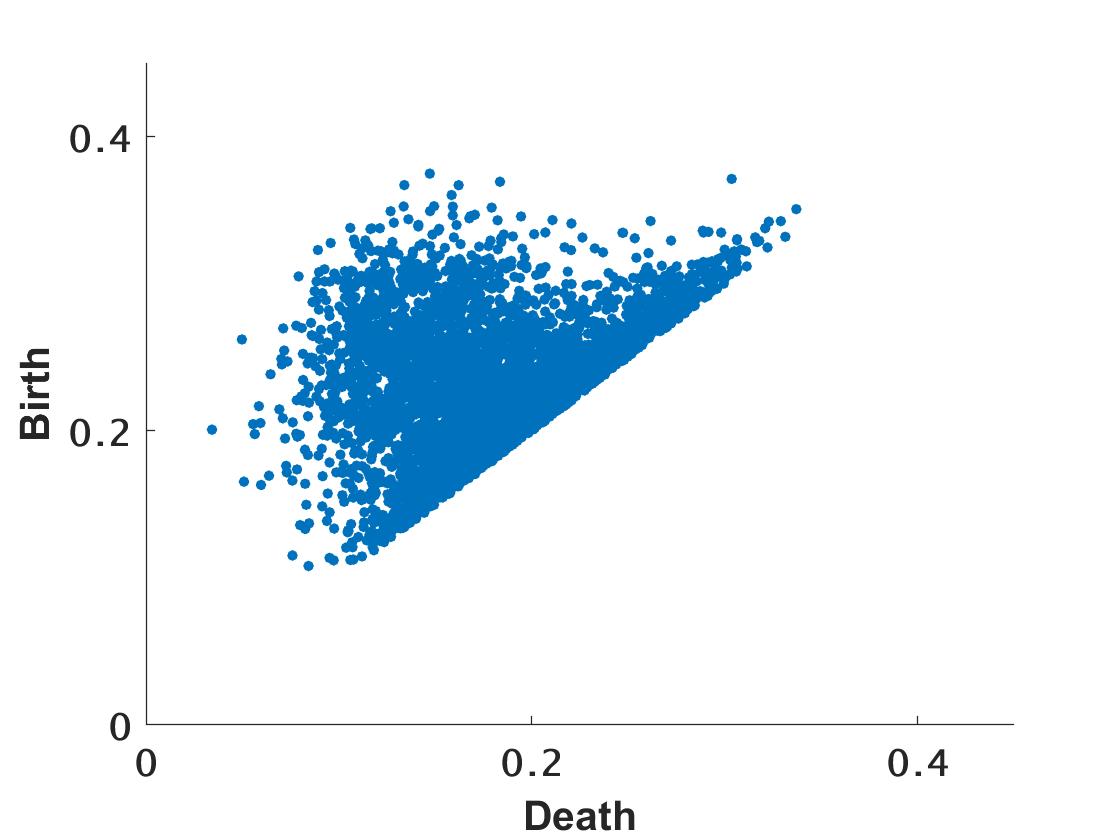}
    }
 \subfigure[]
    {
        \includegraphics[scale=.08]{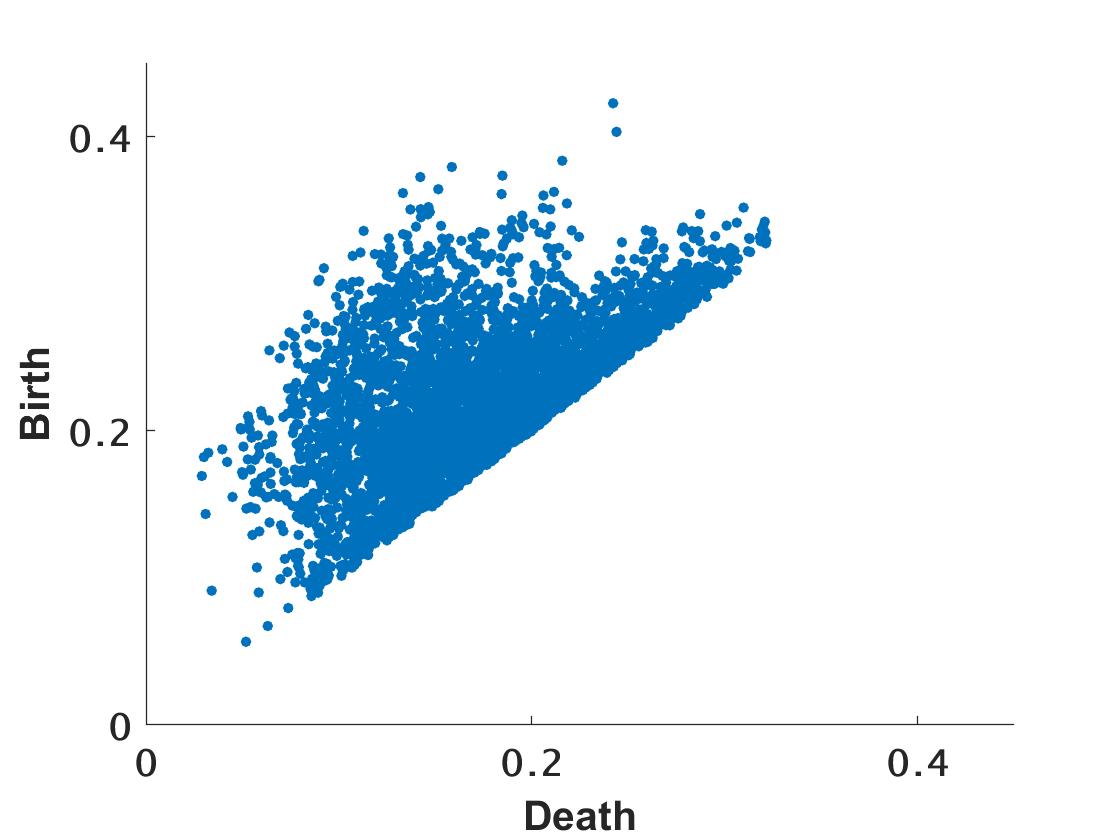}
    }
  \caption{\footnotesize  Superposition of 100 $H_0$ persistence diagrams for the concentric circles example, for the original diagrams  in (a), and then  100 MCMC simulations at (b)  10,  (c) 100, and (d)  500 steps into the simulation.}
 \label{fig:TwoCirclesTotal}
\end{figure}

As for the previous examples, our aim is to model the persistence diagram, in this case only for $H_0$, as there are not enough points in the $H_1$ diagram. Modulo the actual numerical values of the parameters, the procedure and results are essentially the same as for the previous example of non-concentric circles, and so  we shall not report on all the details of the analysis.
Indicative of this  similarity is Figure \ref{fig:TwoCirclesTotal}, which shows the superposition of 100 $H_0$ persistence diagrams, firstly from the original diagrams, based on the sampling procedure, and then from 100 MCMC simulations of these diagrams at  10, 100, and 500 steps into the simulation.

\begin{figure}[h!]
    \centering
    \subfigure[]
    {
        \includegraphics[scale=.09]{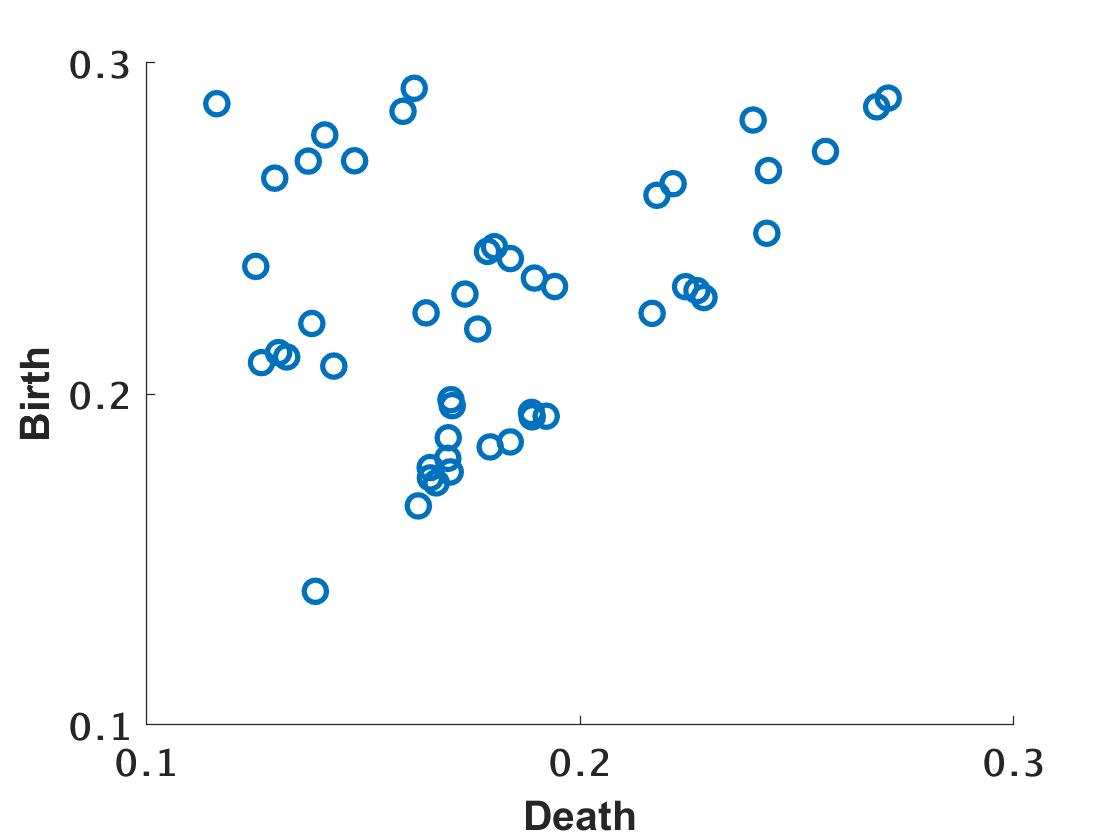}
    }
    \subfigure[]
    {
        \includegraphics[scale=.09]{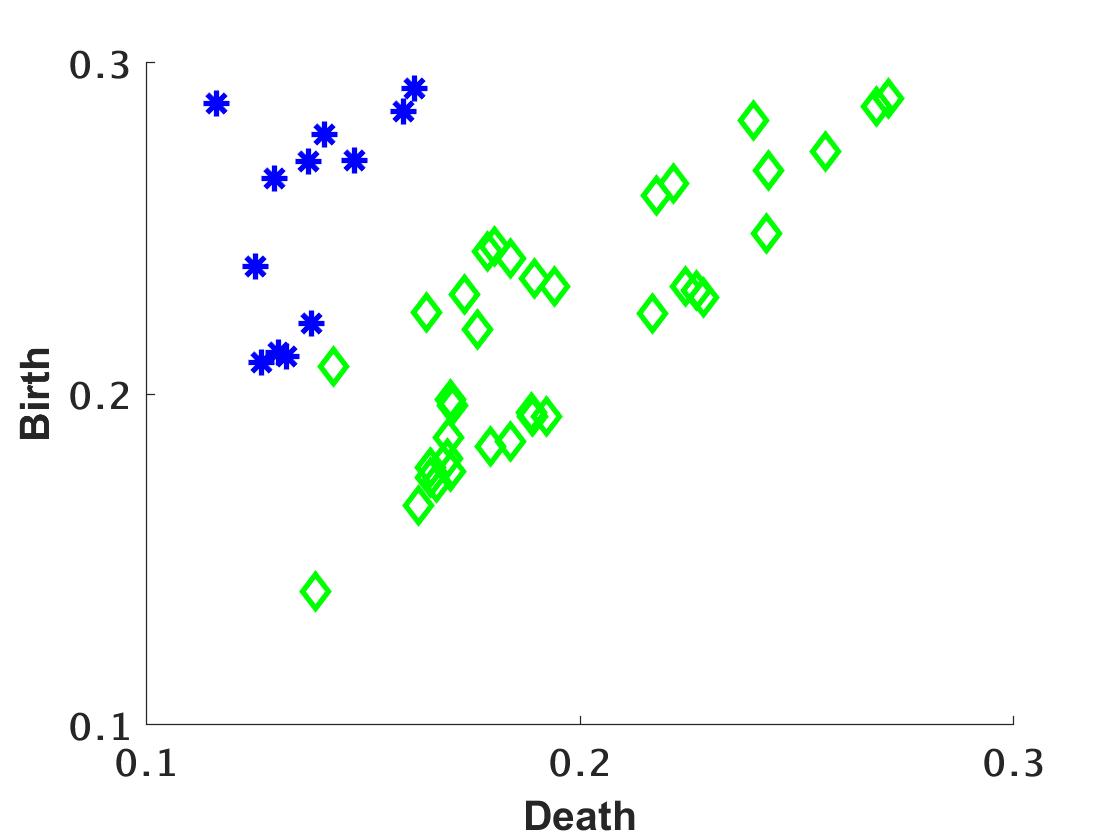}
    }
      \subfigure[]
    {
        \includegraphics[scale=.09]{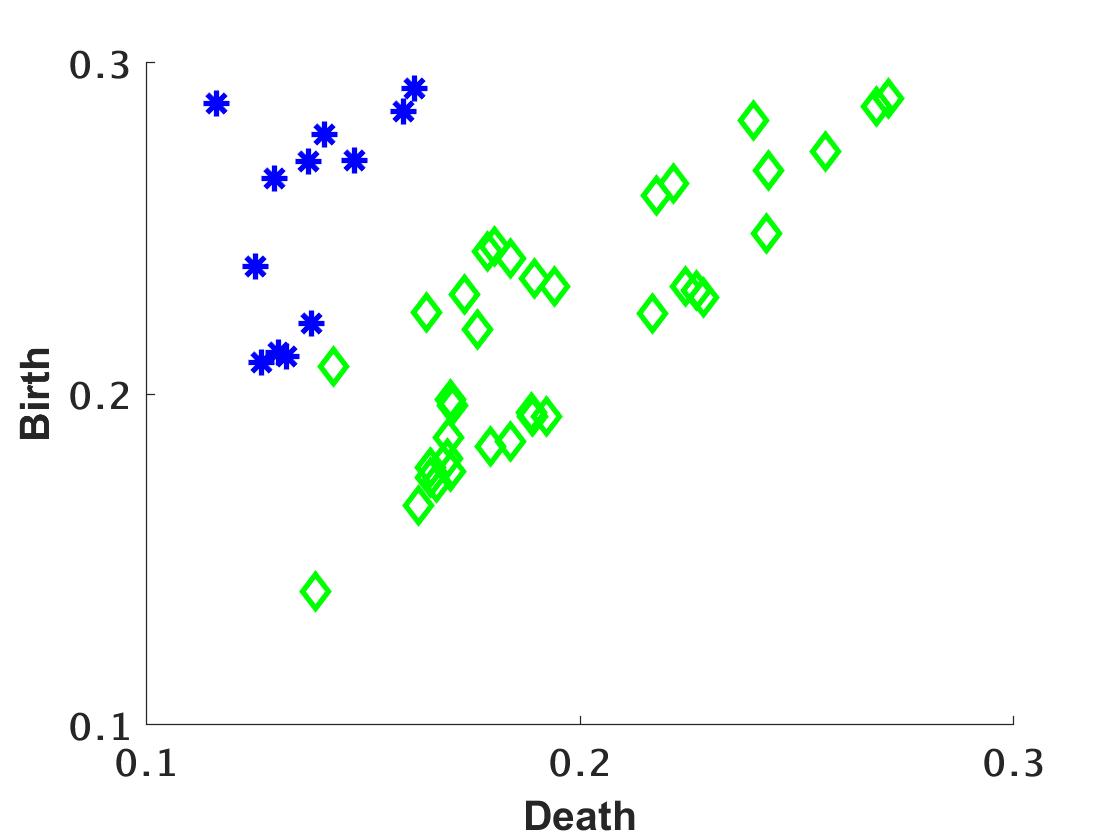}
    }
    \\
 \subfigure[]
    {
        \includegraphics[scale=.09]{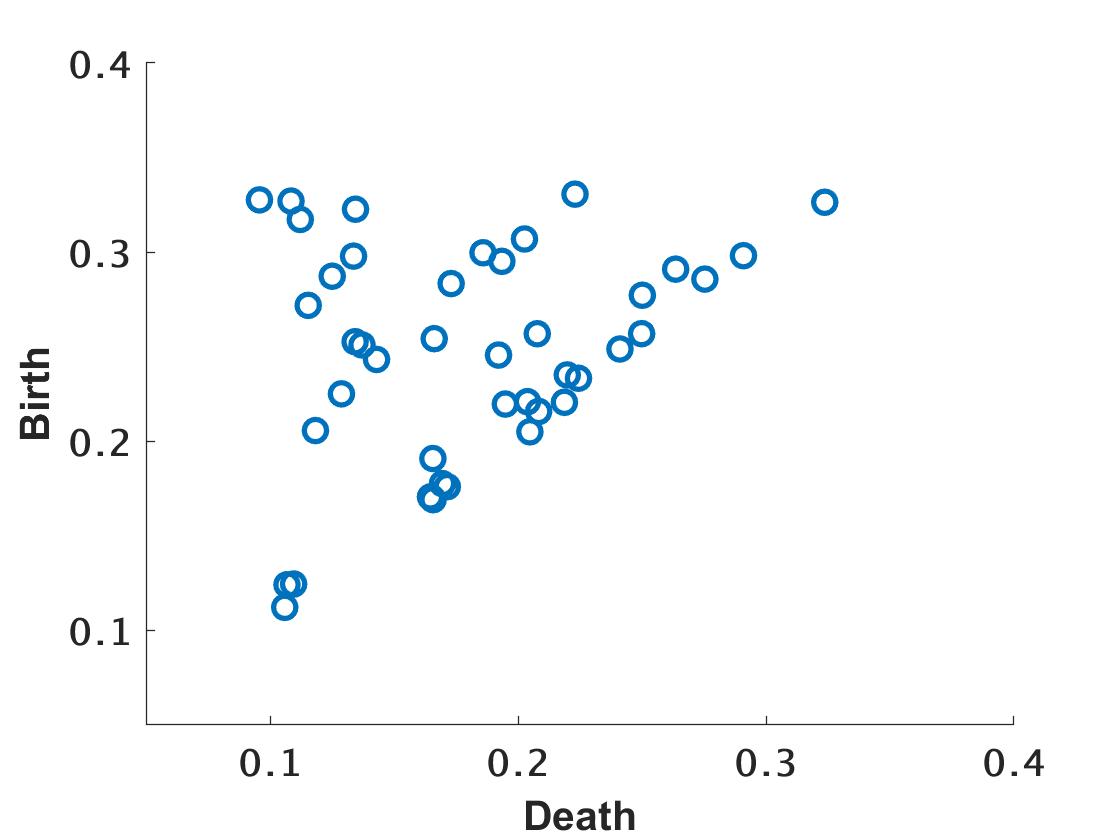}
    }
      \subfigure[]
    {
        \includegraphics[scale=.09]{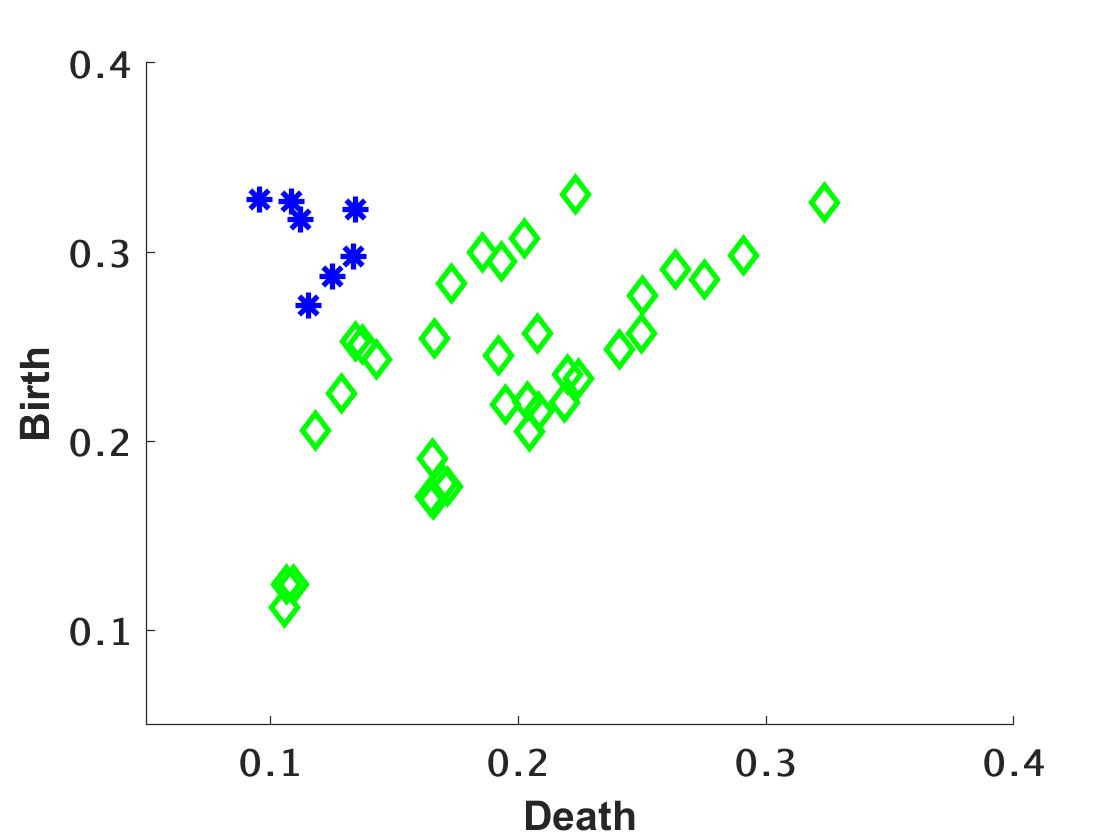}
    }
 \subfigure[]
    {
        \includegraphics[scale=.09]{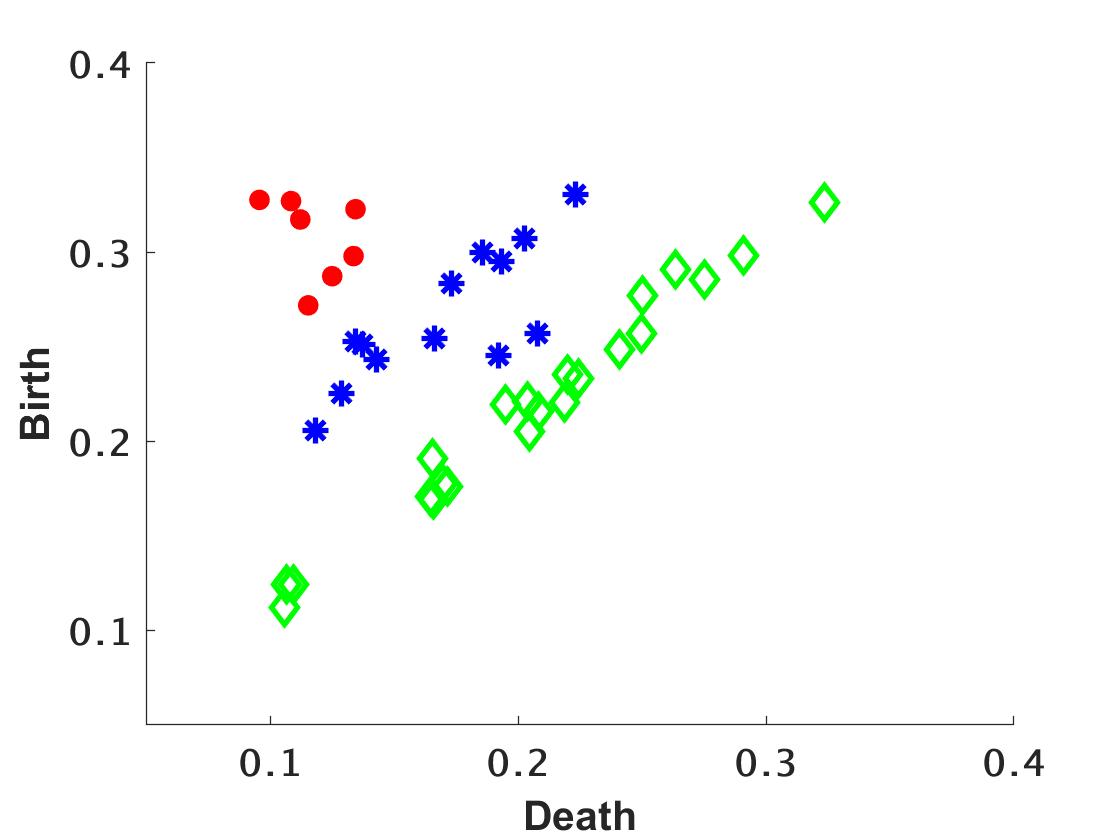}
    }
  \caption{\footnotesize  (a) and (d) Persistence diagrams for concentric circles. (b) and (e) Clustering according to the   75\% criterion. (c) and (f) Clustering according to the  90\% criterion.    In all cases  clusters must contain at least 5\% of the points.}
 \label{fig:TwoCirclesGrouping}
\end{figure}

Despite the fact that there was little different in the analysis of the concentric circles case, there is nevertheless a significant disparity between Figure \ref{fig:TwoCirclesTotal} and Figure \ref{fig:r3}, which shows corresponding results for the non-concentric circles case. In the previous case, the points in the diagram separated into two distinct groups, each one corresponding to the `noise' component for a different circle. There seems to be no corresponding separation in the present case, at least in the examples in Figure \ref{fig:TwoCirclesTotal}. However, consider Figure \ref{fig:TwoCirclesGrouping}.

Figure \ref{fig:TwoCirclesGrouping} (a) and (d) show two $H_0$ persistence diagrams, each coming from a MCMC simulation of diagrams from concentric circles. Apart from the two points to the far left of Figure \ref{fig:twocircle} (a), corresponding to the two connected components (circles) these are `typical' diagrams (although we have chosen them carefully to enhance the exposition).  At first inspection, there seems to be nothing remarkable about either of these diagrams, other than, perhaps, a tendency for them to exhibit `lines' of points parallel to the diagonal.

In order to determine whether this tendency was real or just a visual illusion, we defined, for each point $(d,b)$, a `lifetime' $\ell=b-d$, and used these lifetimes in order to group the points into clusters. 

 The procedure we adopted for this was hierarchical clustering, which groups  data over a variety of scales by creating a  multilevel hierarchy resulting in a cluster tree or dendrogram. 
More specifically, we  exploited the {\tt kmeans} clustering routine of Matlab, using the  cosine distance as the distance metric.  For the roots of the trees we took those points that were at distance greater than the 75-th (or 90-th) percentile of all distances in a specific diagram. All clusters with less than  5\% (or 2.5\%)  of the data points were ignored.  All told, with two choices for the roots and two for cluster sizes, we have the had four choices for clustering, examples of which are shown in  
Figure \ref{fig:TwoCirclesGrouping}. 
 (See, for example, \cite{Hastie} for more information on clustering algorithms in general, and the Matlab site 
 {\tt www.mathworks.com/help/stats/kmeans.html} for details on the specifics.)

In  Figure \ref{fig:TwoCirclesGrouping},   points belonging to same cluster are marked with the same colour. Either two or three large clusters are found, and a little thought shows from where they come. 

Returning to the empirical density of Figure \ref{fig:twocircle} (b), note that while the heights of the two jagged rings there are similar (in fact we choose the sample sizes so that this would be the case) they are not the same. The peaks of the inner ring tend, on average, to be a little higher, and the `lifetime height', or distance between them and the nearby local minima, a little longer than the corresponding numbers for the outer ring. This is what leads to the two clusters in, for example, Figure \ref{fig:TwoCirclesGrouping} (b) and  (e). Note that in both of these cases the lower cluster includes a large number of points close to the diagonal. In (f), using a finer clustering, the points close to the diagonal are recognised as a separate cluster, corresponding to what one might call `fine scale' topological noise.

\begin{figure}[h!]
    \centering
    \subfigure[]
    {
        \includegraphics[scale=.08]{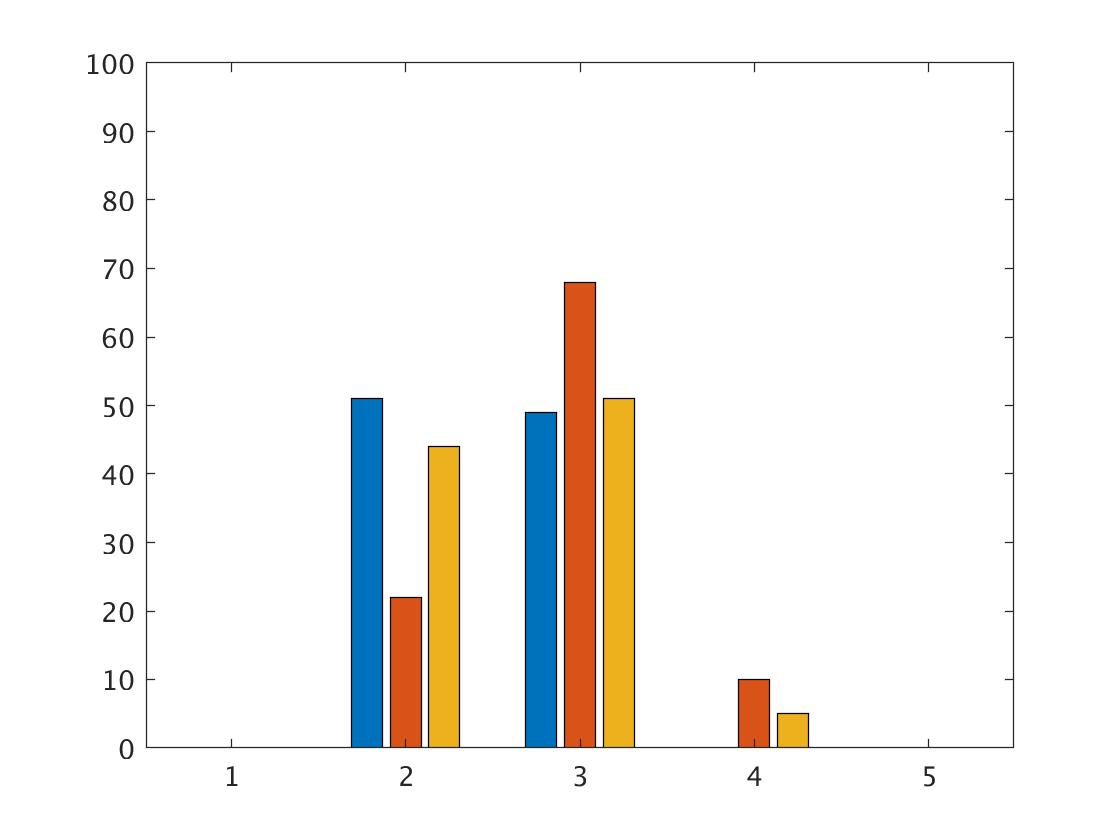}
    }
     \subfigure[]
    {
        \includegraphics[scale=.08]{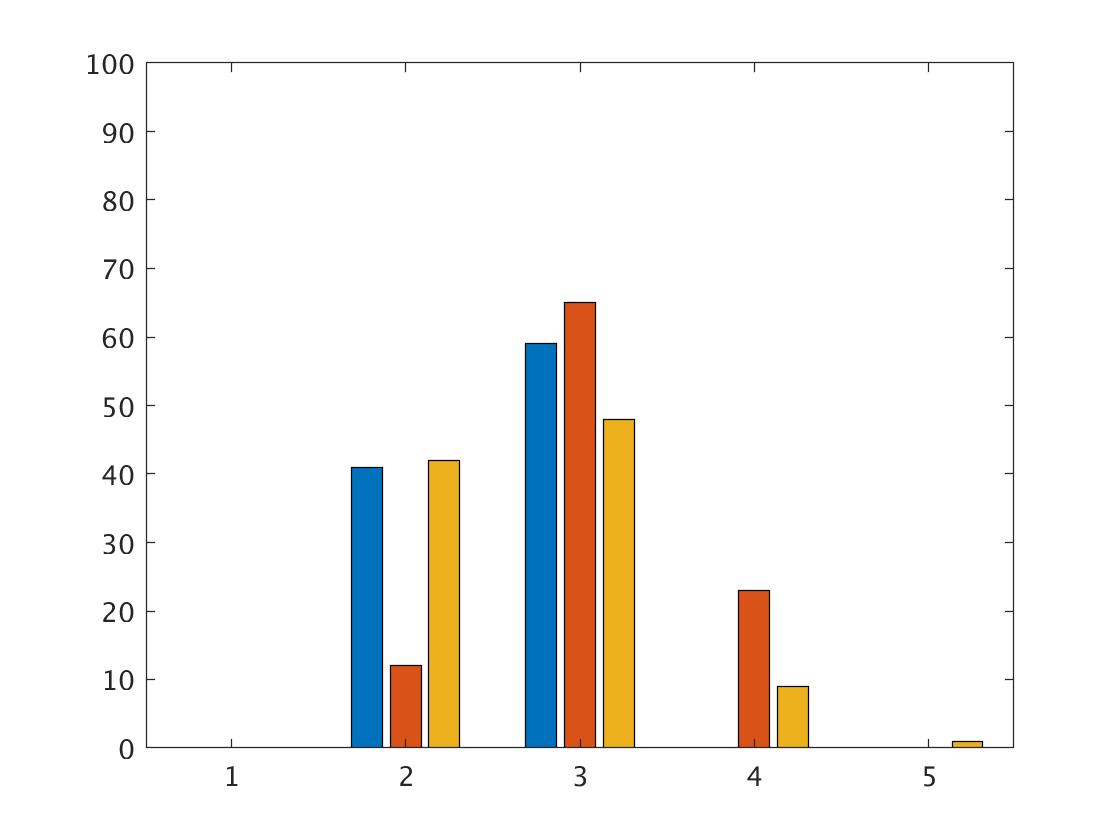}
    }
  \subfigure[]
    {
        \includegraphics[scale=.08]{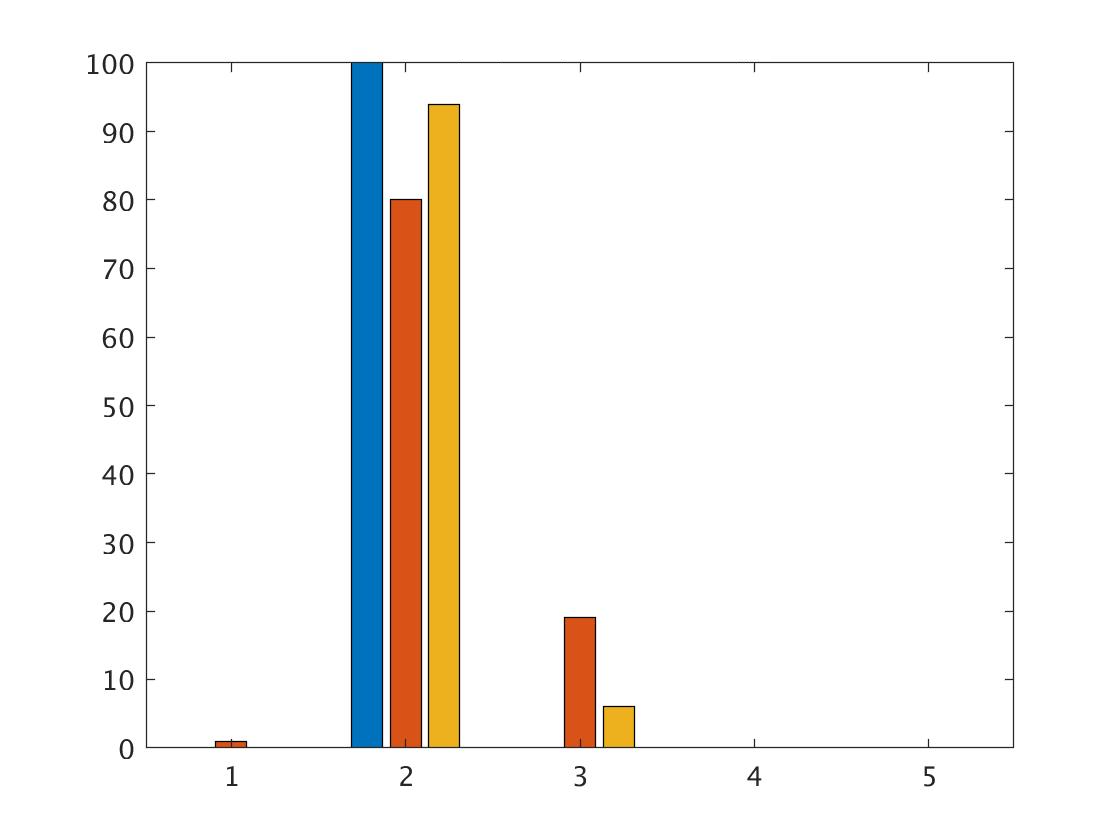}
    }
   \subfigure[]
    {
        \includegraphics[scale=.08]{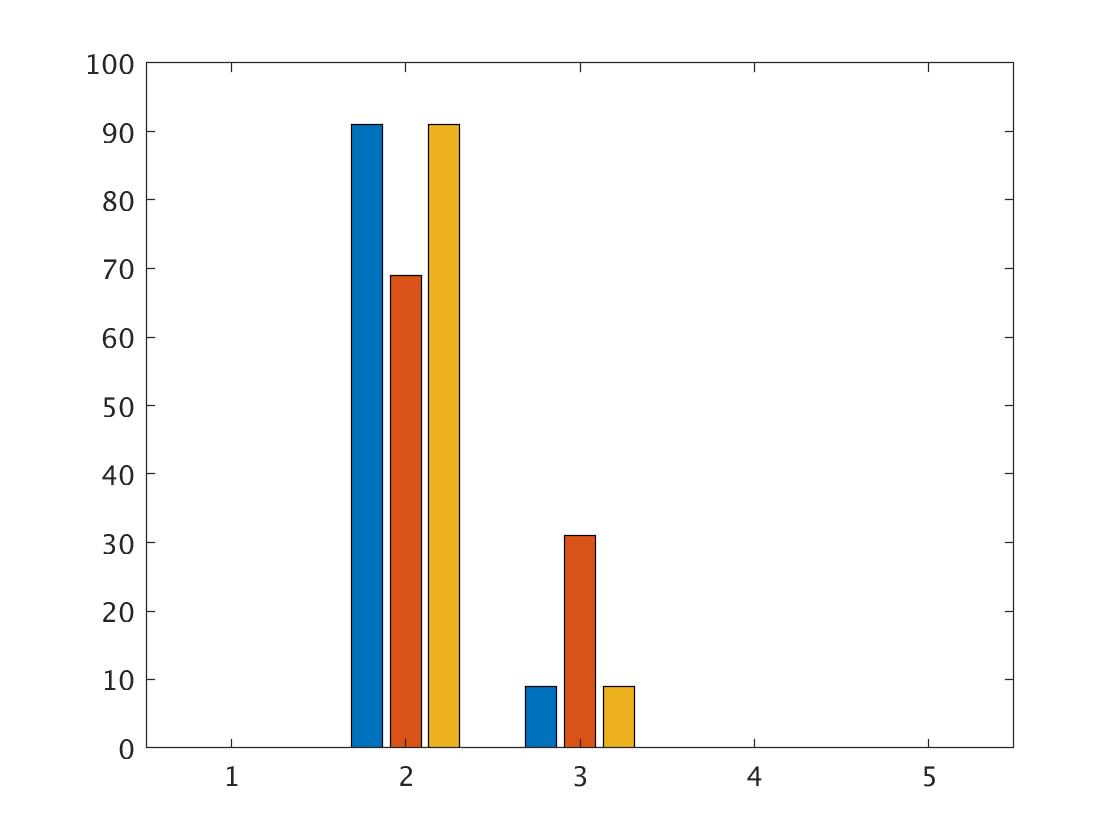}
    }
\caption{\footnotesize
    Distribution of numbers of clusters in the concentric circles persistence diagrams. The frequencies for the observed  diagrams are shown in blue (leftmost in each set of three bars)  with orange (rightmost) for  results of an MCMC simulation of the model of this paper at step 1,000, while  red (middle) corresponds to the model of \cite{PNAS}. In (a) and (b) the initial points of the cluster are taken at the 75-th percentile, and in (c) and (d) at the 90-th. In (a) and (c) only clusters with more than 5\% of the points are kept, and in (b) and (d), 2.5\%. See  text for further details.}
    \label{circlesCompare}
\end{figure}

Figure \ref{circlesCompare} shows that this phenomenon is general, and not restricted to the two persistence diagrams of Figure \ref{fig:TwoCirclesGrouping}. It shows the frequency distribution of  numbers of clusters over different clustering parameters, for 100 `original' persistence diagrams, and for 100 MCMC simulations at step 1,000, one for each original diagram. In addition, we have included the corresponding results for the MCMC simulations based on the model of \cite{PNAS}.

 The success of the current procedure in matching the true cluster number distribution is obvious, as is the improvement gained with the newer Gibbs measure incorporating the empirical density $\bar f^G$.

Given the clarity of the clustering phenomenon seen in Figure \ref{circlesCompare}, it is natural to ask why this was not seen in the superpositions of  Figure \ref{fig:TwoCirclesTotal}. The reason is that the `parallel lines to the diagonal' seen in Figure  \ref{circlesCompare} occur at different distances from the diagonal in different realisations, and so the structure is lost, or blurred, under the superpositions of Figure \ref{fig:TwoCirclesTotal}. 

A take home message from this is the need to study the probabilistic structure of persistence diagrams at different scales, which is precisely what the pseudolikelihood  \eqref{eq:finalpseudo} allows for, by incorporating both local and global behaviour.

%
%
%
%
%
%

\subsection{The two dimensional sphere}
\begin{figure}[h!]
    \centering
    \subfigure[]
    {
        \includegraphics[scale=0.08]{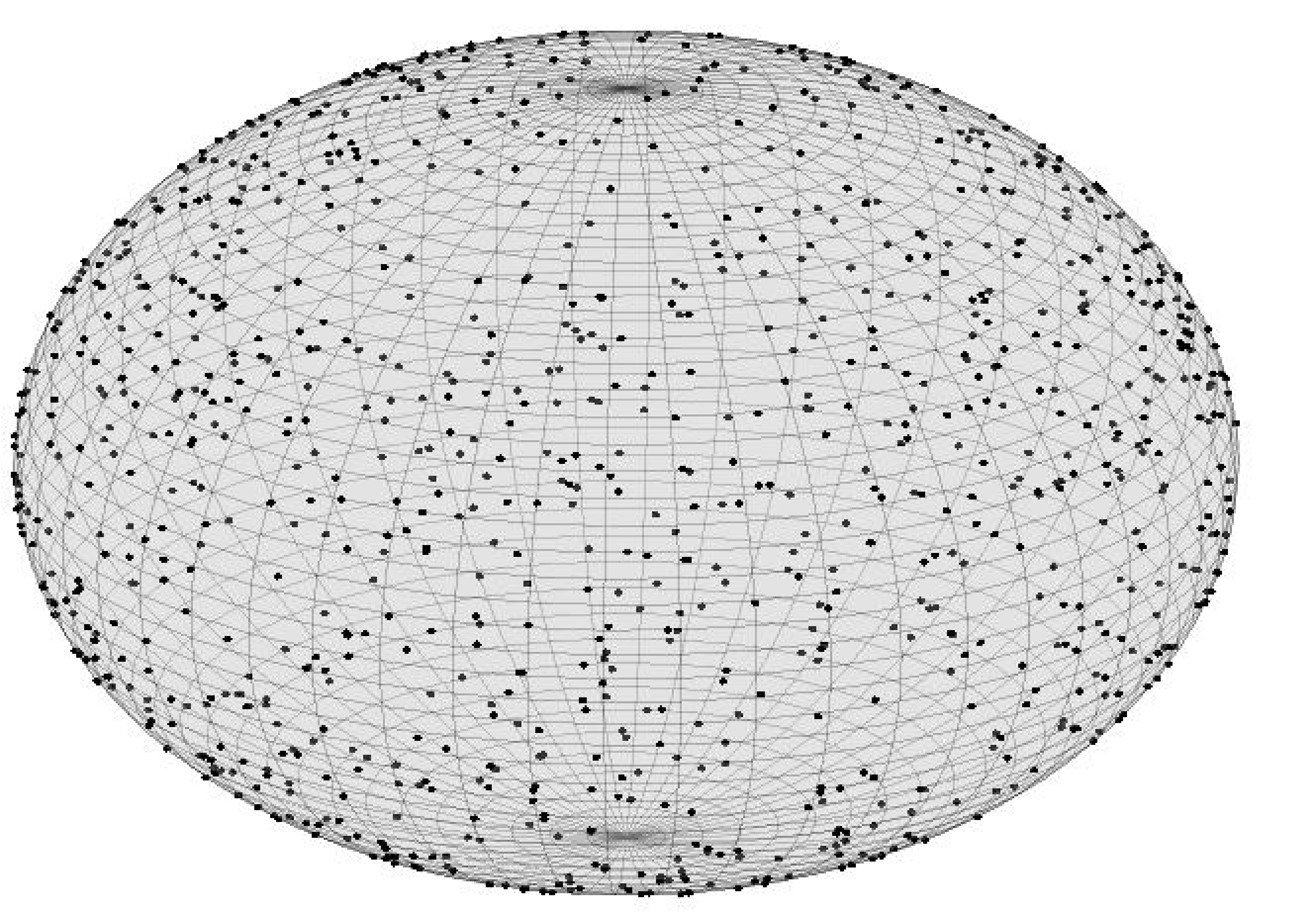}\ \  \ \
    }
       \subfigure[]
    {
        \includegraphics[scale=0.11]{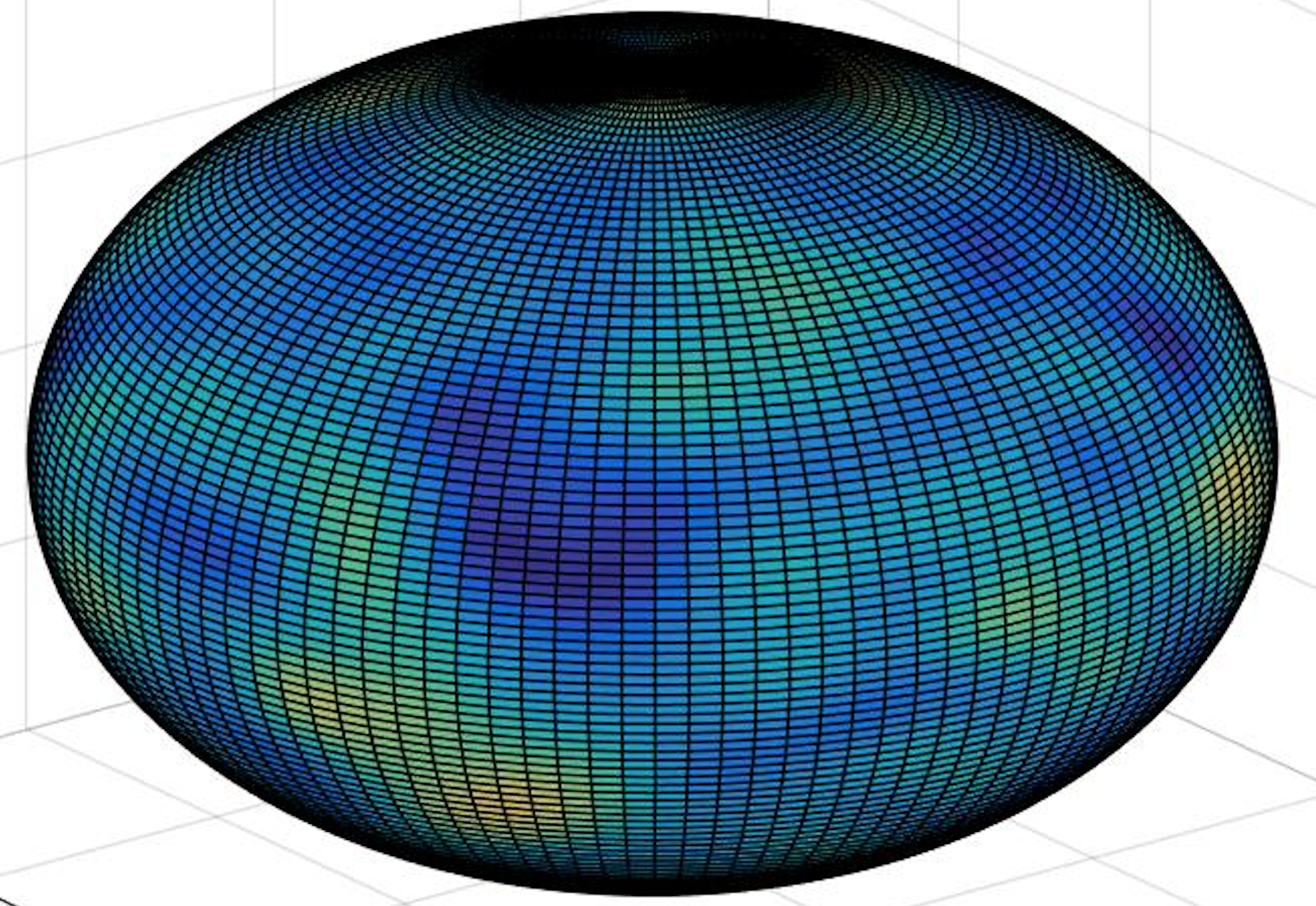} \ \ \ \
    }
       \subfigure[]
    {
        \includegraphics[scale=0.1]{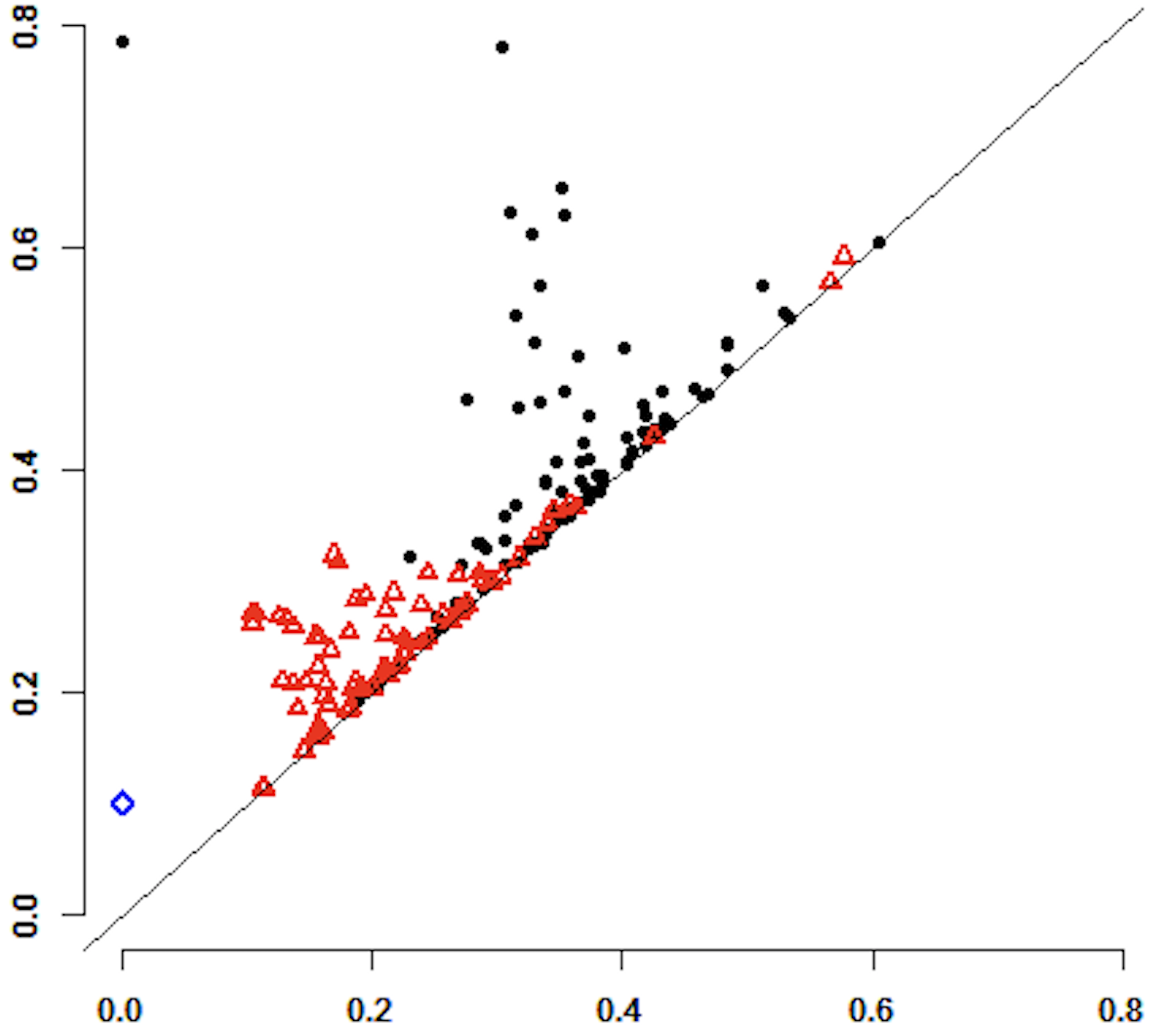}
        
    }
\caption{\footnotesize
 (a) Points sampled from a unit sphere. (b) The corresponding kernel density estimator, shown, for visual clarity, at  only a few quantized levels.  (c)  The corresponding persistence diagram for the upper level sets of the kernel density estimate on the full sphere. Black circles are  $H_0$ persistence points,  red triangles are  $H_1$ points, and the single blue diamond is  the $H_2$ persistence point. Birth times are on the vertical axis.}
\label{fig:sphere}
\end{figure}

Our final example is also described in \cite{agamiadler}, where it is treated according the model of Appendix \ref{sec:appendixglobalshape}. It is  based on  a random sample of $n=1,000$ points from the uniform distribution on the unit  sphere $S^2$ in $\real^3$, and then smoothing the data (in $\real^3$) with a kernel density estimator of bandwidth of $\eta=0.1$. These two steps are shown as Panels (a) and (b) of Figure \ref{fig:sphere}.

Most of the analysis here is parallel to that of the earlier examples, with the main novelties being that the empirical density is  defined over $\real^3$ rather than $\real^2$, and that the $H_1$ persistence diagram contains the same order of magnitude of points as found in the $H_0$ diagram. (74 and 110 points, respectively, for the example of Figure \ref{fig:sphere}  (c).)
 Consequently   there are more than  enough points to fit a spatial model to both diagrams.
 
 A certain (reflected) symmetry is to be expected between the $H_0$ and $H_1$ persistence diagrams. At the level of Betti numbers, 
 if set $\beta_i(u)$ to be the $i$-th Betti number of the excursion set of the empirical density at level $u$,  we must have, for $u>0$, 
 \beqq
 \beta_0(u) \ = 2 + \beta_1(u), 
 \eeqq
 since the difference  $\beta_0(u) - \beta_1(u)$ will always be 2, the Euler characteristic of $S^2$. (For all 
  $u>0$, we must have $\beta_2(u)=0.$)
   
\begin{figure}[h!]
    \centering
        \subfigure[]
    {
        \includegraphics[scale=.08]{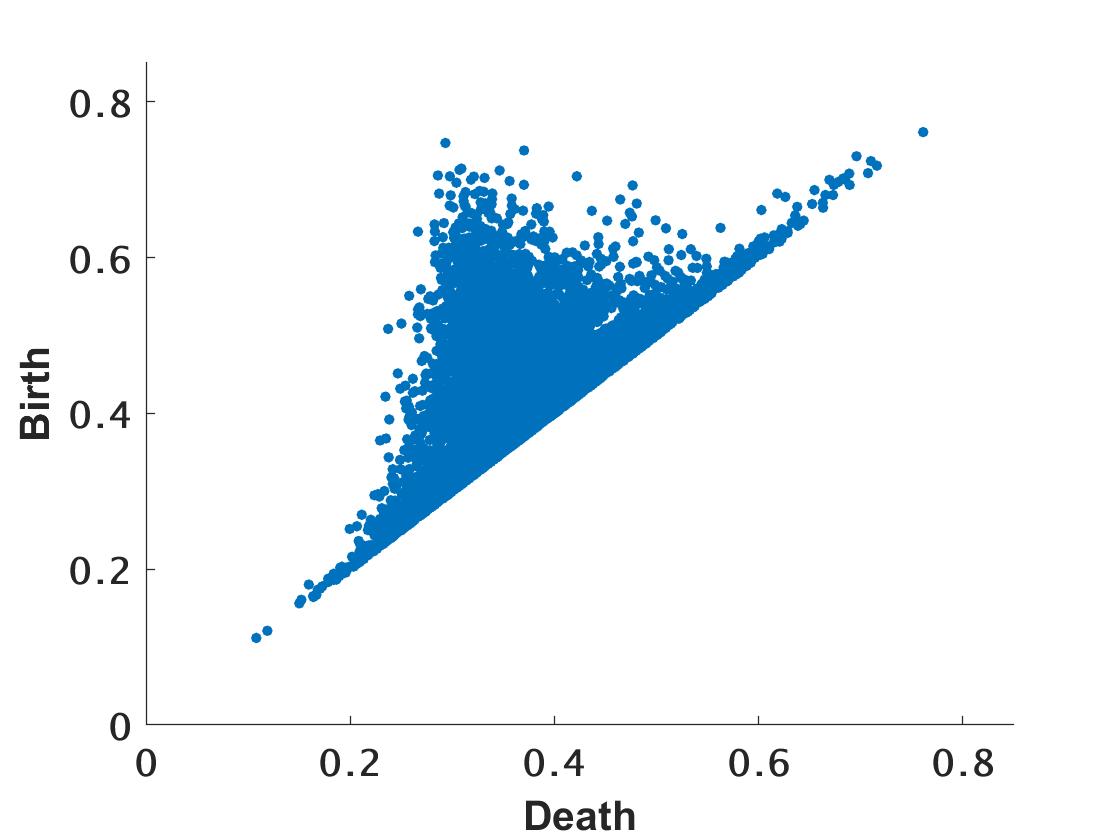}
    }
      \subfigure[]
    {
        \includegraphics[scale=.08]{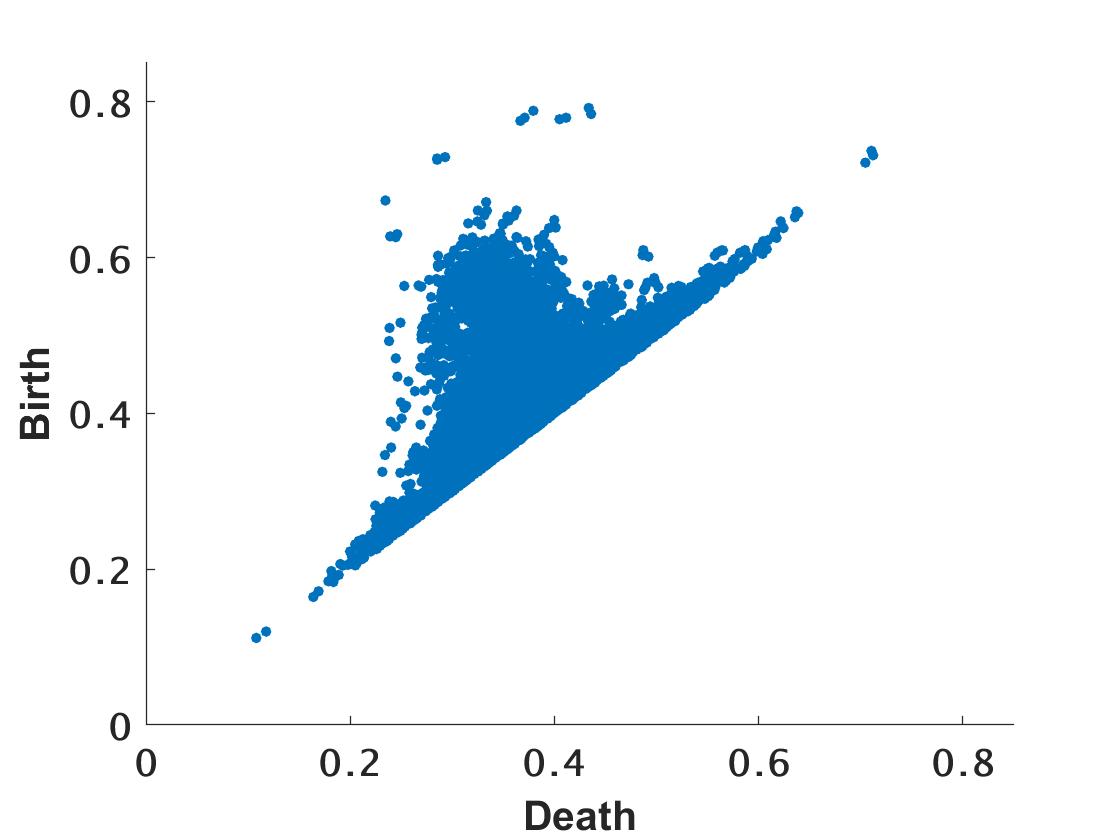}
    }
     \subfigure[] 
 {
        \includegraphics[scale=.08]{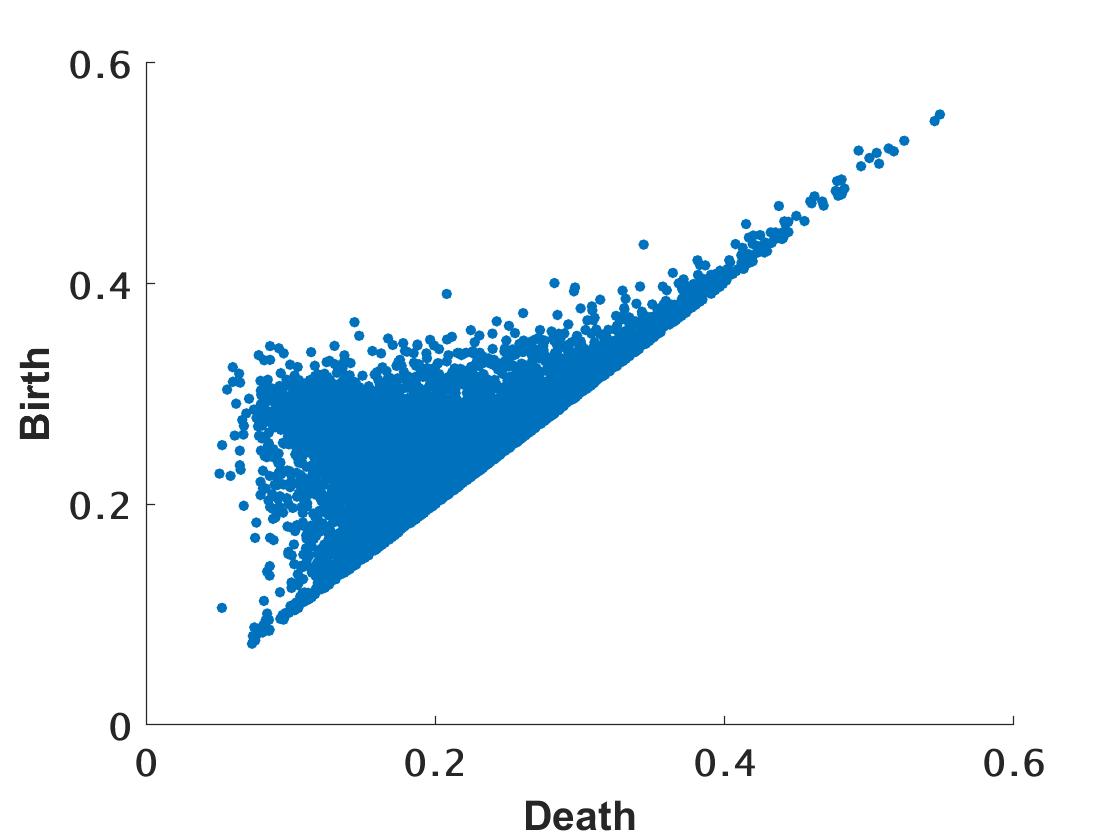}
    }
      \subfigure[]
    {
        \includegraphics[scale=.08]{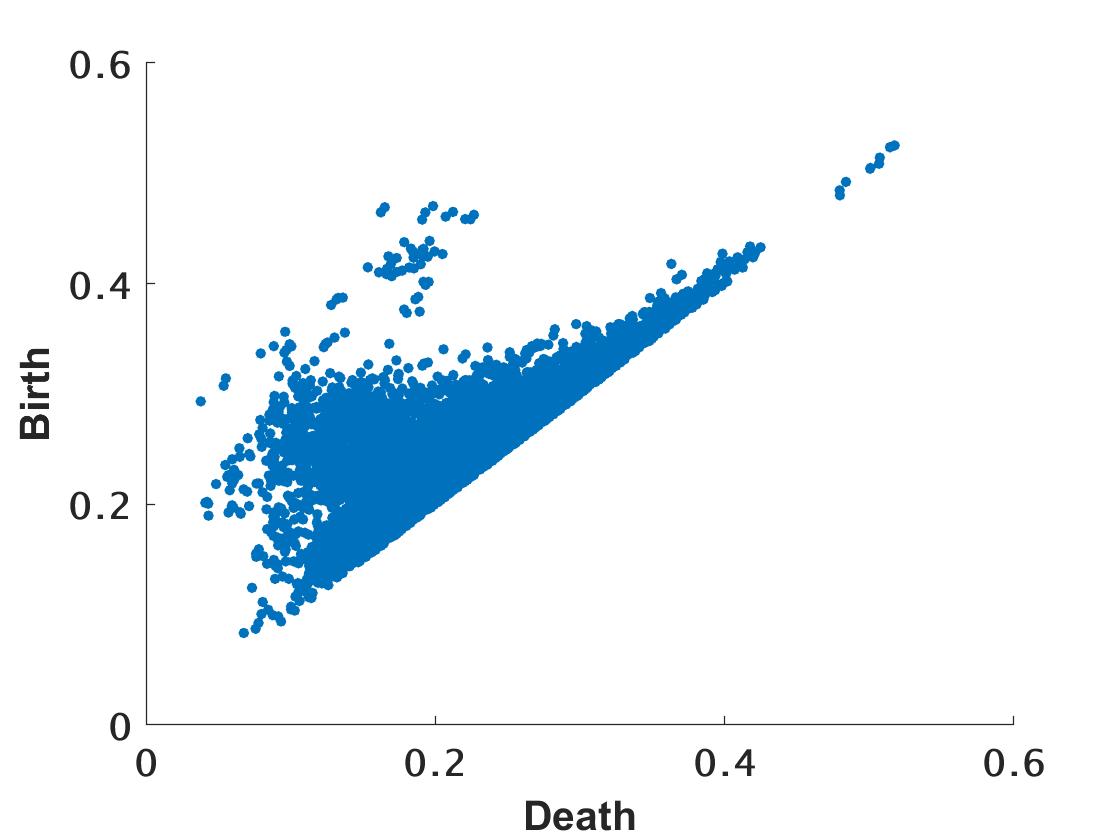}
    }
      \caption{\footnotesize  Superpositions of 100 persistence diagrams for  the sphere example. (a) $H_0$ for the original diagrams. (b) $H_0$ for MCMC simulations at step 500.   (c) $H_1$ for the original diagrams. (d) $H_1$ for MCMC simulations at step 500.}
 \label{fig:SpheresTotal}
\end{figure}

 That this reflected symmetry is  seen in MCMC simulations of fitted models to the $H_0$ and $H_1$ persistence diagrams can be seen from Figure \ref{fig:SpheresTotal} which presents a (by now familiar) summary of persistence diagram behaviour via superpositions of many diagrams. 
 
\begin{figure}[h]
    \centering
       \subfigure[]
    {
        \includegraphics[scale=0.08]{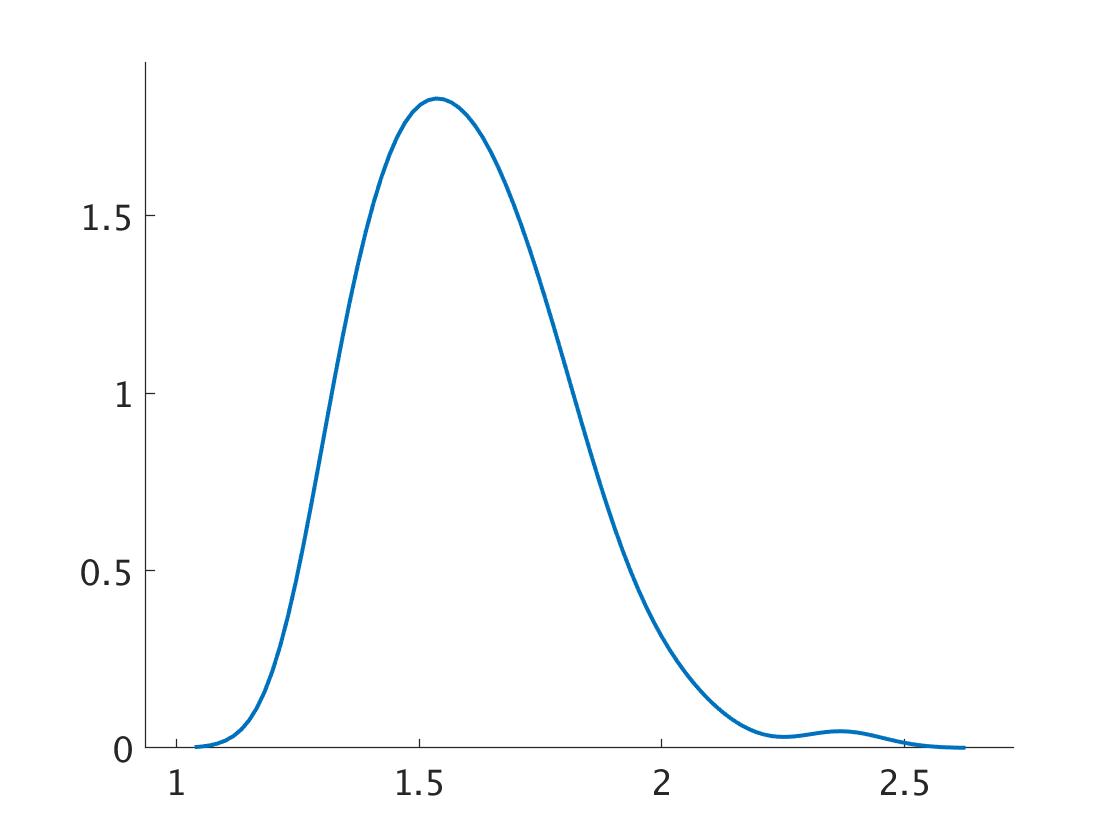}
    }
       \subfigure[]
    {
        \includegraphics[scale=0.08]{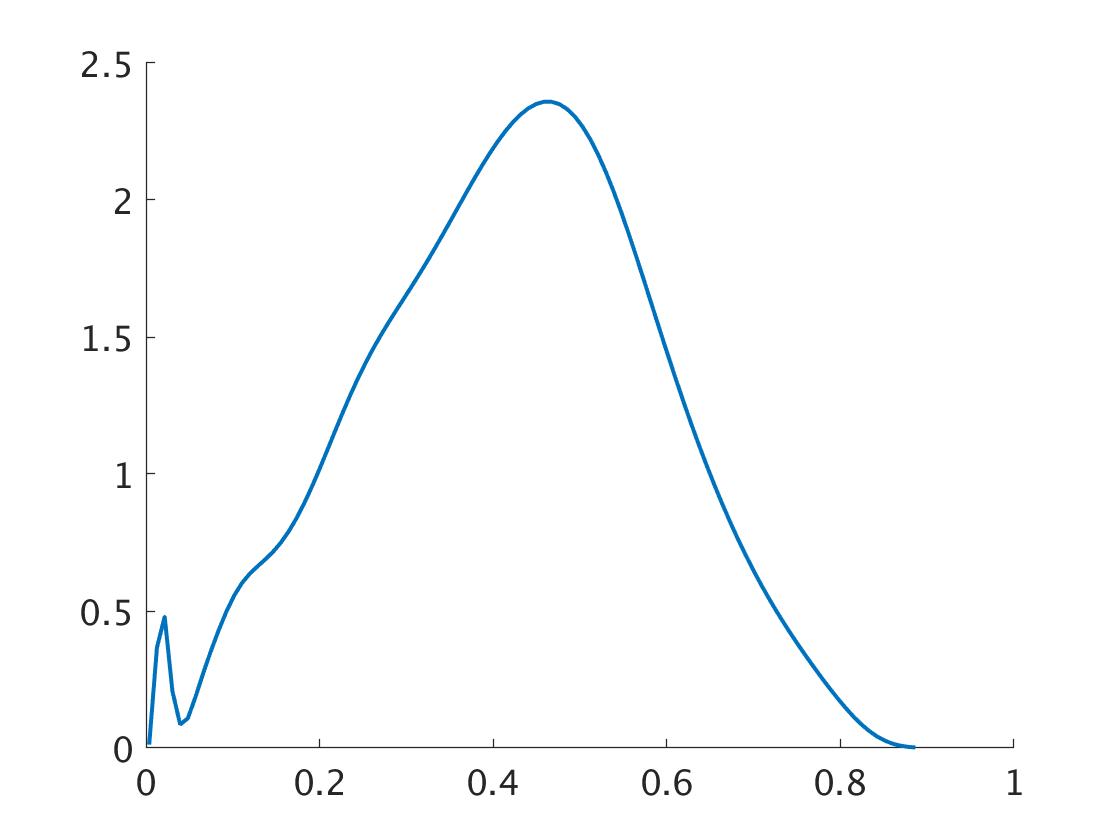}
    }
        \subfigure[]
    {
        \includegraphics[scale=0.08]{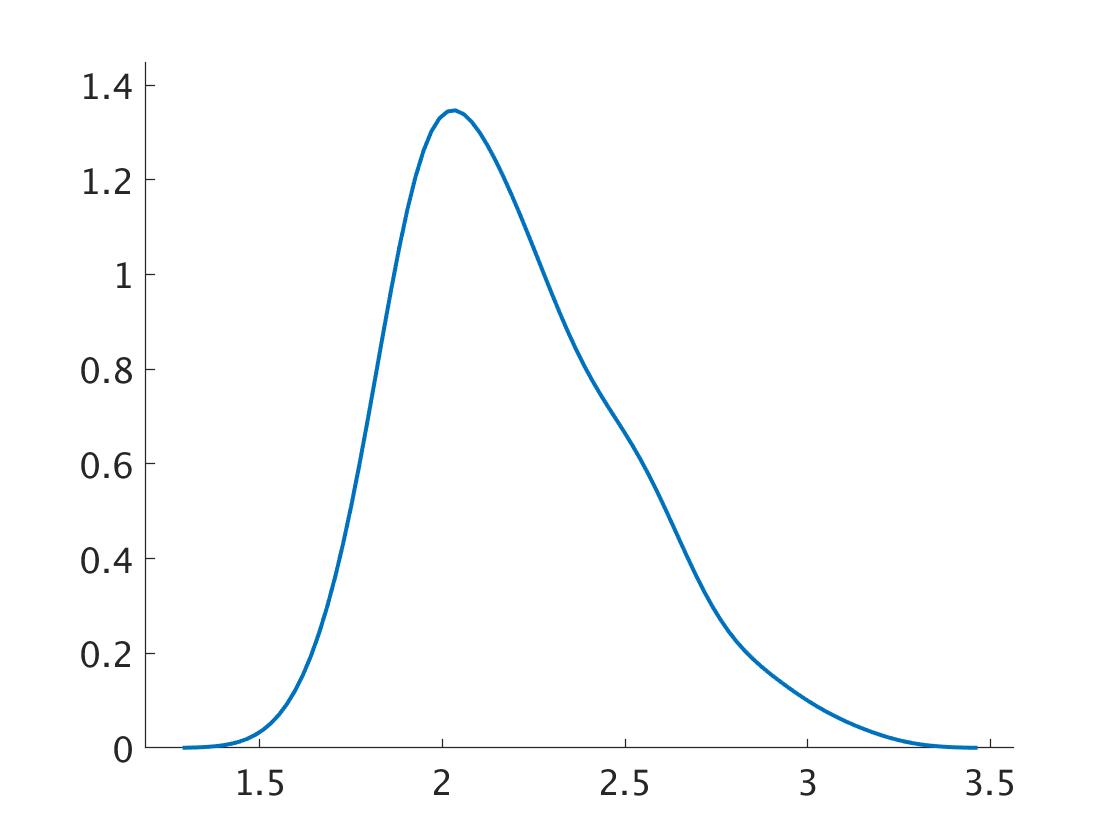}
    }
     \subfigure[]
    {
        \includegraphics[scale=0.08]{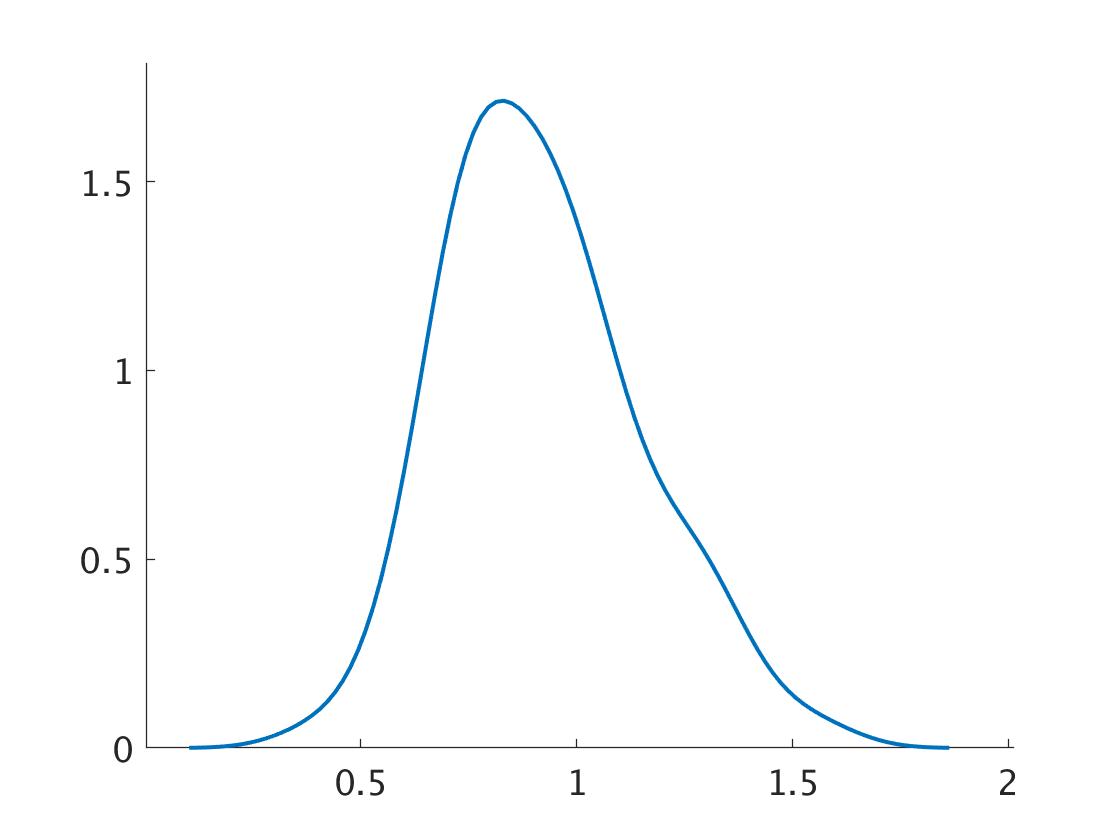}
    }
\caption{\footnotesize
 Smoothed empirical densities for the two normalised  parameter estimates   in the pseudolikelihood \eqref{eq:finalpseudo}  for the
$H_0$ and $H_1$ persistence diagrams, based on 100 simulations of a 2-sphere, in a model with $K=2$.  (a) $H_0$,  $\theta_1$. (b) $H_0$,  $\theta_2$.  (c)
 $H_1$,  $\theta_1$. (d) $H_1$,  $\theta_2$.}
\label{fig:sphereH0est}
\end{figure}

\begin{figure}[h!]
    \centering
    \subfigure[]
    {
        \includegraphics[scale=.08]{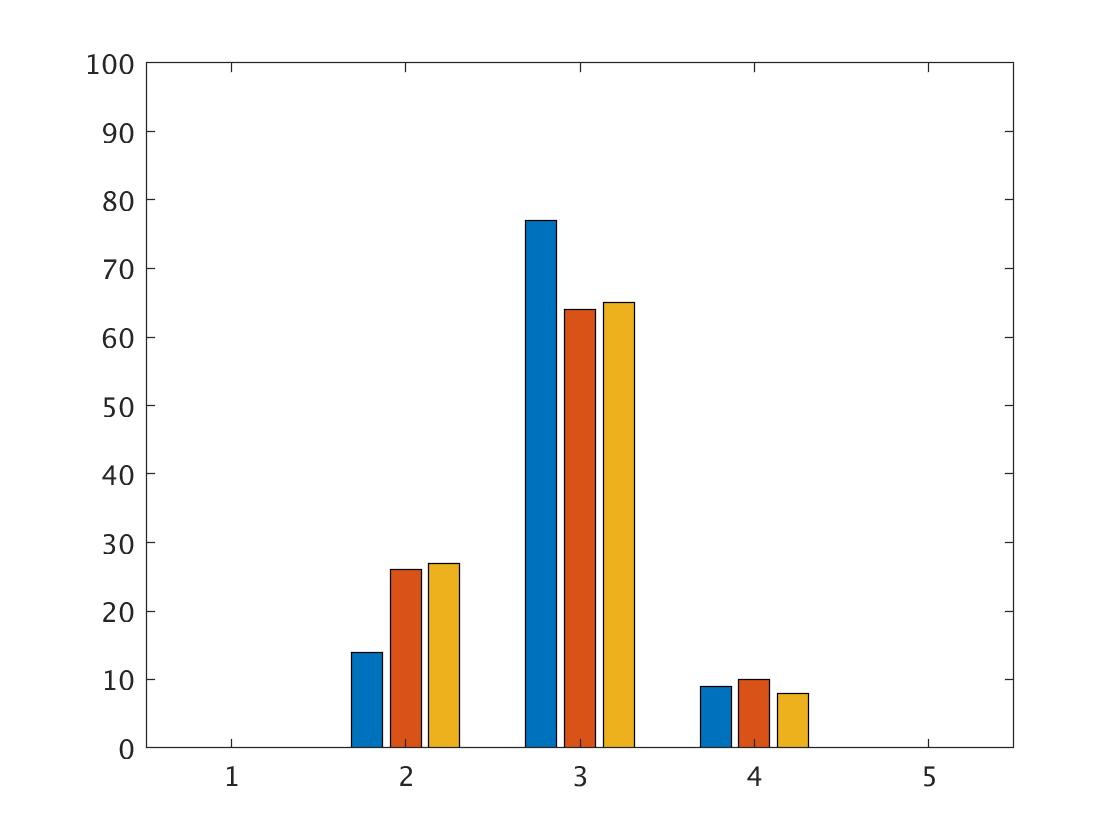}
    }
     \subfigure[]
    {
        \includegraphics[scale=.08]{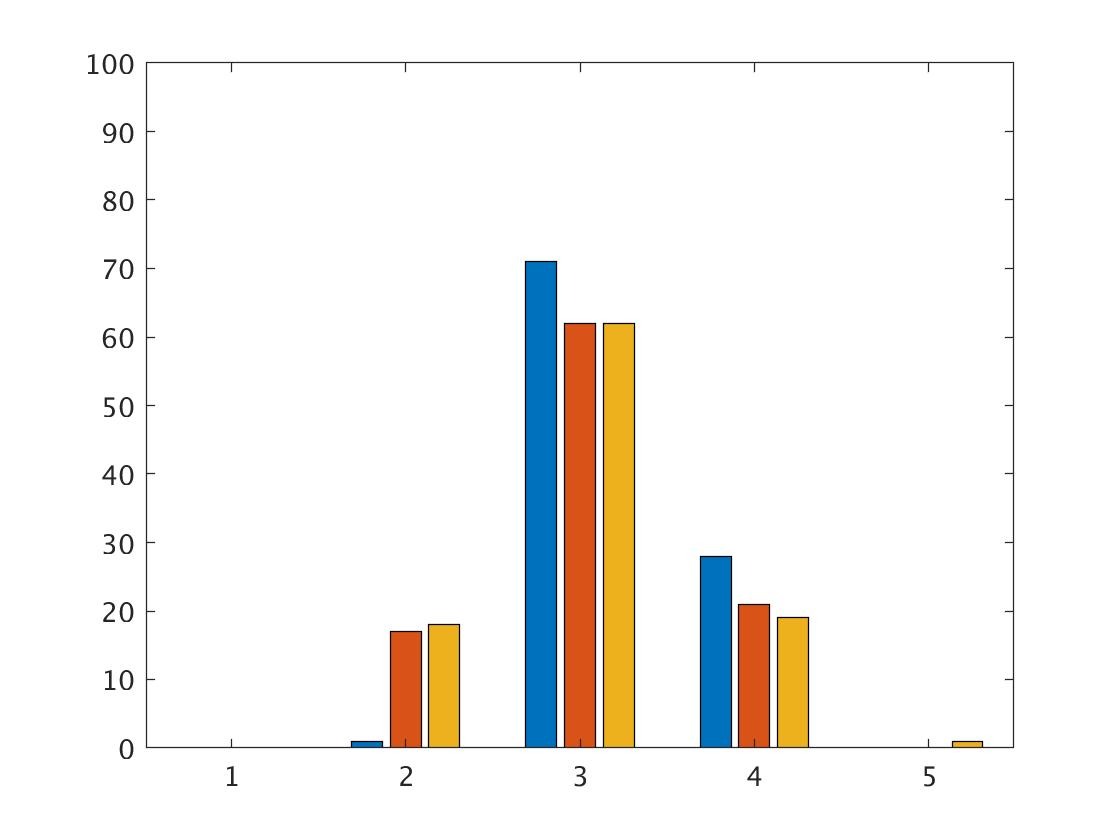}
    }
  \subfigure[]
    {
        \includegraphics[scale=.08]{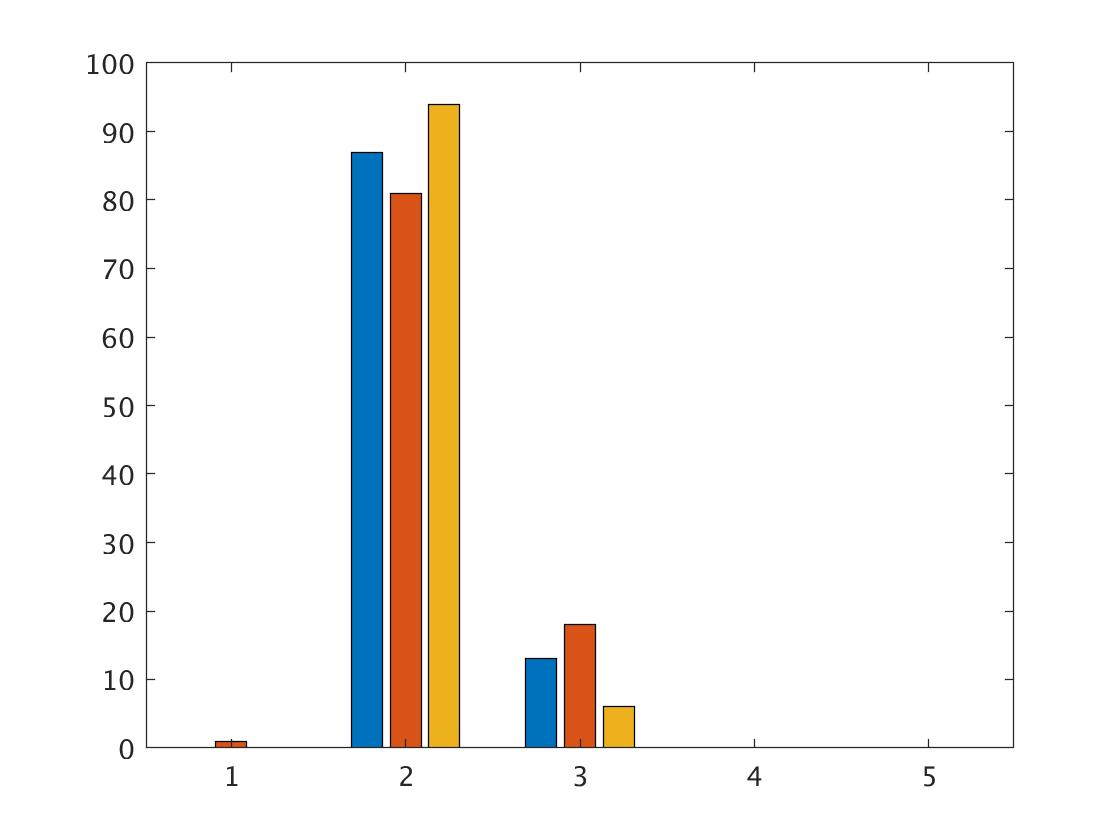}
    }
   \subfigure[]
    {
        \includegraphics[scale=.08]{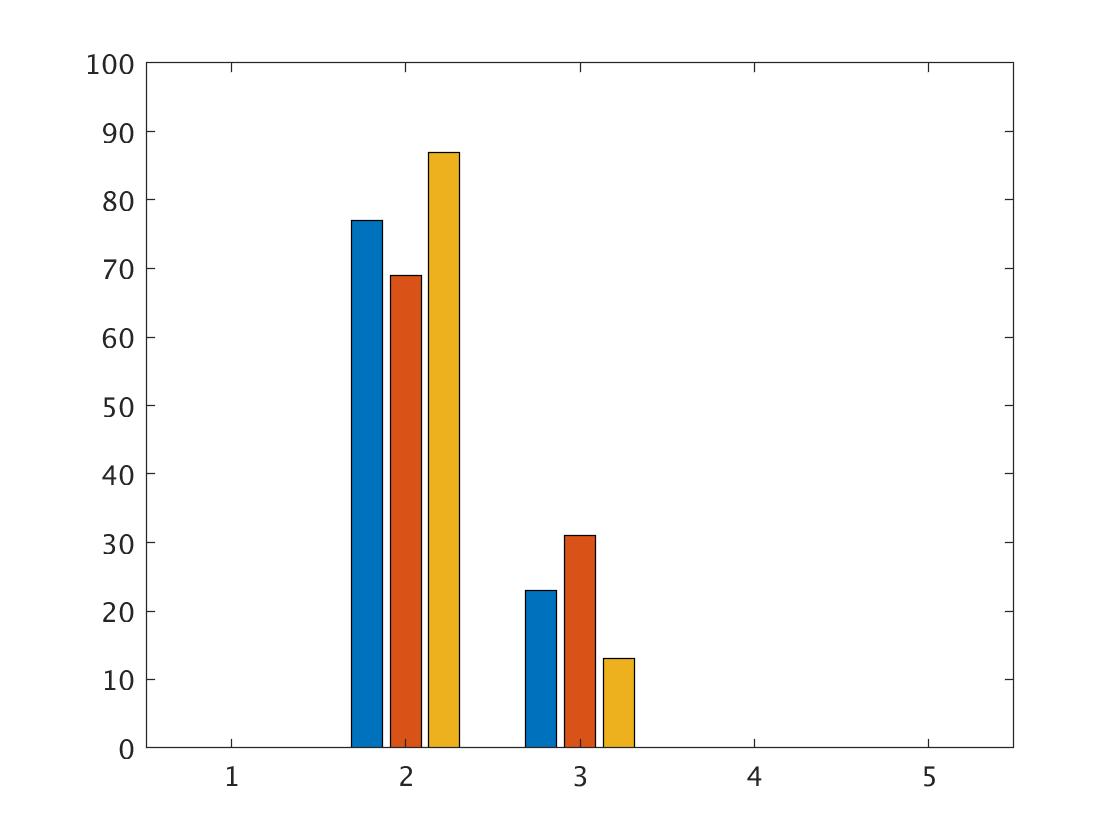}
    }
   \\
    \subfigure[]
    {
        \includegraphics[scale=.08]{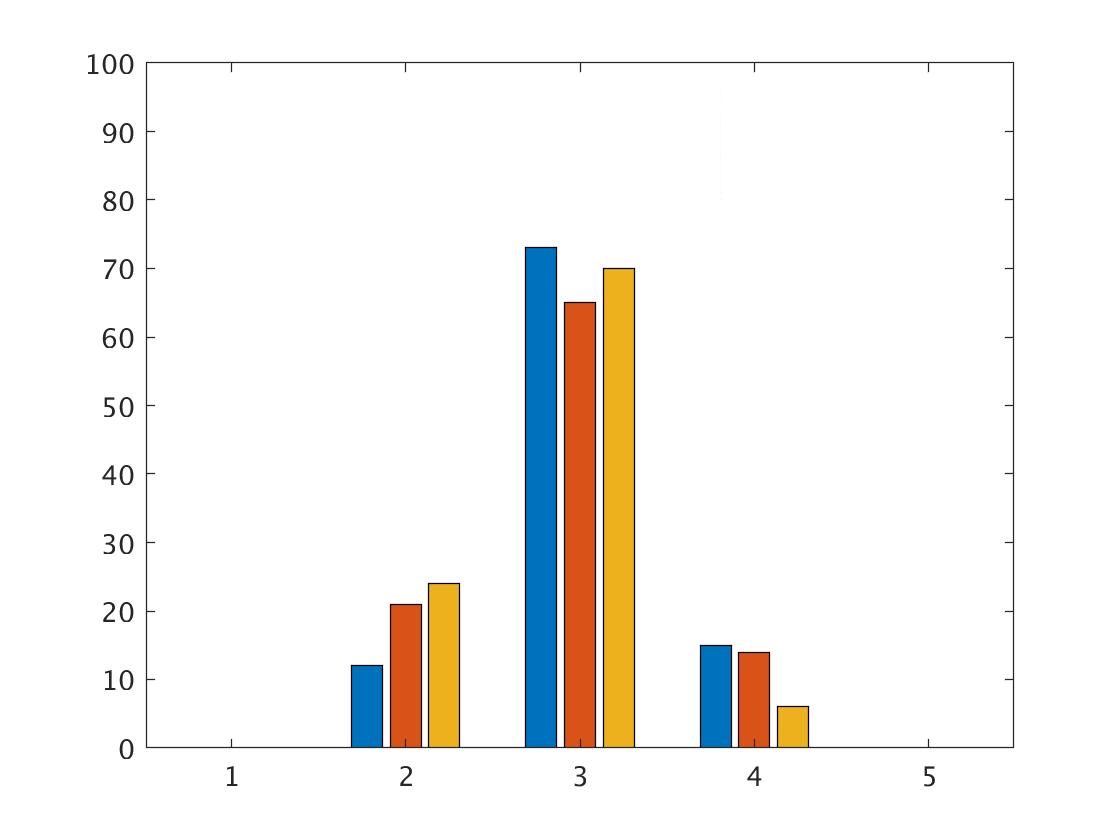}
    }
     \subfigure[]
    {
        \includegraphics[scale=.08]{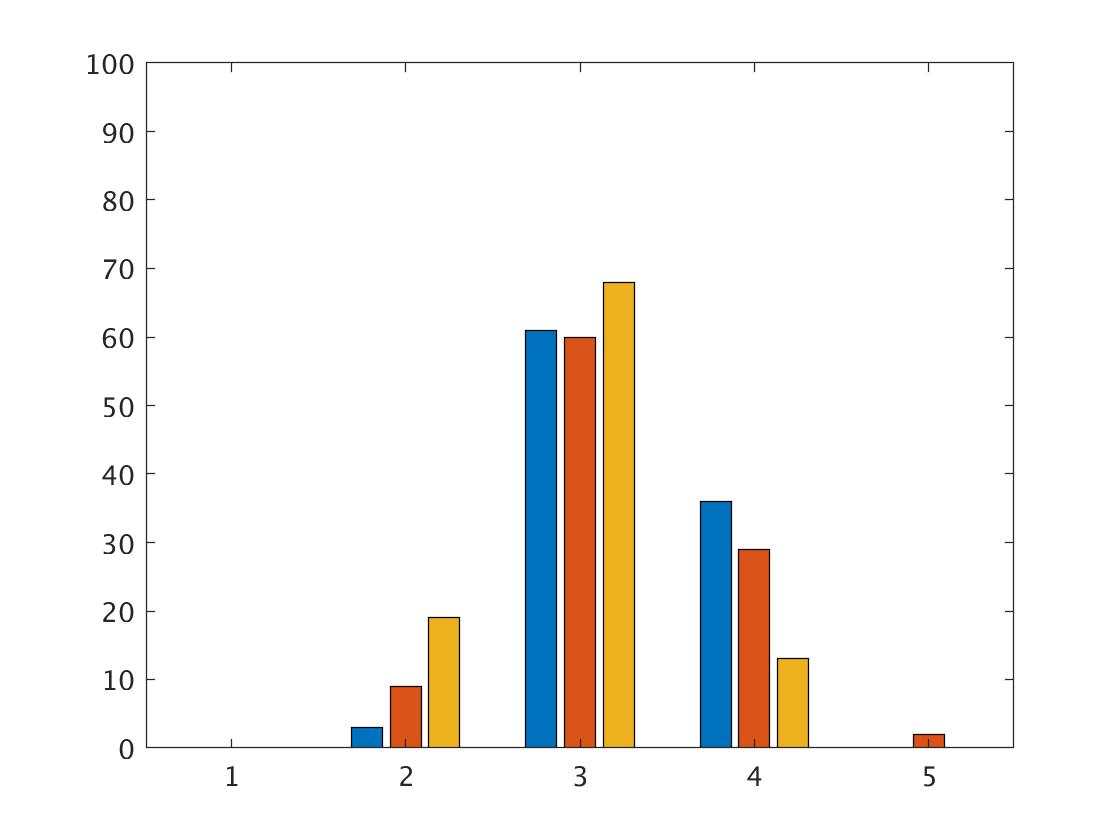}
    }
  \subfigure[]
    {
        \includegraphics[scale=.08]{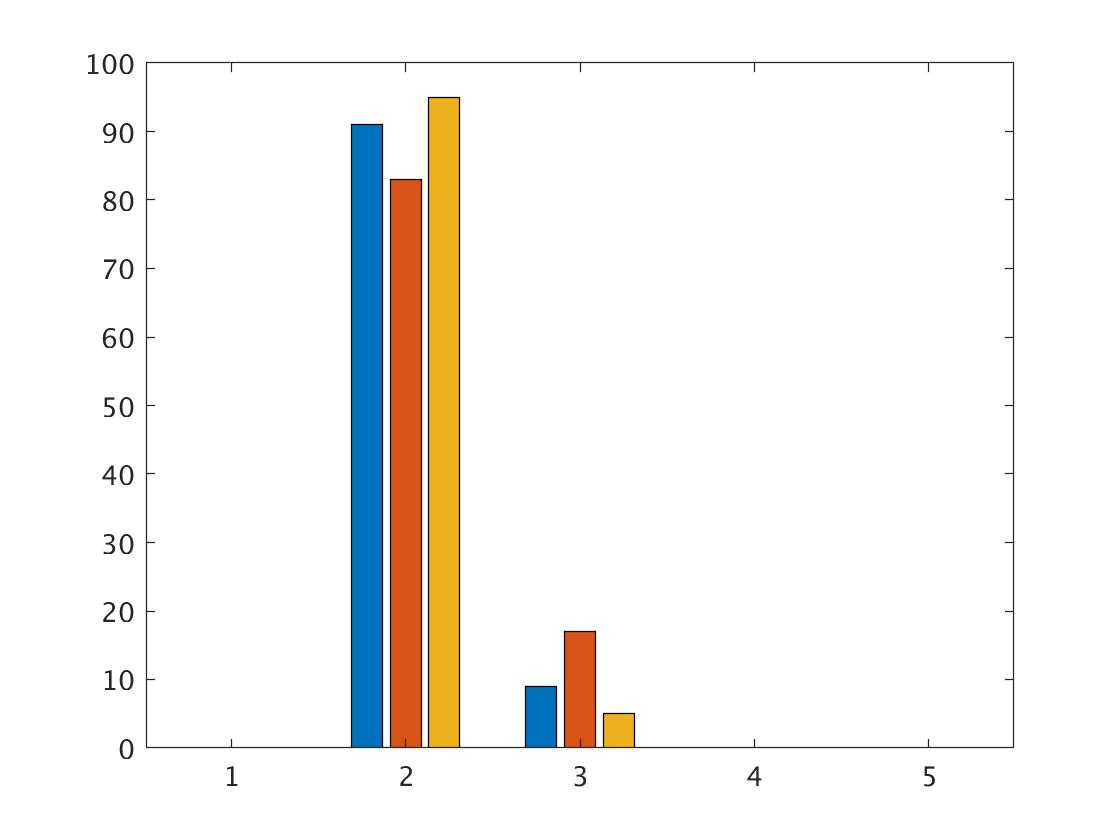}
    }
   \subfigure[]
    {
        \includegraphics[scale=.08]{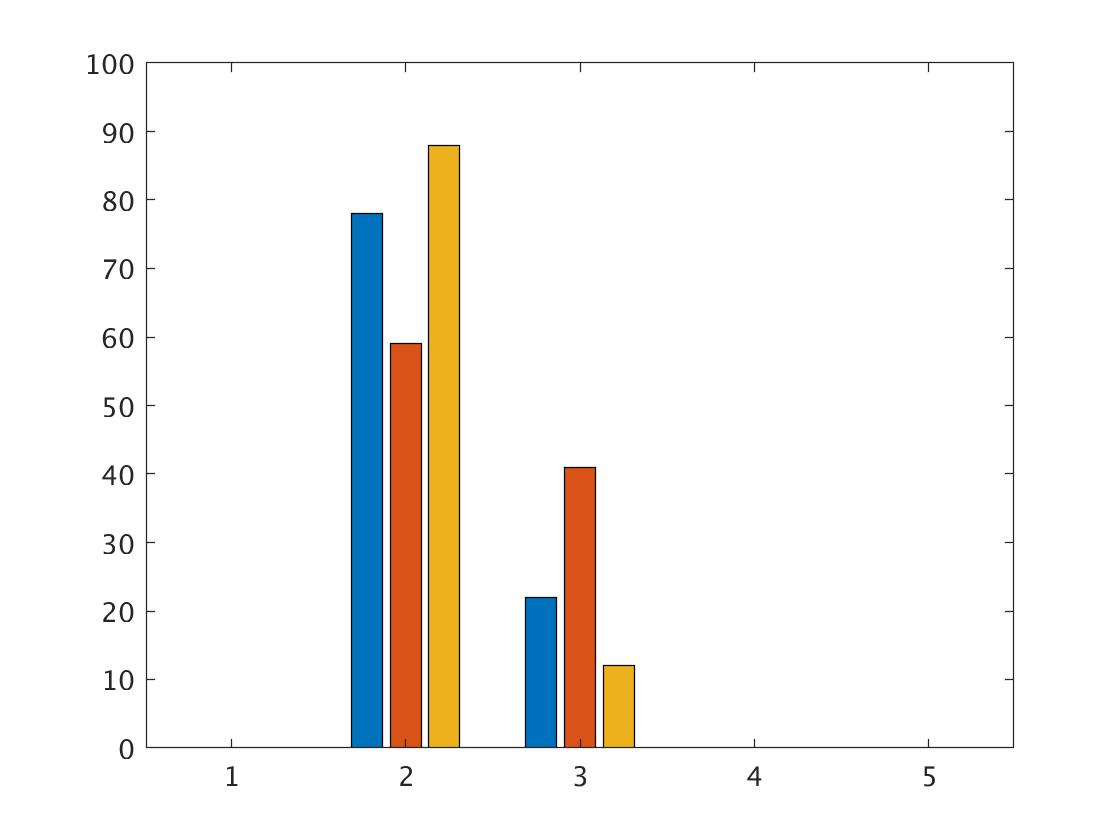}
    }

\caption{\footnotesize
   Distribution of numbers of clusters in the  persistence diagrams  for the spherical example.  
  (a)--(d) are for the $H_0$ diagrams, and (e)--(h) for $H_1$. Otherwise the colouring and ordering of the plots is as in the caption to Figure \ref{circlesCompare}.}
\label{2sphereH1Compare}
\end{figure}

Given this symmetry, it is also to be expected that parameter estimates and correlations will be similar for both homologies. That this is indeed the case can be seen in the graphs of Figure \ref{fig:sphereH0est}, which show the behaviour of the $\Theta$ estimates when a model with $K=2$ is fitted.  
Similarly, Table \ref{table:sphereCorrelations}, based  on a model with $K=3$, shows similar correlations between parameter estimates. Figure 
\ref{2sphereH1Compare} shows information on cluster sizes for this example, as did   Figure \ref{circlesCompare} for the concentric circles example. It is interesting that, once again, the points in the diagrams tend to group into a small number of well defined clusters, although for this example we do not have a clear topological explanation for the phenomenon.  However, as before, the distribution of the number of clusters is well preserved by the MCMC simulations.

\begin{table}[h!]
\begin{center}
\begin{tabular}{|l|ccc|}
\hline
 &$\rho(\theta_1,\theta_2)$ & $\rho(\theta_1,\theta_3)$ & $\rho(\theta_2,\theta_3)$ \\  \hline 
$H_0$  & 0.4844*  &  0.1848 &  -0.2079*\\
$H_1$ &0.2517* & 0.1766 &  -0.1345    
\\ \hline
\end{tabular}
\caption{\footnotesize{Correlations between parameter estimates for a model with $K=2$ for the sphere example and $H_0$ and $H_1$ persistence diagrams. The correlations are calculated  over 100 simulations.  Statistically significant correlations at the  5\% level are starred.}}
\label{table:sphereCorrelations}
\end{center}
\end{table}

\section{RST and bagplots: Detecting  topology}
\label{sec:bagplots}

\subsection{Bagplots}
\label{subsec:bagplots}

The notion of bagplots goes back to  the seminal paper  \cite{Tukey75}, which  heuristically introduced a notion of the `depth' of a particular  point, in a cloud of $n$ points, regardless of the dimension $d$ of the Euclidean space  in which they lie. 
Known for its originator as the `Tukey depth', it is the the minimum number of  points in the full collection which lie on one side of a hyperplane through the point. The `Tukey median' of the full set is a point maximizing the Tukey depth, and a `Tukey centerpoint'  is a point of depth at least $n/(d + 1)$. A  median must be a centerpoint, but not every centerpoint is a  median. Neither need be unique. 

The most useful application  of this notion  of depth is that it provides points which have, naturally, only the partial ordering of $\real^d$ with, up to ties,  a total ordering. In particular, the deeper a point is relative to the  cloud, the higher its  depth value. This ordering has proved extremely powerful in extending univariate tools  related to signs and ranks, order statistics, quantiles, and   outliers  to the multivariate setting in a unified way. (See \cite{Serfling} for a general overview.)

Tukey depth led naturally, in \cite{Bagplot}, to the tool now know as the bagplot (sometimes called a starburst plot) which  is a method in robust statistics for analysing high dimensional data, analogous to the one-dimensional box plot. 
In dimensions $d=2$ and $d=3$ it also lends itself to easy visualisation, as in Figure \ref{fig:bagplotexamples} below. Analogously to the box plot, in these lower dimensions it allows for the visualisation of the location, spread, skewness, and outliers of the point cloud. 

A bagplot plot, for two dimensional point sets, consists of three nested, convex, polygons, called the `bag', the `fence', and the `loop'. The points are first ordered according to Tukey depth,  from highest to lowest depth, with ranking associated at random among ties.  The  inner polygon, called the bag, contains (no more than) the first 50\% of the points according to this ranking, and is defined by the convex hull of its members. 

The outermost of the three polygons, called the fence, is not drawn as part of the bagplot, but is used to construct it. It is formed by inflating the bag by an `inflation   factor' $c>0$. (The default, based on experience, is to take $c=3$). The inflation is carried out with respect to the Tukey median, if it is uniquely defined, or, otherwise, with respect to the center of gravity of the bag.

 Observations outside the fence are flagged as outliers. The observations that are not marked as outliers are surrounded by the loop, the convex hull of the observations within the fence.

Higher dimensional bagplots can be defined similarly, although it is only the three dimensional ones that are visualisable. In this case the bagplot consists of an inner and outer bag, with the outer bag  drawn in transparent colors so that the inner bag remains visible.

Bagplots are invariant under affine transformations of Euclidean space, and generally robust against outliers. For details of the Matlab  implementation that we used, see  \newline {\tt physionet.org/physiotools/ecg-kit/common/LIBRA/bagplot.m}, and for general implementation details see \cite{HubertStephan,LiuZuo}.

Figure \ref{fig:bagplotexamples} contains four  examples of bagplots for persistence diagrams. The examples were chosen  to demonstrate a number of phenomena of relevance to TDA and the analysis of the next subsection.

  \begin{figure}[h!]
    \centering
    \subfigure[]
    {
        \includegraphics[scale=0.08]{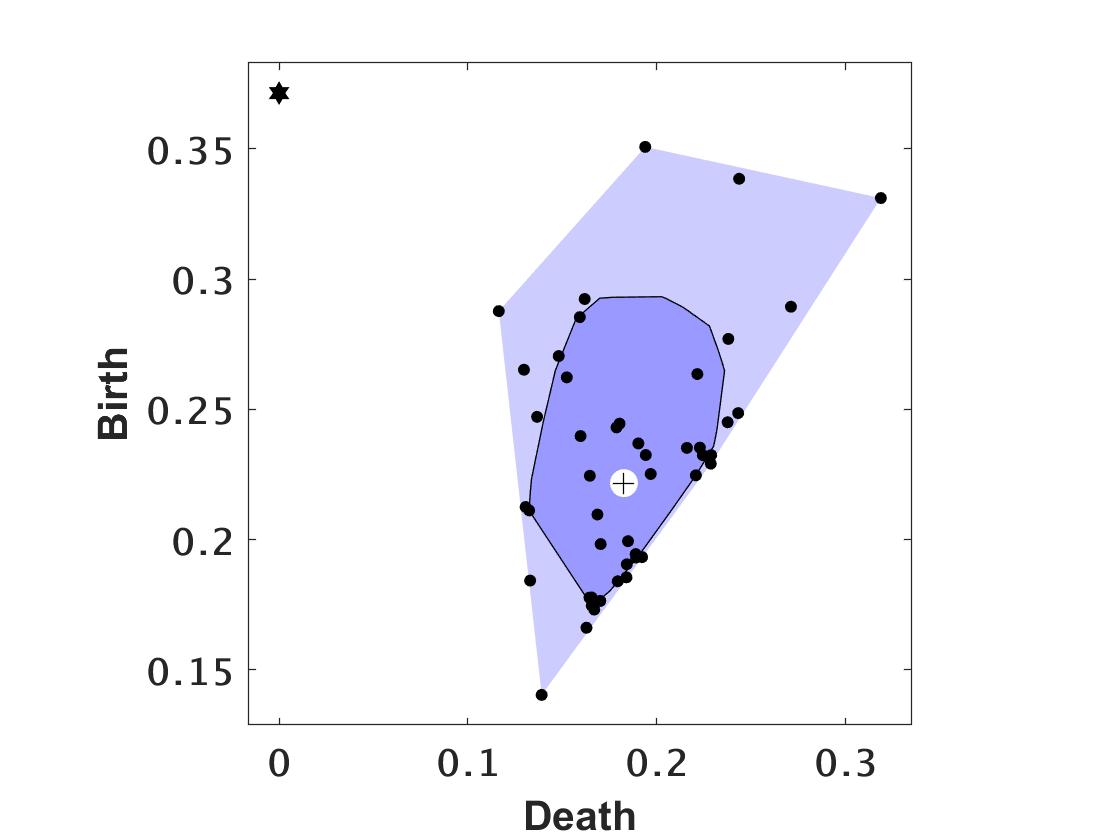}
    }
       \subfigure[]
    {
        \includegraphics[scale=0.08]{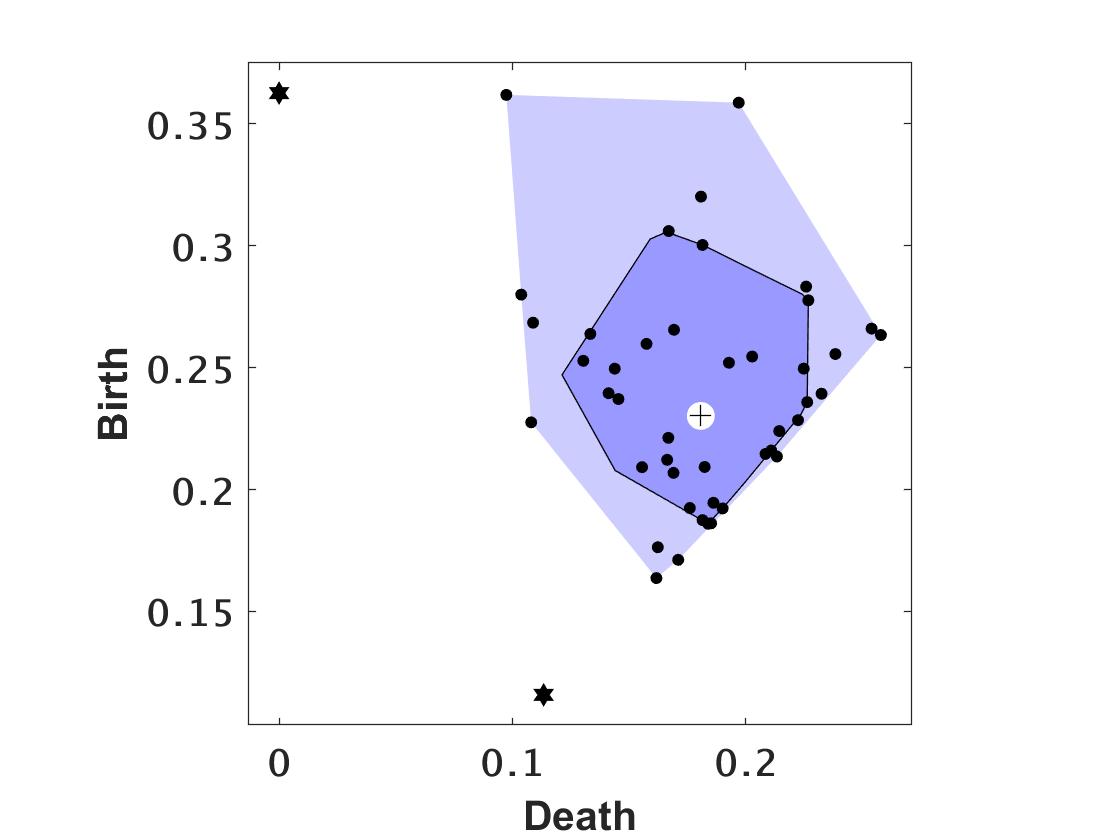}
    }
     \subfigure[]
    {
        \includegraphics[scale=0.08]{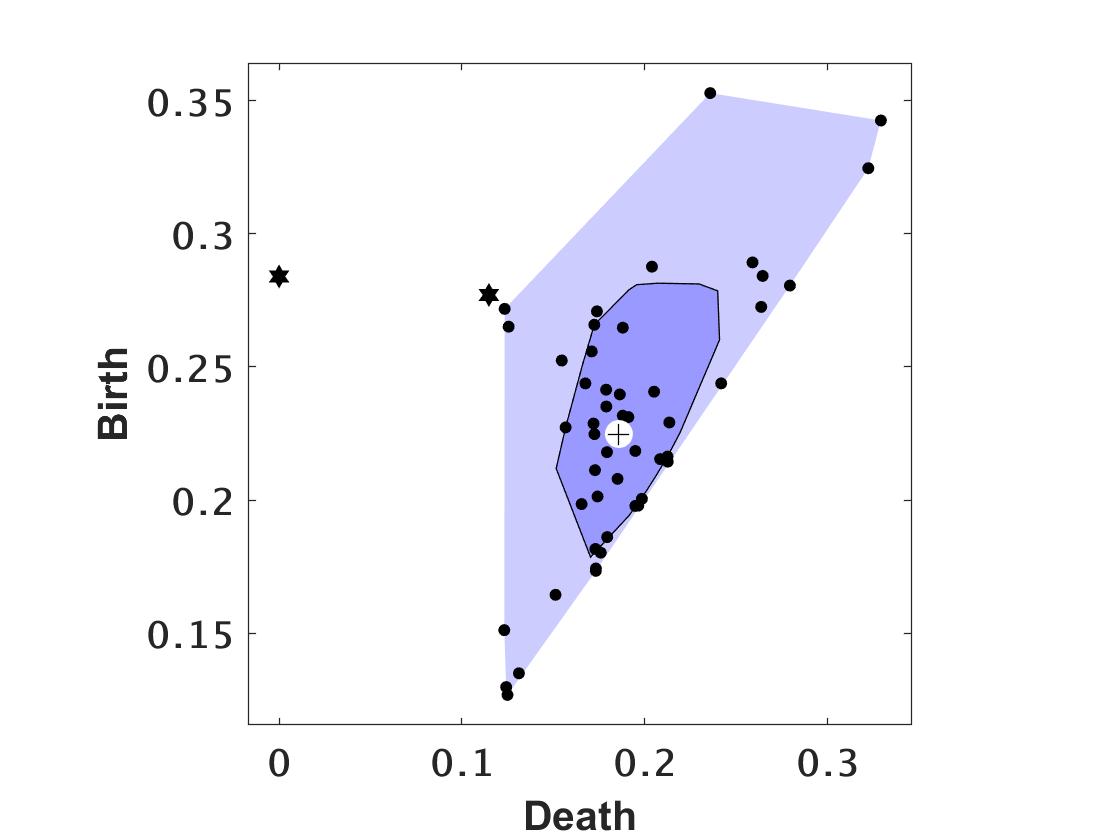}
    }
     \subfigure[]
        {
        \includegraphics[scale=0.08]{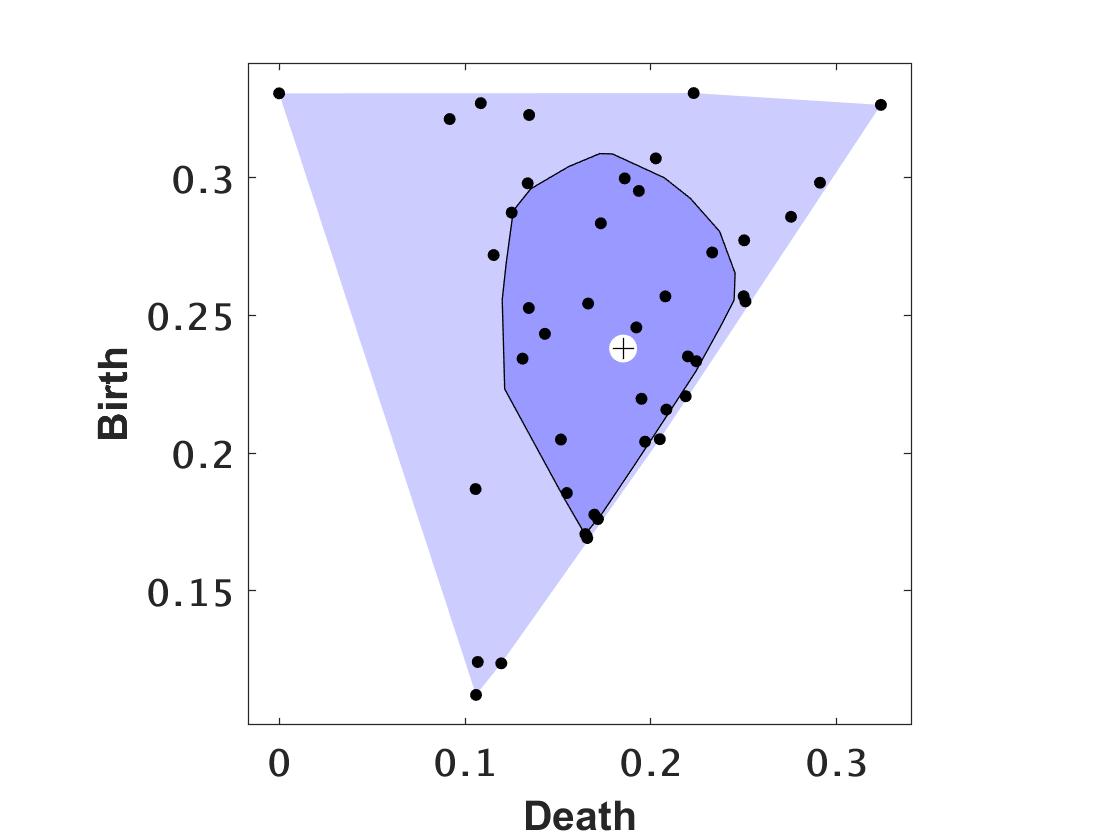}
    }
    \caption{\footnotesize  Bagplots of persistence diagrams from the concentric circles example. See text for details.}
\label{fig:bagplotexamples}
\end{figure}

The bags and loops are all clear in Figure \ref{fig:bagplotexamples}, with the former a darker color than the interior of the latter. Outliers are marked with stars, and the Tukey median with a plus. 
All the persistence diagrams are $H_0$ diagrams and come from the concentric circles example of the previous section. Since the point at infinity has been removed, we would hope to be able to identify one point of the diagram as signalling a component, and  classify the remaining points as noise. The four examples here do this with varying degrees of success.

Example (a) is what we would hope for. The single outlier at death time 0 is clearly identified, and identifies the component we are looking for. Example (b) is reasonable. Although there are two outliers here, the additional one is so close to the diagonal that it is clearly a consequence of noise, and no statistical test is needed to know it is of no consequence. In general, points of this kind, which arise due to very small, very  local, fluctuations in the empirical density of the diagram, are unavoidable but meaningless. Example (c) also has an additional point identified as an outlier, which `only just' makes it into this category. Looking at the diagram, one would be inclined to ignore it, but any automated statistical procedure would classify it as significant, as indeed it should. (This, after all, is what Type I error rates are all about in Statistics.) Example (d) is disappointing, as a rather unusual arrangement of points at the top of the diagram lead to the `signal point' being included within the fence.

The summary from Figure \ref{fig:bagplotexamples} is therefore that only in one out of four cases do bagplots correctly identify the `signal point'  as an outlier. However, it turns out that situations such as  (c) and (d) actually occur rarely (i.e.\ with low probability) in the sampling problems being discussed here, and so they fall, respectively,  into the classes of Type I and Type II  Errors; viz.\ random fluctuations that can lead to false discoveries or to failing to identify true phenomena. Situation (b) occurs more often, but the position of the outlier, along the diagonal, shows little in the way of consistency across different samples. Furthermore, as we already noted, these points are so close to the diagonal that they are easily discounted. 

To see  how bagplotting works in practice, we now look at a specific example, combining the technique  with the replication methodology developed throughout the paper.

\subsection{Bagploting for topology: An example of RST in action}

We continue now with the example of Section \ref{sec:concentric}, based on sampling from two concentric circles. Our aim is to see how well we can detect the single `real' outlier in the persistence diagrams such as those in Figure  \ref{fig:bagplotexamples} and, more importantly, how we can assign levels of statistical significance to detections.

Some notation for this, in a very general setting (i.e.\ not restricted to concentric circles): Let  ${\cal P} = \{P_1,\dots,P_{m}\}$ be a collection of $m$  original persistence diagrams. In most cases, unfortunately, we will have $m=1$. Using the techniques of earlier sections, let 
$\widehat {\cal P}_k=\{\hat P_{k,1},\dots \hat P_{k,n}\}$, $k=1,\dots,m$, be a collection of  $n=n1+n_2$  MCMC simulations associated with the original diagram $P_k$, after model fitting.  Let  ${\cal B}^c = \{B^c_1,\dots,B^c_{m}\}$ be  the collection of $m$ bagplots for the diagrams of $\cal P$,  with inflation factor $c$, and     $\widehat {\cal B}^c_k=\{\hat B^c_{k,1},\dots \hat B^c_{k,n}\}$ the  corresponding bagplots for the diagrams of  $\widehat {\cal P}_k$.

In most practical situations, of course, one would have only one sample persistence diagram  ($m=1$), which we denote by $P$,  one set of MCMC simulations, $\widehat {\cal P}$, and one set of bagplots,   $\widehat {\cal B}=\{\hat B^c_{1},\dots \hat B^c_{n}\}$, and so we shall first develop a general statistical procedure for identifying topological signal points in this situation. 

We first need to fix some  external parameters, chosen from general considerations for the problem at hand; viz.\ 
\begin{itemize}
\item[(i)] Fix a $p^*\in [0,1]$ as an acceptably small probability of making initial bad guesses about how many topologically significant points there might be in the diagram. We found that $p^*=0.05$ was a reasonable choice.
\item[(ii)] A range $[C_*,C^*]$ of possible values for the bag inflation factors, along with a step size $\delta>0$, giving a collection 
${\cal C}_\delta =\{C_*,C_*+\delta,C_*+2\delta,\dots,C^*\}$ of inflation factors.
\item[(iii)] An integer $A$, denoting an upper bound (not necessarily tight) to the number of topologically significant points that, {\it a priori}, might be expected.
\end{itemize}
With these choices, the procedure is as follows:

\begin{algorithm}
\caption{Bagplot detection of topological signal in persistence diagrams.}
\label{Boxplot:algorithm}
\begin{algorithmic}[1]
\State  For each $c\in {\cal C}_\delta$, and each $\hat B^c_{1},\dots \hat B^c_{n_1}$, compute the number of outliers, denoted by
$O_1^c,\dots,O^c_{n_1},$ in each  bagplot.   
\State Choose an inflation factor $C^*$. One way to do this, which we found useful in practice, is to take 
\beq
C^* \ = \ \min\big\{c\in{{\cal C}_\delta}: \  n_1^{-1} \sum_{i=1}^{n_1} {\bf 1}_{O^c_i \geq  A} \ \leq \ p^*.
\big\}.
\label{eq:star}
\eeq
\State  For each of the  bagplots  $\hat B^{C^*}_{n_1+1},\dots \hat B^{C^*}_{n_1+n_2}$
 find the points in the original diagram $P$ which lie beyond the fence; i.e.\ are outliers. 
 \State To each point $x$ in the original diagram $P$, let $f(x)$ be the proportion of times that $x$ is classified as an outlier in the previous step.
 Declare $1-f(x)$ to be the `$p$-value associated with $x$'. That is, if $x$ is now classified as `topologically significant',  then 
  $1-f(x)$  can be  interpreted as an estimate of the probability that this classification is incorrect. 
\end{algorithmic}
\end{algorithm}

The above algorithm is designed to handle a single original persistence diagram. In the fortunate situation that replications of the underlying experiment led to a number of diagrams, there are at least two natural approaches for adapting this algorithm. The first is to simply apply  it individually to each diagram, identifying $p$-values throughout, and then combine the results with any of the standard methods of multiple testing. The second is to exploit the fact that there are multiple diagrams, and use the original diagrams themselves, or  go a short distance into the MCMC to take mild perturbations of them, rather than the first $n_1$ MCMC simulations as in Algorithm \ref{Boxplot:algorithm},  for computing the $C^*$ of \eqref{eq:star}. 

The reader will certainly be able to think of other variations of the general procedure for the multiple diagram case, and, indeed, natural adaptions of Algorithm \ref{Boxplot:algorithm}. Our aim here is not to present a fully polished and tested approach, but rather to introduce the idea of using bagplots, and to show how this fits naturally into a RST framework.

To see how well these ideas work in practice, we carried out 
a numerical experiment, using the 100 samples of 800 points on concentric  circles as in Section  \ref{sec:concentric}, and for each one computed a $H_0$ persistence diagram for its empirical density. We then followed the procedure of Algorithm  \ref{Boxplot:algorithm},
with the adaption of using the 100 original diagrams to estimate $C^*$, as described above. With $A=2$ and  $p^*=0.05$, this led to $C^*=2.92$, a little smaller than the standard default value of 3. 

In order to summarise the results, in each of the 100 original diagrams, $P_k$, we ordered  the frequencies $f(x)$ from largest to smallest, so that 
\beqq
f^*_k(1)  &= \max_{x\in P_k} f(x),  \qquad & p^*_k(1) = \text{argmax}_{x\in P_k} f(x),\\
f^*_k(j+1)  &= \max_{x\in P_k\setminus \cup_{i=1}^j  p^*_k(i)   } f(x),   
&p^*_k(j+1) = \text{argmax}_{x\in P_k\setminus \cup_{i=1}^j  p^*_k(i)} f(x),\qquad j\geq1.
\eeqq

\begin{figure}[h!]
    \centering
       \subfigure[$\varepsilon=0$, First outlier] {
        \includegraphics[scale=.11]{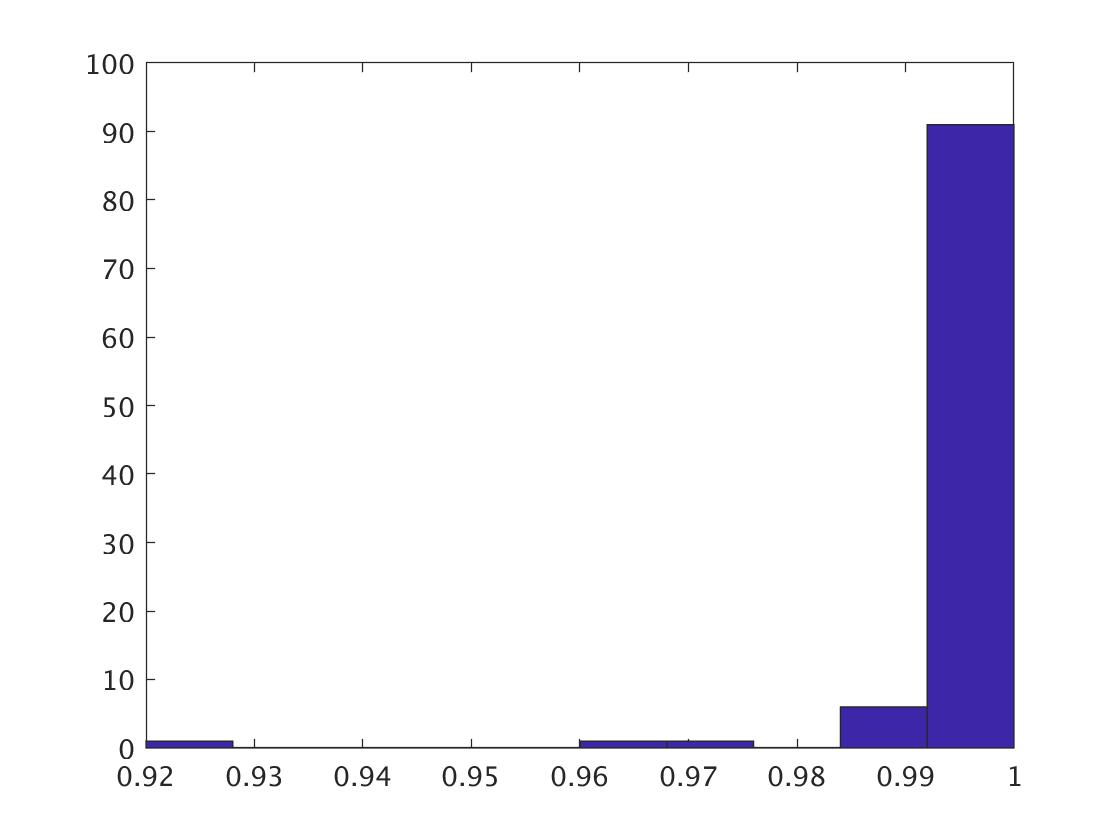}
    }
       \subfigure[$\varepsilon=0$, Second outlier]
    {
        \includegraphics[scale=.11]{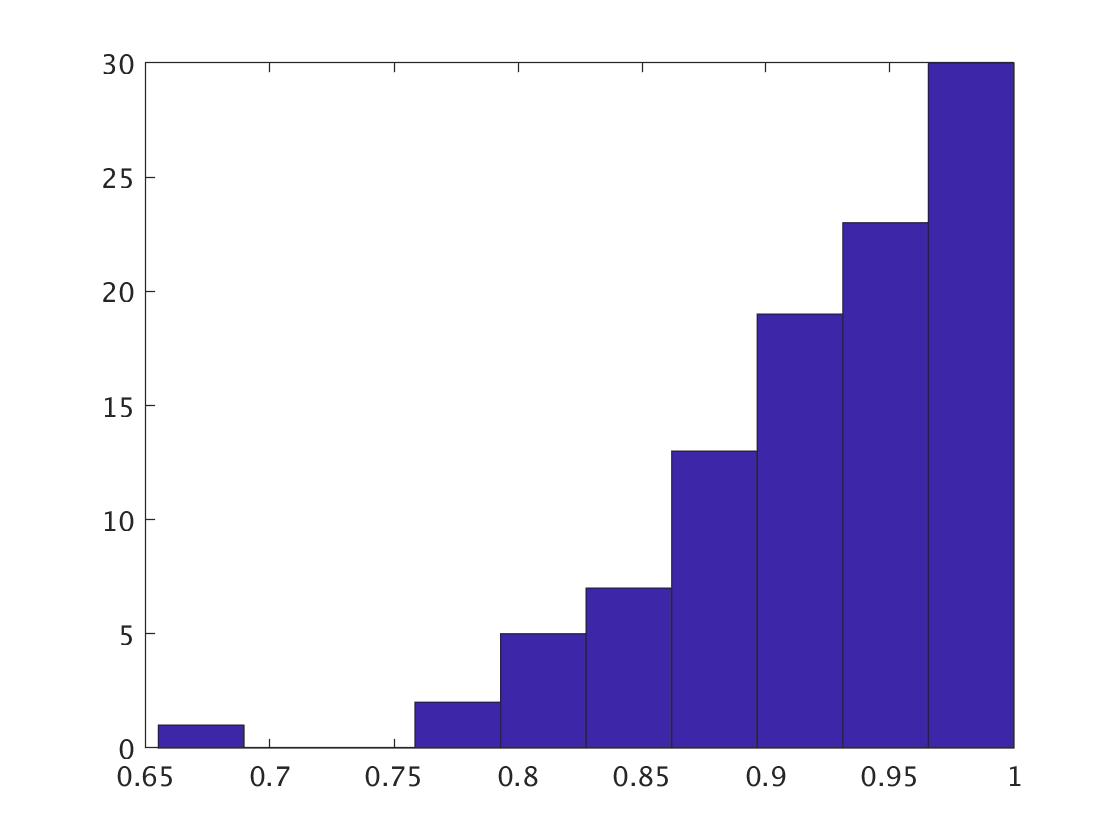}
    }
     \subfigure[$\varepsilon=0$,  Third outlier]
    {
        \includegraphics[scale=.11]{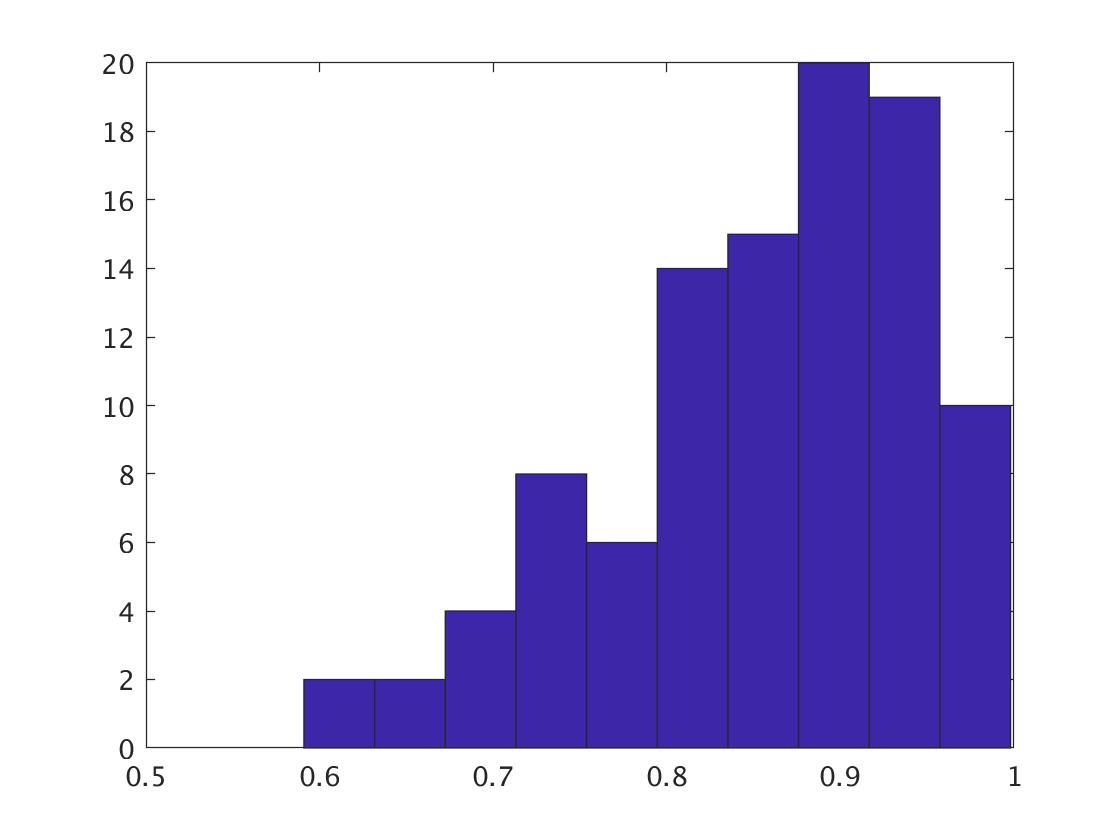}
    }
       \subfigure[$\varepsilon=0.001$, First outlier] {
        \includegraphics[scale=.11]{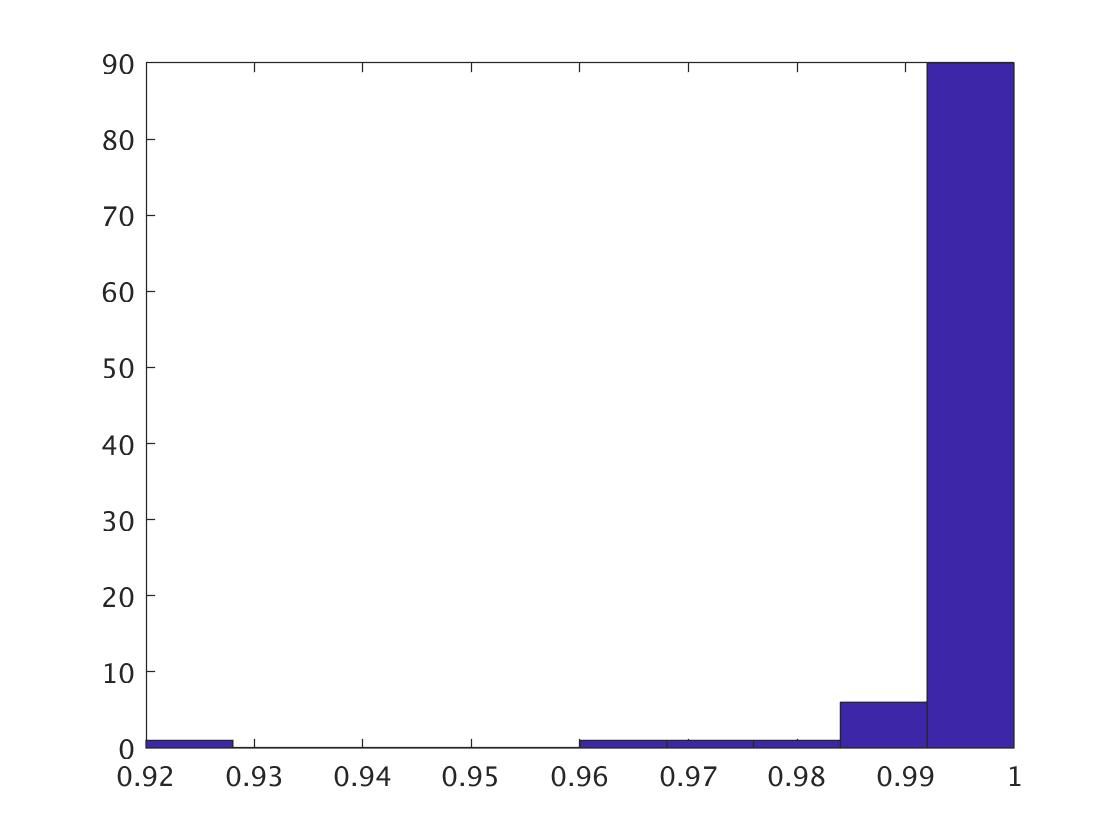}
    }
       \subfigure[$\varepsilon=0.001$, Second outlier]
    {
        \includegraphics[scale=.11]{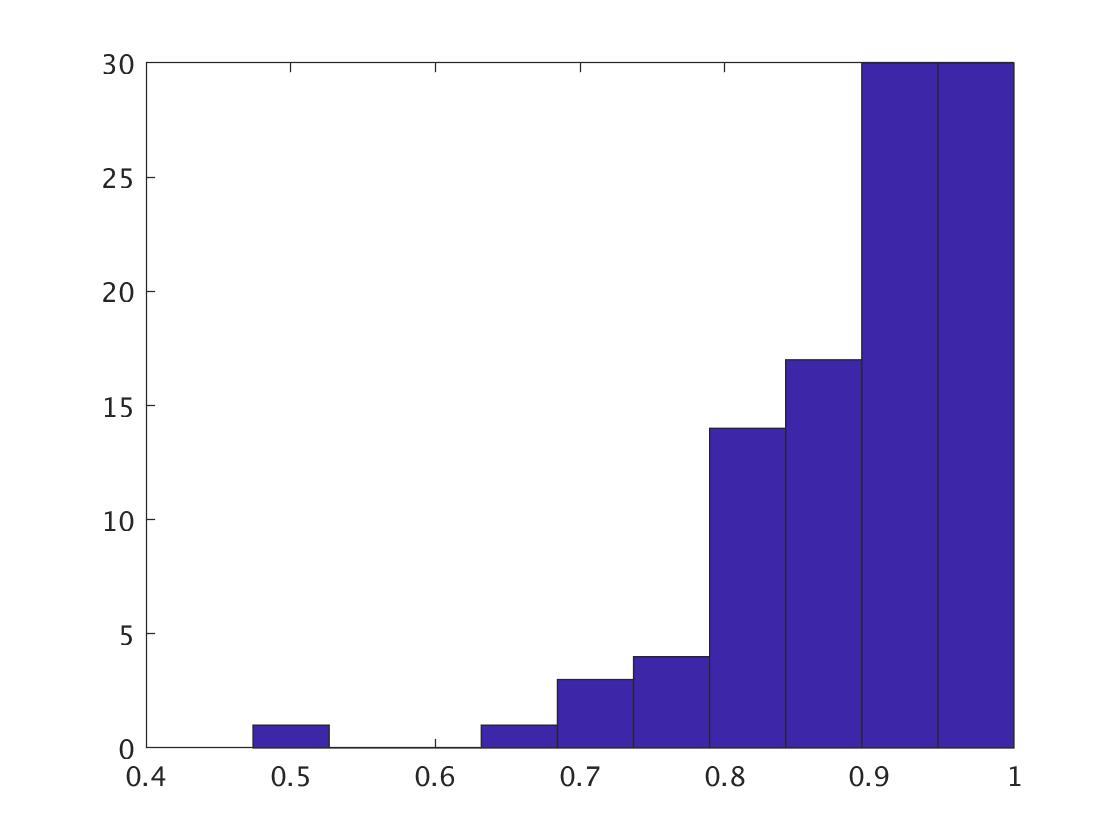}
    }
     \subfigure[$\varepsilon=0.001$,  Third outlier]
    {
        \includegraphics[scale=.11]{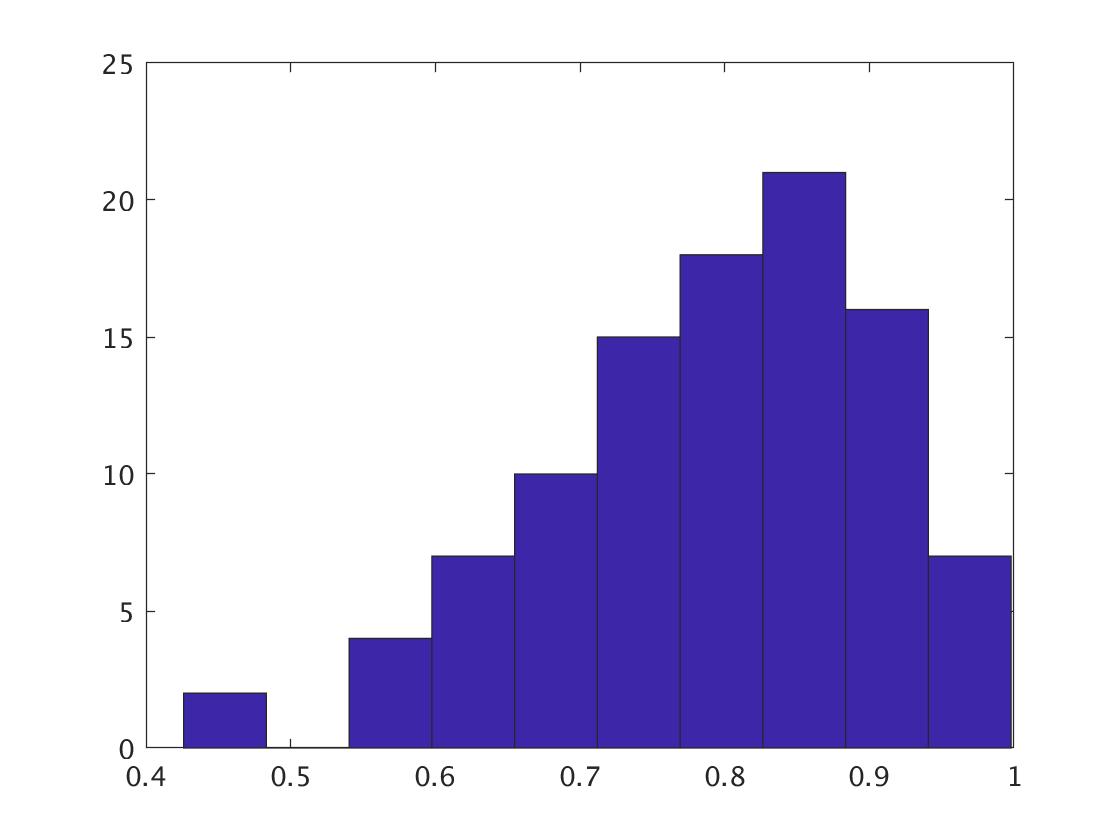}
    }    
    
       \subfigure[$\varepsilon=0.002$, First outlier] {
        \includegraphics[scale=.11]{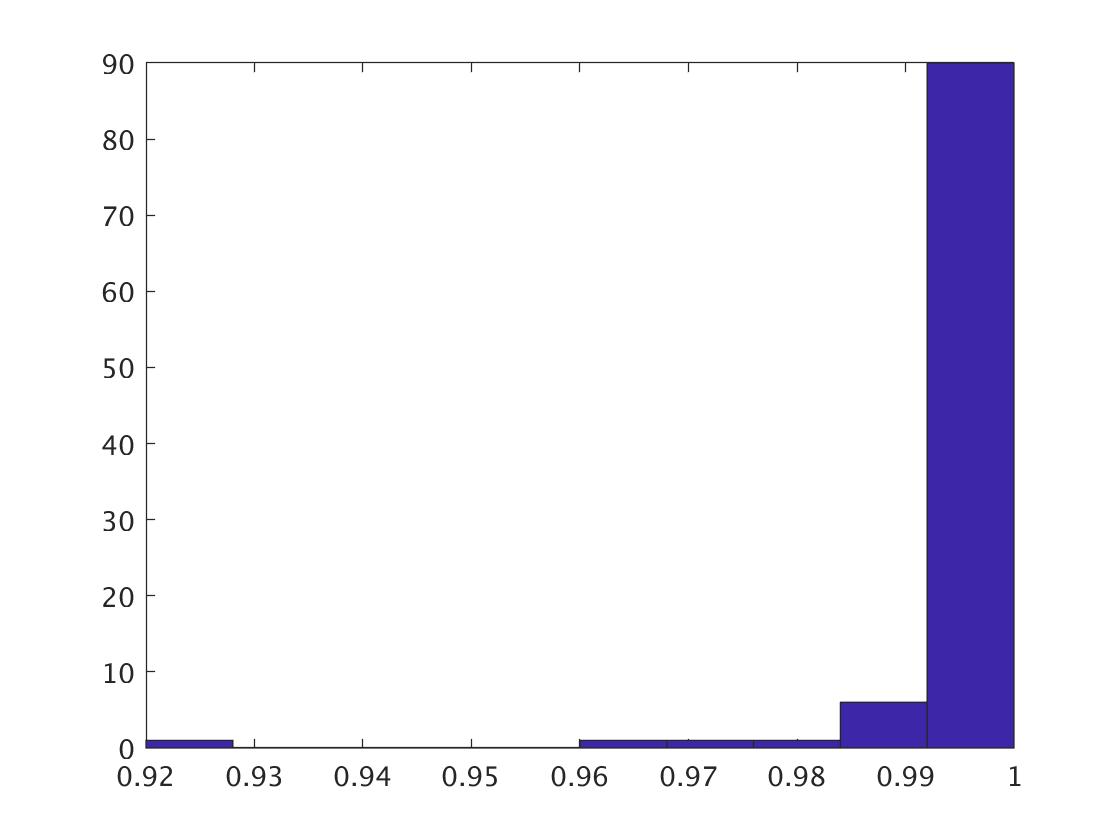}
    }
       \subfigure[$\varepsilon=0.002$, Second outlier]
    {
        \includegraphics[scale=.11]{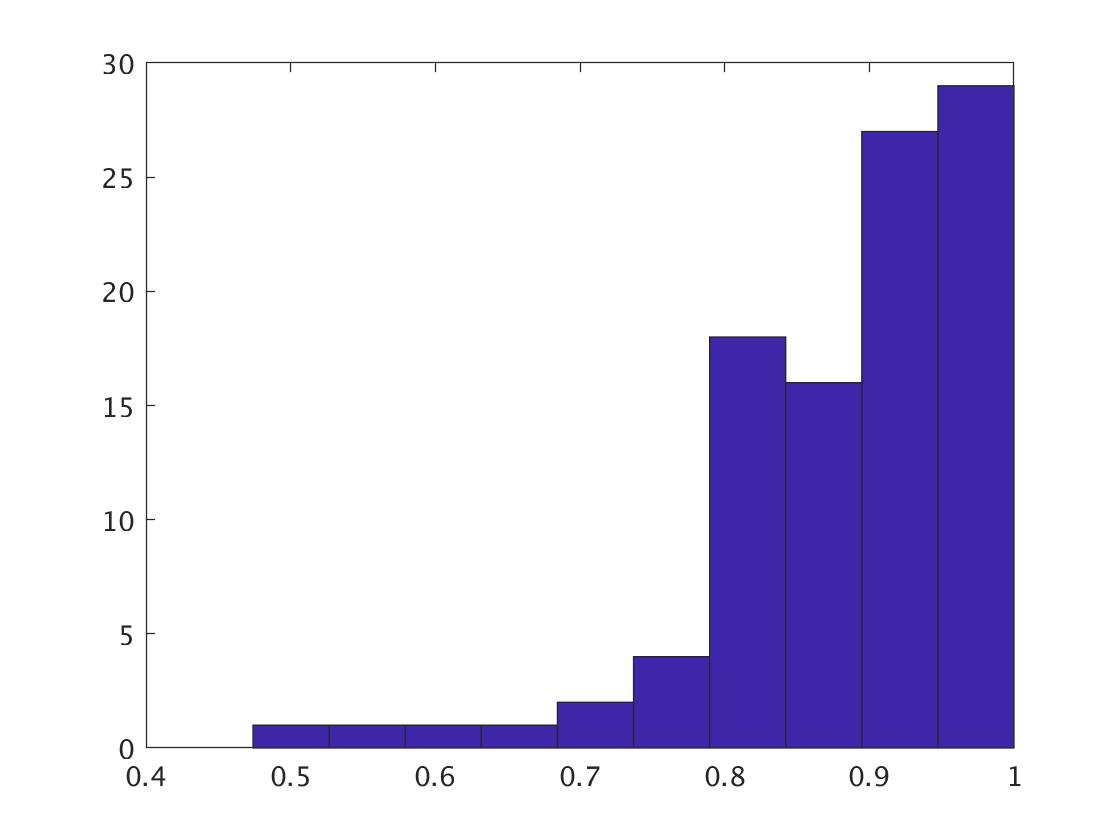}
    }
     \subfigure[$\varepsilon=0.002$,  Third outlier]
    {
        \includegraphics[scale=.11]{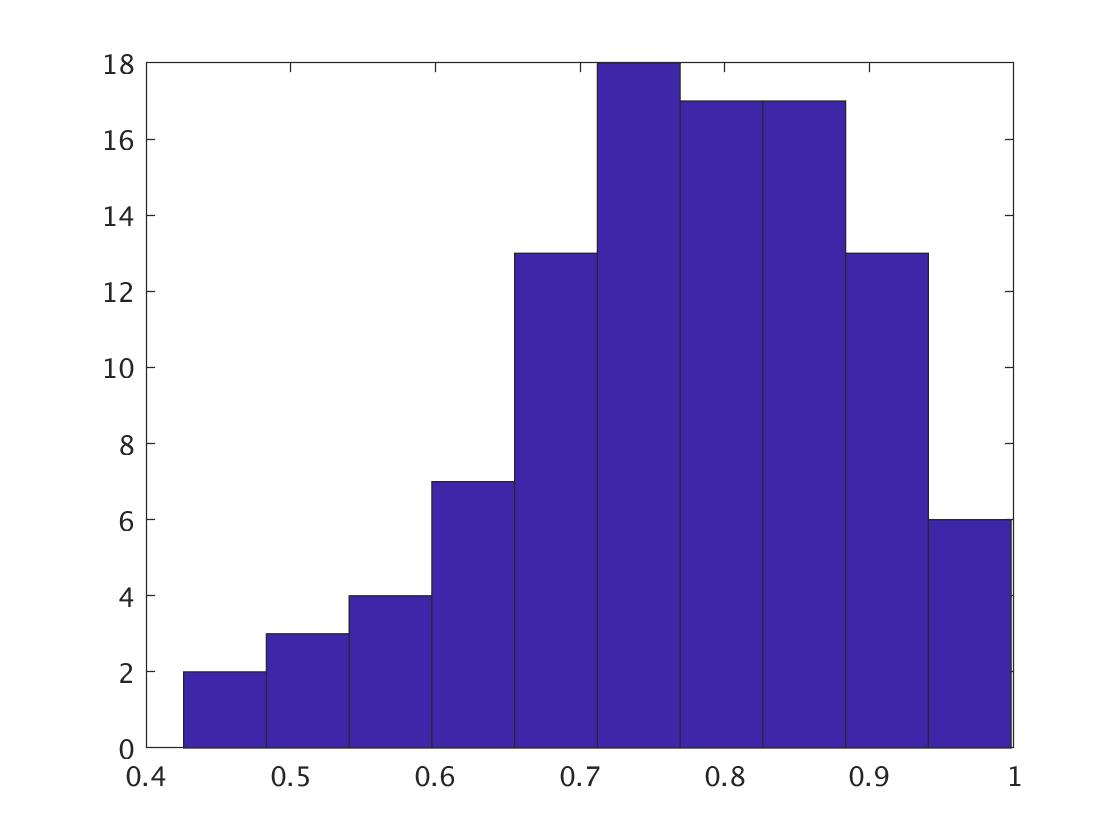}
    }    
       \subfigure[$\varepsilon=0.005$, First outlier] {
        \includegraphics[scale=.11]{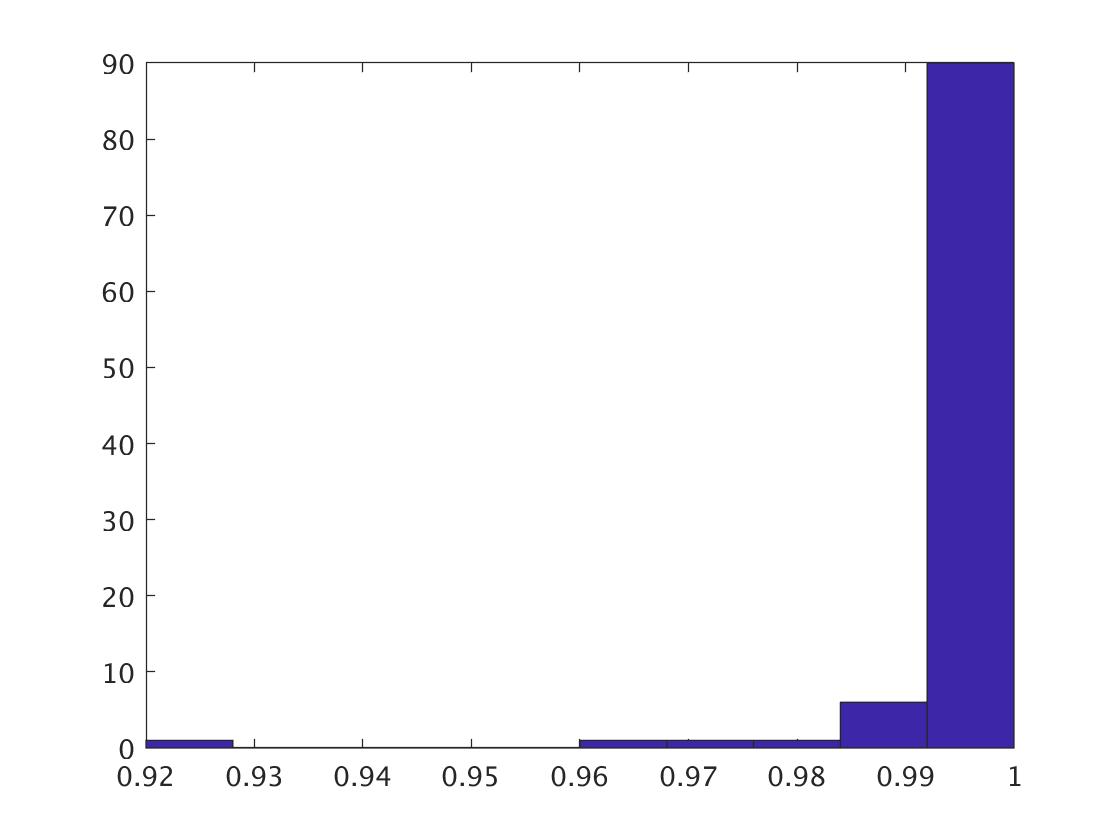}
    }
       \subfigure[$\varepsilon=0.005$, Second outlier]
    {
        \includegraphics[scale=.11]{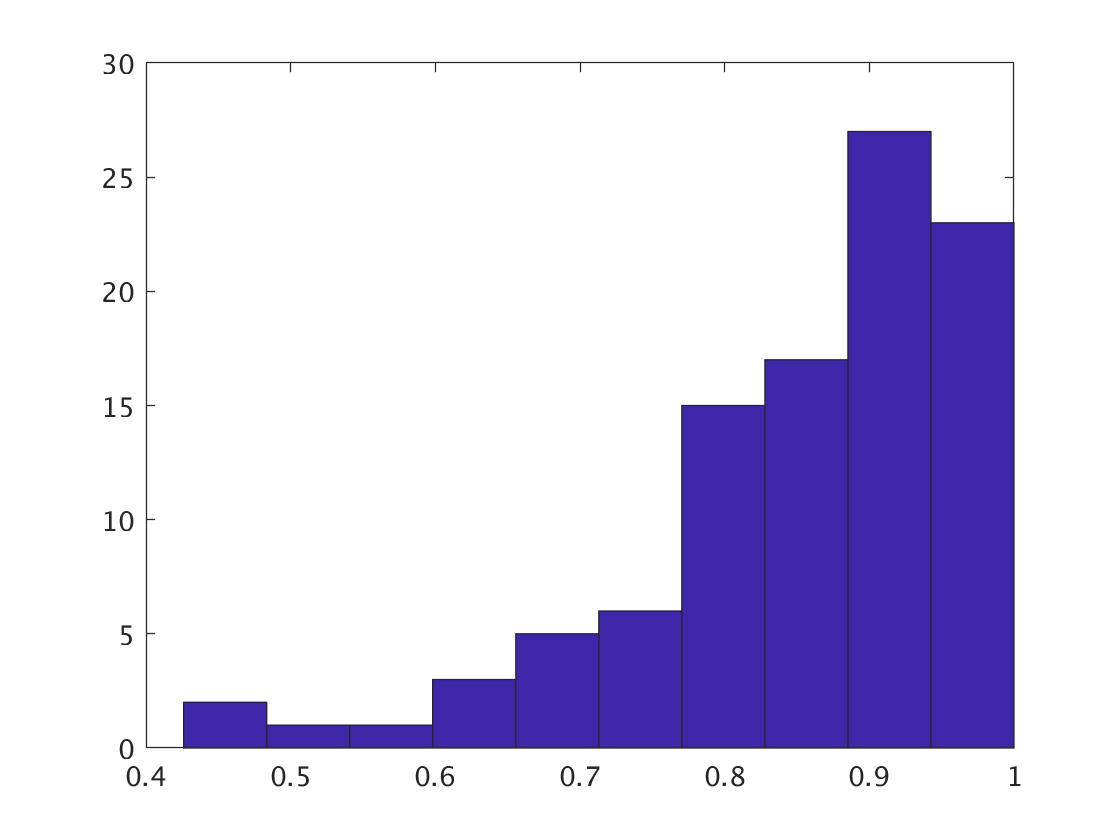}
    }
     \subfigure[$\varepsilon=0.005$, Third outlier]
    {
        \includegraphics[scale=.11]{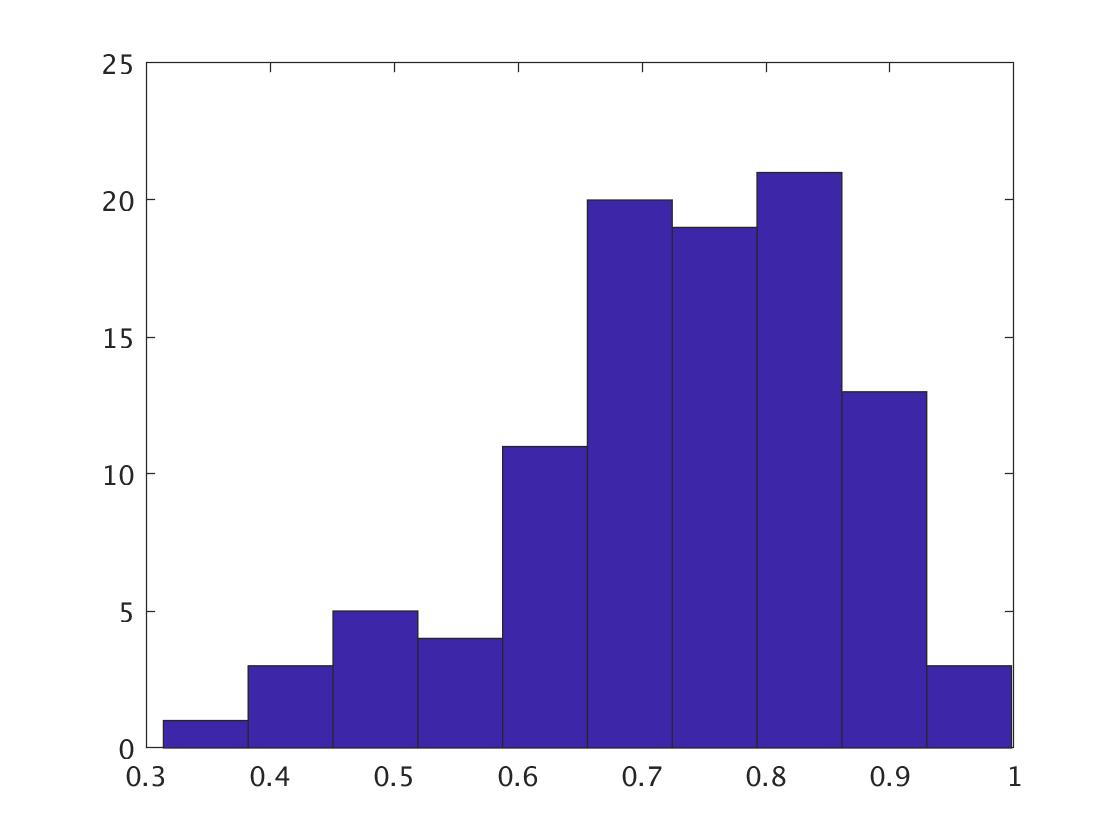}
    }     
\caption{\footnotesize
    Histograms for frequencies (1 minus the $p$-values) of  first, second, and third most significant points in the 100 persistence diagrams from the concentric circle example. See text for details.
    }
\label{fig:bagplotfrequencies}
\end{figure}

Histograms for each of the sets $\{f^*_k(j)\}_{k=1}^{100}$ are shown in  Panels (a)--(c) of Figure \ref{fig:bagplotfrequencies}. Alternatively, the first three rows of Table \ref{table:pvalues} show more succinct, and  easier to absorb, information. For example, looking at the first block of results (those with ``$\varepsilon=0$", a parameter that will be explained in a moment)  in the table, we see that,  in 93\% of the cases, the most extreme point in the original persistence diagram (which we know to be `real') was assigned a $p$-value of less that 0.01, and in 99\% of the cases a $p$-value of less than 0.05.  These are excellent results.

\begin{table}[h!]
{\footnotesize
\begin{center}
\begin{tabular}{|c|c|rrr|}
\hline
$\varepsilon$ & Outlier & & $p$-value &\\
 & number&  $\leq .01$ &   $\leq 0.05$ & $\leq 0.10$ \\  \hline 
0.000 & 1 & .93 &  .99  & 1.00 \\
 & 2 & .16&.39 &.72 \\
 & 3 & .03 & .11 & .38 \\   \hline  
 0.001 & 1 & .92 &  .99  & 1.00 \\
 & 2 &.14 &.28 & .58\\
 & 3 & .02& .06&.19 \\   \hline 
0.002 & 1 & .92 &  .99  & 1.00 \\
 & 2 & .12 & .27& .54 \\
 & 3 &.02 & .05& .16\\   \hline  
 0.005& 1 & .92 &  .99  & 1.00 \\
 & 2 & .11& .21& .47\\
 & 3 &.01 & .02 & .09\\ \hline
 \end{tabular}
\caption{\footnotesize{Frequency at which $p$-values were found for  the  first, second, and third most significant points in the 100 persistence diagrams from the concentric circle example. See text for details.  }}
\label{table:pvalues}
\end{center}
}
\end{table}

On the other hand, the second and third most extreme points (which we know to be `noise') were assigned low $p$-values an uncomfortable number of times. For example, the second point was assigned $p$-values of less than 0.01 and 0.05 in 16\%  and  39\% of the cases, respectively. This is far too often to be reasonable, from a practical point of view.

Studying the actual boxplots one sees immediately from where come these relatively high $p$-values: from outliers of that are perhaps removed from the main bulk of points, but extremely close to the diagonal, as in the outlier of Panel (b) of Figure \ref{fig:bagplotexamples}. These are clearly points that should not be considered informative, and so we repeated the above analysis, but, this time, ignored all outliers that were a distance of $\varepsilon$ or less from the diagonal. Since there is no way to know what a `good value' of $\varepsilon$ should be, we took three different values, 0.001, 0.002 and 0.005. These represent, on average, approximately 0.25\%, 0.5\%, and 1.25\% of the distance from the diagonal to the furthest point in the diagram. That is, all are      so close to the diagonal that no statistical test is needed to be comfortable with a decision to ignore them.

However, it turns out that the specific value of $\varepsilon$ chosen does not matter very much, as long as it is non-zero. Panels (d)--(l) of Figure 
\ref{fig:bagplotfrequencies} contain histograms corresponding to (a)--(c), arising from a repetition of the experiment described above, but with outliers close to the diagonal removed, for each of the three non-zero values of $\varepsilon$. Similarly, the last three blocks of Table \ref{table:pvalues} summarise the decisions one would reach based on these results.  The result is that the rate at which the true outlier is identified is essentially unchanged,  the rate at which `false' outliers are identified as `true topological signal' is significantly  reduced, for both the second and third outlier, all $p$-values, and all values of $\varepsilon$. Perhaps most importantly, from the point of view of practice, is that the results are basically insensitive to the specific value of $\varepsilon$.

In summary,  despite the fact that the bagplot based  RST procedure here is rather involved, with quite a lot of wishful thinking along the way, in practice it seems to work remarkably well.

\appendix
\section{Appendices}

\subsection{Simulating from $\bar f^G$}
\label{sec:appendixglobalshape}
While the MCMC calculations of the paper are well summarised in Algorithm 1 of Section \ref{sec:replications}, sampling from the distribution $\bar f^G(x)$  of \eqref{eq:kernel-new} itself requires some care. While reasonably straightforward, there are  numerical subtleties, and so for completeness we describe the full procedure here.

Recall that $\bar f^G(x)$ is just the kernel density estimate  $\hat f^G$, restricted to $\real\times\real_+$, and normalised. To sample from it, we first  denote by $R$ the  smallest  rectangular subset of the half plane which includes the 
set $\{x\in\real\times\real_+: \bar f^G(x)>\varepsilon\}$, for some $\varepsilon>0$ that is  case specific.
Then divide $R$ into $I_1\times  I_2$ equal sized rectangles $I_{ij}$, where $I_1$ and $I_2$ are typically of the order of 100, but, again, case specific.

The second step involves assigning probabilities to these rectangles, which, a prior, could be done by integrating $\bar f^G$ over each one. However, noting the original empirical density $\hat f^G$ comes from a  Gaussian kernel, considerable computational time is saved by first defining its integrated version
\beqq
\hat F^G (x) \ \definedas \ \frac{1}{n} \sum _{i=1}^{n}\Phi_{\Sigma} (x-x_{i}),
\eeqq
where $\Phi_\Sigma$ is the Gaussian (cumulative) distribution function corresponding to the Gaussian kernel in the definition of $\hat f^G$. Extend $\hat F^G$ to a measure on rectangles in the usual way, and define the probabilities (which now sum to 1) 
\beqq
p_{ij} \ = \ \frac{\hat F^G(I_{ij})}{\sum_{i=1}^{I_1}\sum_{j=1}^{I_2}\hat F^G(I_{ij})}.
\eeqq
By taking any linear enumeration of the indices $(i,j)$ it is now trivial to chose a rectangle at random, according to these probabilities, by the inverse transform method (e.g.\  \cite{RobertCasella,Handbook}). 

Having chosen a rectangle, we now chose a point uniformly, at random, from it. This is the value $x^*$ taken for Step 3 of Algorithm 1.

\subsection{The model of \cite{PNAS}}
\label{sec:appendixearliermodel}

The model originally developed in \cite{PNAS}, as with the one used in  the current paper, was a Gibbs distribution, and so can be described through its Hamiltonian, as below, retaining the notation of Section \ref{sec:modelling-new}. We shall do this only for projected persistence diagrams, so that each point $x$ in the diagram is of the form $x=(x^{(1)},x^{(2)})\in\real\times\real_+$.

Define 
 \beqq
\sigma_H^2 =  \sum_{x\in\tx_N} \big(x^{(1)} - \bar x^{(1)}\big)^2,\ \ \
\sigma_V^2 = \sum_{x\in\tx_N} \big(x^{(2)} \big)^2,
\eeqq
where $\bar x^{(1)}=N^{-1} \sum_{i=1}^N x_i^{(1)}$, so  
$\sigma_H^2$ is  the variance of the horizontal points. On the other hand,  $\sigma_V^2$ is square of the $L_2$ norm  of the vertical points, rather than the centred variance  (because of the non-negativeness of  the $x^{(2)}$). 

For integral $K> 0$, a collection $\Theta = (\theta_H,\theta_V,\theta_1,\dots\theta_K)$ of $\real$-valued parameters, and a $\delta>0$, define the Hamiltonian
\beq
H_{\delta,\Theta}^K(\tx_N)
= \theta_H \sigma^2_H
+\theta_V  \sigma^2_V
+ \delta^{-2} \sum_{k=1}^K \theta_k 
\sum_{i=1}^N \sum_{z\in\cN_k(x_i)} \|z-x_i\| \, \mathbbm{1}_{\{\|z-x \|\leq \delta \}} .
\label{eq:HamiltonianOld}
\eeq
With this Hamiltonian replacing the one defined by \eqref{eq:TheGibbs} the remaining development in \cite{PNAS} -- in particular that of an appropriate pseudolikelihood model -- is parallel to that in Section \ref{sec:modelling-new}.

We note though the main differences between the models. The first is the parameter $\delta$, which limits nearest neighbour interactions only to those neighbours that are closer than $\delta$. This had a mild numerically stabilising effect in the model defined by \eqref{eq:HamiltonianOld}, that, for reasons that are not entirely clear, disappeared in the model of the current paper. Consequently, we no longer use it. The first two terms in the Hamiltonian, involving second moments, were intended to play  the role that the empirical density $\bar f^G$ plays in the current paper; viz.\ they controlled the overall shape of the random diagrams, and worked ``against" the control resulting from   the nearest neighbour interactions. However, as shown by most of the examples in Section \ref{sec:examples-new}  -- in particular the Gaussian excursion set and non-concentric circles examples -- these terms were not able to capture many of the subtleties found in persistence diagrams. Furthermore, as the MCMC simulations progressed, the simulated diagrams had a tendency to move towards the diagonal, in a fashion that was inconsistent with their overall use.

\bibliographystyle{abbrvnat}


\bibliography{PDModelling}

\end{document}